\documentclass[11pt]{article}

\newtheorem{theorem}{Theorem}[section]
\newtheorem{lemma}[theorem]{Lemma}
\newtheorem{proposition}[theorem]{Proposition}
\newtheorem{corollary}[theorem]{Corollary}

\newtheorem{Rk}[theorem]{Remark}

\usepackage{fourier}
\usepackage{enumerate}
\usepackage{graphicx,amsmath,amsfonts,hyperref}
\usepackage{authblk}
\usepackage{geometry}
\newcommand{\NN}{{\mathcal N}}
\newcommand{\TT}{{\mathcal T}}
\newcommand{\FF}{{\mathcal F}}
\newcommand{\E}{{\mathbb E}}
\newcommand{\N}{{\mathbb N}}
\newcommand{\Z}{{\mathbb Z}}
\newcommand{\R}{{\mathbb R}}
\newcommand{\ind}[1]{\mathbf{1}_{#1}\,}
\newcommand{\PPP}[1]{{\mathbb{P}}\left(\,#1\,\right)}
\newcommand{\EEE}[1]{{\mathbb{E}}\left[\,#1\,\right]}
\newcommand{\A}[1]{{A\big(n\mathbf{1}_{I(#1)}\big)}}
\newcommand{\PP}{{\mathbb{P}}}
\newcommand{\EE}{{\mathbb{E}}}
\newenvironment{prooft}[1]{\vskip 2mm\noindent {\it{ Proof of  #1.}}} {\hfill  $\square$ \vskip 2mm \noindent}
\newcommand{\conv}[2][n]{\underset{#1\rightarrow #2}{\longrightarrow}}

\title{The bi-dimensional Directed IDLA forest}

\author[,1]{\textsc{Nicolas Chenavier} \thanks{\texttt{nicolas.chenavier@univ-littoral.fr}}}
\author[,2]{\textsc{David Coupier} \thanks{\texttt{david.coupier@imt-lille-douai.fr}}}
\author[,3]{\textsc{Arnaud Rousselle} \thanks{\texttt{arnaud.rousselle@u-bourgogne.fr}.}}
 \affil[1]{Universit\'e du Littoral C\^ote d'Opale, UR 2597, LMPA, Laboratoire de Math\'ematiques Pures et Appliqu\'ees Joseph Liouville,
62100 Calais, France.}
\affil[2]{Institut Mines T\'el\'ecom Lille Douai, Cit\'e Scientifique, 59655 Villeneuve d'Ascq, France.}
\affil[3]{Institut de Math\'ematiques de Bourgogne, UMR 5584, CNRS, Universit\'e Bourgogne Franche-Comt\'e,
 F-21000 Dijon, France.}
\date{}
\begin{document}
\maketitle

\begin{abstract}
We investigate three types of Internal Diffusion Limited Aggregation (IDLA) models. These models are based on simple random walks on $\Z^2$ with infinitely many sources that are the points of the vertical axis $I(\infty)=\{0\}\times\Z$. Various properties are provided, such as stationarity, mixing, stabilization and shape theorems. Our results allow us to define a new directed (w.r.t.\,the horizontal direction) random forest spanning $\Z^2$, based on an IDLA protocol, which is invariant in distribution w.r.t.\,vertical translations.
\end{abstract}




\section{Introduction}

The Internal Diffusion Limited Aggregation (IDLA) is a random growth model first introduced for chemical applications in 1986 by Meakin and Deutch \cite{MD86} and then, in a mathematical framework, by Diaconis and Fulton in \cite{DF91}. In this model, the aggregate is recursively defined by adding to the aggregate the first site out of the current aggregate visited by a random walk starting from some source point. The classical IDLA model is constructed in $\Z^d$ as follows. We start with $A_0=\emptyset$. At step $n$, a simple symmetric random walk starts from the origin $0$ until it exits the current aggregate $A_{n-1}$, say at some vertex $z$, which is added to $A_{n-1}$ to get  $A_n = A_{n-1}\cup\{z\}$. In the classical IDLA model (and also in this paper), the word \textit{particle} is used to refer to the random walk which is stopped when it exits the current aggregate $A_{n-1}$, and settled on the new vertex $z$.

A first shape theorem was established by Lawler, Bramson and Griffeath in \cite{LBG92} for the classical IDLA model. It asserts that the aggregate $A_n$ (when it is suitably normalized) converges a.s.\,to an Euclidean ball as $n$ goes to infinity, with fluctuations (w.r.t.\,the limit shape) which are at most linear. Since then, several papers (by Lawler \cite{Lawler95}, Asselah and Gaudilli\`ere \cite{AG13long,AG13short,AG14} and Jerison, Levine and Sheffield \cite{JLS12,JLS13,JLS14}) have improved the bounds for fluctuations which are known to be logarithmic in $2D$ and sublogarithmic in higher dimensions.

Recently, many variants of this problem have been considered. In particular, IDLA on discrete groups with polynomial or exponential growth have been studied in \cite{B04,BlB07}, on non-amenable graphs in \cite{H08}, with multiple sources in \cite{LP10}, on supercritical percolation clusters in \cite{DLYY, Shellef}, on comb lattices in \cite{AR16,HS12}, on cylinder graphs in \cite{JLS14b,LS,S19}, constructed with drifted random walks in \cite{L14} or with uniform starting points in \cite{BDCKL}.

\begin{center}
\begin{figure}
\begin{center}
\includegraphics[height=6cm]{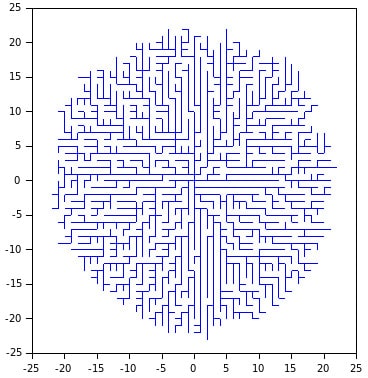}
\end{center}
\caption{{A realization of $\mathcal{T}_{1500}$. \label{fig:tree}} }
\end{figure}
\end{center}

A random infinite tree $\TT_\infty$ can be associated with the sequence of IDLA models $(A_n)_{n\geq 0}$ defined above in a very natural way. To our knowledge, this object has not been introduced in the literature. The tree $\TT_1$ only consists of the root $0$. By induction, $\TT_n$ is obtained by adding to $\TT_{n-1}$ the new vertex $z$ such that $A_n=A_{n-1}\cup\{z\}$ and the edge used by the $n$-th particle to reach $z$ from $A_{n-1}$ (see Figure~\ref{fig:tree}). Hence, we can define a.s.\,a random graph
\[
\TT_\infty = \bigcup_{n\geq 1} \! \uparrow \TT_n,
\]
which actually is a tree (since each vertex of $\Z^2$ may only be added once) rooted at the origin. The lower bound for the shape theorem specifies that its edge set spans the whole set $\Z^2$.

The first question about the random tree $\TT_\infty$ concerns the existence of (many) infinite branches with asymptotic directions {(it is clear that $\TT_\infty$ contains at least one infinite branch since the degree of each of its vertices  is a.s.\,bounded)}. In \cite{HN01}, Howard and Newman have developed an efficient strategy leading to such results. The key point would be to prove that the fluctuations of the branch in $\TT_{\infty}$ joining the origin to a given vertex $z\in\Z^2$ w.r.t.\,the segment $[0,z]$ (in $\R^2$) are negligible w.r.t.\,$| z|_2$, where $|\cdot|_2$ denotes the Euclidean norm in $\R^2$. In this case, the tree $\TT_\infty$ is said \textit{straight}. Although one might strongly conjecture such a result (especially because fluctuations in the shape theorem are logarithmic), it is difficult to prove it for several reasons. First, any branch $\gamma$ of the IDLA tree $\TT_{\infty}$ is not produced by a single particle but by many particles, each of them adding exactly one edge depending on the shape of the current aggregate. Secondly, the radial character of $\TT_\infty$ (its branches are directed to the origin) prevents its distribution to satisfy {many} useful invariance properties. 

A way to overcome this (second) obstacle would be to consider a directed forest  w.r.t.\,to some vector $u\in\R^2$ (i) whose distribution admits some invariant translation (w.r.t.\,any orthogonal vector to $u$) properties making its study easier than the one of $\TT_\infty$ and (ii) which approximates locally and far from the origin, {\it i.e.}\,on the ball $B(-nu,R)$ with $n\gg 1$ and $R$ constant, the distribution of the random tree $\TT_\infty$ so that we could transfer results about this directed forest to $\TT_\infty$. This strategy relies on the following remark: the radial character of the branches of $\TT_\infty$ restricted to $B(-nu,R)$ should fade away as $n\to\infty$ (with $R$ constant) so that the branches should be directed w.r.t.\,the vector $u$ and not toward the origin. In particular, this strategy has been used successfully by Baccelli and Bordenave in \cite{BB07} to approximate the \textit{Radial Spanning Tree} by the \textit{Directed spanning Forest} and to show that it is straight. See also \cite{C18} in which the author exploits this link between trees and the associated forests to quantify the density of infinite branches of these trees. Furthermore, the associated forest generally is much more than a tool: in \cite{CSST}, for example, the authors prove that the Directed spanning Forest mentioned above converges in distribution to the Brownian web.

This strategy is the original motivation for this work.\\

One of our main results is the construction of a new random forest $\FF_\infty$, called the \textit{directed infinite-volume IDLA forest}, and directed w.r.t.\,the horizontal vectors $u=(\pm 1,0)$ (see Section~\ref{subsec:forestidla}). This construction relies on the use of three IDLA processes (and their properties) whose limits are denoted by $A_n[\infty]$, $A_n^{\ast}[\infty]$ and $A_n^\dag[\infty]$, and based on  infinitely many sources given by the sites of the vertical axis $I(\infty)=\{0\}\times \Z$. Thus we show that the distribution of $\FF_\infty$ is invariant (and even mixing) w.r.t.\,vertical translations (Theorem~\ref{th:forestIDLA}). Finally, we think that the directed IDLA forest $\FF_\infty$ is an interesting mathematical object with several conjectures (finiteness of its trees, straightness of its branches, scaling limit). In particular,  we conjecture that $\FF_\infty$ approximates in distribution, locally and far from the origin, the IDLA tree $\TT_\infty$ (see Conjecture 1 in Section~\ref{sect:conjectures}).

\begin{center}
\begin{figure}
\begin{center}
\includegraphics[height=6cm]{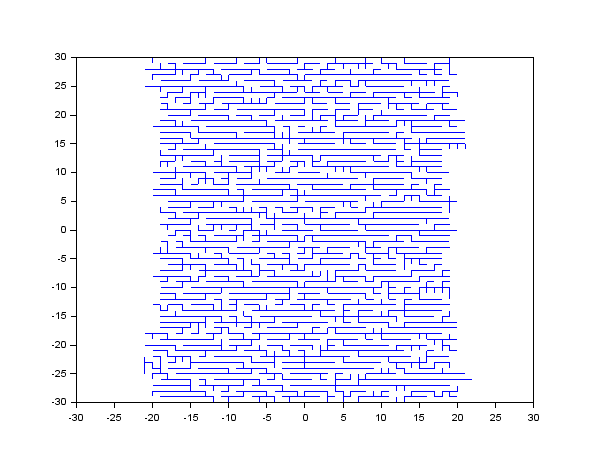}\end{center}
\caption{{A realization of ${\FF_{40}^\dag[200]}$ observed in $\Z_{30}$.
 \label{fig:forest}} }
\end{figure}
\end{center}

Let us note that Berger, Kagan and Procaccia proposed in \cite{BKP19} a random forest model based on an IDLA protocol. Nevertheless, their construction is based on oriented random walks which certainly simplifies the construction but also prohibits U-turns that are possible in $\TT_\infty$ and also in our directed IDLA forest $\FF_\infty$ (see e.g. the edge with vertex $(8,-8)$ in Figure \ref{fig:forest}). Hence, their model does not seem to be a good candidate to mimic the infinite IDLA tree $\TT_\infty$ and to capture its properties.\\

{More recently, in the context of External DLA process, an infinite stationary DLA on the upper half planar lattice growing from an infinite line has been introduced on \cite{PYZ20} (see also \cite{MPZ19} for a DLA model on a long line segment). The approach developped in \cite{PYZ20} provides an interesting track to define a stationary random forest in the Internal DLA context. However the strategy developped in \cite{PYZ20} requires a big preparative work on a stationary harmonic measure (\cite{PYZ21,PZ19}). To overcome this issue, we will strongly use instead the so-called Abelian property which holds for the Internal protocol but not for the external one.}

Let us describe the construction of the three infinite aggregates, namely $A_n[\infty]$, $A_n^{\ast}[\infty]$ and $A_n^\dag[\infty]$, { as random subsets of $\mathbb{Z}^2$}. We say that a particle is \textit{emitted} or \textit{sent from level $i$} if the underlying random walk starts from the source $(0,i)$. First, we consider finite aggregates $A_n[M]$, with $M\geq 0$, in which $n$ particles are sent from each site of $I(M)=\{0\}\times\llbracket -M,M\rrbracket,$ where $\llbracket a,b\rrbracket=[a,b]\cap\Z$ for any $a\leq b$. The first $n$ ones are sent from level $0$, then the following $n$ particles from level $1$, next  from level $-1$ and so on up to last particles sent from levels $M$ thus $-M$. In the rest of the paper, this specific order for sending particles is referred to as the {\it usual order}. The sequence $(A_n[M])_{M\geq 0}$ is increasing and allows us to define a limiting infinite aggregate $A_n[\infty]$. Note that this model is not covered by \cite{LP10}.

The second infinite IDLA {model}, denoted by $A_n^{\ast}[\infty]$, is defined in the same spirit as $A_n[\infty]$ but this time by launching a random number $N_i$ of particles from each level $i$ w.r.t.\,the usual order. The $N_i$'s are i.i.d.\,Poisson random variables with parameter $n$ and are independent of the underlying random walks. 

Although $A_n[\infty]$ and $A_n^{\ast}[\infty]$ are not identically distributed (see Proposition~\ref{prop:DifferentLaw}), they have in common the usual order in which particles are sent and for this reason they are not conducive to define a translation invariant IDLA forest. However, this usual order will be  particularly useful to derive the stabilization properties for $A_n[\infty]$ and $A_n^{\ast}[\infty]$. A natural way to get back an invariance property in distribution w.r.t.\,vertical translations consists in sending particles \textit{uniformly} from the vertical axis $I(\infty)$. This idea motivates the definition of the third infinite IDLA {model} $A_n^\dag[\infty]$. So let us consider a family $(\NN_{i})_{i\in\Z}$ of i.i.d.\,Poisson point processes (PPP's) in $\R^{+}$ with intensity $1$ (and also independent of the underlying random walks). Each PPP $\NN_{i}$ is attached to the level $i$ in the sense that particles are sent from the source $(0,i)$ according to the clocks given by $\NN_{i}$. Let us specify that the trajectories of particles are instantaneously realized. Hence, for any $M$, the aggregate $A_n^\dag[M]$ is built by sending particles from the source set $I(M)$ according to the clocks given by the corresponding PPP's up to time $n$. Remark that this construction ensures that, at each time, the next particle (if it exists) is sent from a source chosen uniformly on $I(M)$. Thus, $A_n^\dag[\infty]$ is defined as the increasing union of the $A_n^\dag[M]$'s.

Defining an IDLA forest {$\FF_n^\dag[M]$} from the (finite) aggregate $A_n^\dag[M]$ is easy to do (see Section~\ref{sect:ForetnM} and Figure~\ref{fig:forest}). But taking the limit $M\to\infty$ in the sequence $({\FF_n^\dag[M]})_{M\geq 0}$ needs to take some precautions. Indeed, given $M'>M$, any particle starting from a level $M'\geq |i|>M$ may generate a set of discrepancies between both forests ${\FF_n^\dag[M]}$ and ${\FF_n^\dag[M']}$ through a tricky phenomenon that we have called a \textit{chain of changes} and described in Section~\ref{sect:ChainChanges}. The existence of arbitrarily long chains of changes in the sequence $({\FF_n^\dag[M]})_{M\geq 0}$ is the main obstacle to define the directed IDLA forest $\FF_\infty$. {We stress that one of the difficulties in this paper is the loss of the so-called Abelian property. Such a property claims that the distribution of any aggregate based on an IDLA protocol does not depend on the order in which the particles are sent {(see Section \ref{sec:abelian} for a precise statement)}. Unfortunately, this property does not hold for the forests. }\\ 

To overcome {these} obstacle{s}, we proceed as follows. We first establish two stabilization results for both infinite aggregates $A_n[\infty]$ and $A_n^{\ast}[\infty]$. The first one (Theorem~\ref{th:stabilize})  asserts that $A_n[\infty]$ and $A_n^{\ast}[\infty]$, restricted to a neighborhood of the origin, are not sensitive to particles coming from far levels. On the opposite, the second one (Theorem~\ref{cor:StabWRTO}) claims that $A_n[\infty]$ and $A_n^{\ast}[\infty]$, restricted to high levels, are not sensitive to particles sent from a neighborhood of the origin. These fruitful tools imply that the infinite aggregates $A_n[\infty]$ and $A_n^{\ast}[\infty]$ are mixing w.r.t.\,vertical translations. Combining with Proposition~\ref{prop:DifferentLaw}, we obtain that $A_n^{\ast}[\infty]$ a.s.\,avoids an infinite number of horizontal lines which means that $A_n^{\ast}[\infty]$ is made up with infinitely many finite and disjoint connected components. The same holds for $A_n^\dag[\infty]$ since $A_n^\dag[\infty]$ and $A_n^{\ast}[\infty]$ are equally distributed thanks to the Abelian property. This is this latter statement which prevents the existence of arbitrarily long chains of changes in the sequence $({\FF_n^\dag[M]})_{M\geq 0}$ and allows us to take first the (vertical) limit $M\to\infty$ in the sequence $({\FF_n^\dag[M]})_{M\geq 0}$ and thus the (horizontal) limit $n\to\infty$ to finally define the directed infinite-volume IDLA forest $\FF_\infty$. 

Finally, we prove {shape theorems (Section \ref{sec:shapeTh})} for the three aggregates $A_n[\infty]$, $A_n^{\ast}[\infty]$ and $A_n^\dag[\infty]$ restricted to the strip $\Z_{n^{\alpha}} = \Z\times\llbracket -n^{\alpha} , n^{\alpha}\rrbracket$ (for any $\alpha>0$) as $n$ tends to infinity. Adapting to our context the strategy developed by Asselah and Gaudilli\`ere \cite{AG13long,AG13short}, we prove that, with probability $1$ and for any $n$ large enough, the fluctuations (w.r.t.\,the Hausdorff distance) between $A_n[\infty] \cap \Z_{n^{\alpha}}$ and the rectangle $\left\llbracket-\frac{n}{2},\frac{n}{2}\right\rrbracket\times\left\llbracket - n^{\alpha},n^{\alpha}\right\rrbracket$ are at most logarithmic {(Theorem \ref{th:ShapeTh}). A shape theorem also holds for $A_n^{\ast}[\infty]$ and $A_n^\dag[\infty]$ (Theorem \ref{th:ShapeThPoisson})}. As a consequence, the vertex set of the directed IDLA forest $\FF_\infty$ fulfills the whole set $\Z^2$.\\ 

The rest of this paper is organized as follows. In Section~\ref{sec:def}, we define the infinite aggregates $A_n[\infty]$, $A^*_n[\infty]$ and $A^\dag_n[\infty]$ and we state their first properties. {The proofs of the properties concerning $A_n[\infty]$ and $A^*_n[\infty]$  mostly rely on the particular order in which the particles are sent.} Then, in Section~\ref{sec:hightpart}, we prove that a.s.\,particles emitted above some random level contribute to the aggregate before visiting the strip $\Z_M$ and that $A_n[\infty]$ and $A_n^*[\infty]$ are not equally distributed. In Section~\ref{sec:stabOrigine}, we show that the aggregates above some random levels a.s.\,do not depend on  particles which are sent arround the origin. The mixing properties of the aggregates are discussed in Section~\ref{sec:mixing} and used to deduce that $A_n^*[\infty]$ and $A_n^\dag[\infty]$ a.s.\,avoid infinitely many lines ${\Z}\times\{i\}$, $i\in \Z$.  Section~\ref{sec:shapeTh} is devoted to the shape theorems. Finally, Section~\ref{sec:forest} contains the construction of the infinite volume IDLA forest that spans $\Z^2$ and whose distribution is invariant w.r.t.\,vertical translations. Our paper ends with three conjectures on the random forest $\FF_\infty$ and on the random tree $\mathcal{T}_\infty$. 

\section{Construction and properties of $A_n[\infty]$, $A_n^*[\infty]$ and $A_n^\dag[\infty]$\label{sec:def}}   

Let $n\geq 1$ be fixed. In this section, we introduce three infinite random aggregates with sources on $I(\infty)=\{0\}\times \Z$. Then we give some properties of these aggregates. {Before introducing them, we state  the Abelian property since it plays a crucial role.}

\subsection{Abelian property}
\label{sec:abelian}
{Let $A\subset \Z^2$ and $z\in \Z^2$. First, we recall the definition of the Diaconis-Fulton smash sum $A \oplus \{z\}$: 
\begin{itemize}
\item if $z\not\in A$, then $A \oplus \{z\} = A\cup\{z\}$;
\item if $z\in A$, then $A \oplus \{z\}$ is the random set obtained by adding to $A$ the vertex in which a simple random walk exits $A$. 
\end{itemize} 
The \textit{Abelian property} (see \cite{DF91}, p.\,97)   claims the following: given any finite set $\{z_1,\ldots, z_k\}$ of vertices of $\Z^2$, and a set $A\subset \Z^2$, the distribution of 
\[((A\oplus \{z_1\}) \oplus \{z_2\})\oplus \cdots \oplus \{z_k\}\]
does not depend on the order of the $z_i$'s, {\it i.e.} if $\sigma:\{1,\ldots, k\}\rightarrow \{1,\ldots, k\}$ is a permutation of $\{1,\ldots, k\}$, then we have the equality in distribution:
\[((A\oplus \{z_1\}) \oplus \{z_2\})\oplus \cdots \oplus \{z_k\} \overset{\text{law}}{=}  ((A\oplus \{z_{\sigma(1)}\}) \oplus \{z_{\sigma(2)}\})\oplus \cdots \oplus \{z_{\sigma(k)}\}.\]
}

\subsection{Construction of the aggregates}

\subsubsection{Construction  of $A_n[\infty]$}
\begin{center}
\begin{figure}
\begin{center}
\includegraphics[height=6cm]{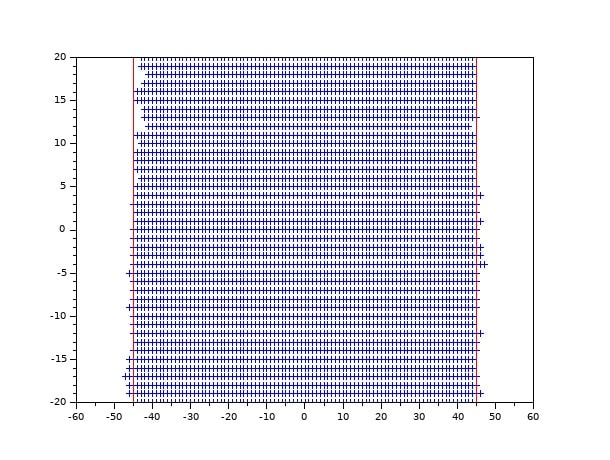}
\end{center}
\caption{{A realization of the aggregate $A_{90}[200]\cap \Z_{20}$ based on 90 particles per site $(0,i)$, with $|i|\leq 200$, and intersected by the strip $\Z_{20}$. \label{fig:shapetheorem}} }
\end{figure}
\end{center}

The aggregate $A_n[\infty]$ is a natural extension of the standard IDLA cluster  (see {\it e.g.} \cite{LBG92}) with sources on the $I(\infty)$ instead of the origin. To construct it, we first introduce a family of finite random aggregates $A_n[M]$, $M\geq 0$. When $M=0$, the random set $A_n[0]$ is the standard IDLA cluster with volume $n$. It is obtained inductively as follows. One by one, particles perform independent simple symmetric $2$-dimensional random walks. Each particle starts from the origin and moves until it reaches a site that has not been visited previously, at which point it stops. Then $A_n[0]$ is the cluster of occupied sites after the $n$-th particle stops. In a similar way, given a realization of $A_n[M-1]$, we send $n$ particles from the site $(0,M)$ then $n$ particles from the site {$(0,-M)$}. The set $A_n[M]$ denotes the aggregate which is produced by these $2n$ particles and by the aggregate $A_n[M-1]$. As an illustration, Figure~\ref{fig:shapetheorem} gives a realization of $A_{90}[200]$ when it is observed in the strip $\Z_{20}=\{0\}\times \llbracket -20, 20\rrbracket$.

By construction, the sequence of random aggregates $(A_n[M])_{M\geq 0}$ is increasing and the number of sites of $A_n[M]$ is 
\[\# A_n[M] = (2M+1)n.\]
The infinite random aggregate $A_n[\infty]$ is then defined as:
\begin{equation*}
\label{eq:defAn}
  A_n[\infty] = \bigcup_{M\geq 0}A_n[M].
  \end{equation*}

\subsubsection{Construction of $A_n^*[\infty]$}
\label{subsec:defpoisson}
The aggregate $A_n^*[\infty]$ is constructed in the same spirit as above but this time the number of particles which are sent from each site of $I(\infty)$ is no longer equal to $n$ but random. To define it, let $(N_{i})_{i\in\Z}$ be a family of independent Poisson random variables with parameter $n$. The random variable $N_{i}$ is the number of particles starting from $(0,i)$. Then, for each $j\in \Z$, we consider a family of simple random walks, starting from $(0,j)$, which are independent of $(N_i)_{i\in \Z}$. Given a realization of $(N_i)_{i\in \Z}$, we denote by $A_n^*[M]$, $M\geq 0$, the IDLA {model} obtained by sending particles from levels $|i|\leq M$ in the usual order, {\it i.e.}  first by sending the $N_0$ particles from $(0,0)$, thus the $N_1$ particles from $(0,1)$, thus the $N_{-1}$ particles from level $(0,-1)$ and so on. By construction, the cardinality of the set $A_n^*[M]$ is random and can be equal to zero with positive probability, so that $A_n^*[M]$ does not have the same distribution as $A_n[M]$. However, notice that the mean size of $A_n^*[M]$ is  \[\EEE{\# A_n^*[M]} = (2M+1)n.\]

Similarly to the above, the sequence of random aggregates $(A_n^*[M])_{M\geq 0}$ is increasing and the  infinite random aggregate $A_n^*[\infty]$ is defined as:
\[A_n^*[\infty] = \bigcup_{M\geq 0}A_n^*[M].\]
This new aggregate will be used to derive some properties of a third infinite random aggregate that we introduce below.

\subsubsection{Construction of $A_n^\dag[\infty]$}
As in the previous constructions, we first define an increasing family of finite random aggregates. The number of particles from each site is still Poisson but this time the order for which they are sent is modified: the particles are not sent w.r.t.\,the usual order but w.r.t.\,a family of random clocks. To do it, let $(\NN_{i})_{i\in\Z}$ be a family of independent and identically distributed Poisson point processes (PPP's) in $\R^{+}$, with intensity $1$. Each PPP $\NN_i$ provides an increasing sequence $\left(\tau_{i,j}\right)_{j\geq 1}$ of random clocks. Then, we attach to the collection $\left\{\tau_{i,j} : i\in\Z , j\geq 1\right\}$ a family of independent and identically distributed symmetric random walks $\left\{S_{i,j} : i\in\Z , j\geq 1\right\}$ which are also independent of the PPP's. In other words, at time $\tau_{i,j}$, the $j$-th particle from level $i$ starts and its trajectory, associated with  $S_{i,j}$, is instantaneously realized and adds a new site to the current aggregate. 

Now, let $M\geq 0$. We  denote by $A_{n}^\dag[M]$ the IDLA {model} obtained, with the same protocol as above, by sending particles from levels $|i|\leq M$ according to the PPP's $(\NN_{i})_{i\in\Z}$ until time $n$. As in Section~\ref{subsec:defpoisson}, the random number $N_i=\#\NN_{i}([0,n])$ of particles starting from level $i$ in $A_{n}^\dag[M]$ is Poisson with parameter $n$. Notice that, conditional on the r.v.'s $N_i=\#\NN_{i}([0,n])$, with $|i|\leq M$, only the order of the particles is changed between the aggregates $A_{n}^\dag[M]$ and $A^*_n[M]$. {As mentioned above, a remarkable property of (finite) IDLA {models} is the  Abelian property}. Such a property implies that, for each $M\geq 0$, 
\begin{equation} \label{eq:samemarginal} A_n^\dag[M]\overset{\text{law}}{=}A_n^*[M]. 
\end{equation} The above equality will be used to derive some properties of $A_n^\dag[M]$ from those which deal with $A_n^*[M]$. 

The aggregates $A_{n}^\dag[M]$, with $n,M\geq 0$, are based on the same family of PPP's  and on the same random walks. The following lemma shows that the sequences $(A_{n}^\dag[M])_{n\geq 0}$ (with $M$ fixed) and  $(A_{n}^\dag[M])_{M\geq 0}$ (with $n$ fixed) are increasing.

\begin{lemma}
\label{lem:inclusion}
For any integers $n,M\geq 0$, the following inclusions hold a.s. 
\begin{enumerate}[(i)]
\item $A_{n}^\dag[M] \subset A_{n}^\dag[M+1]$;
\item $A_{n}^\dag[M] \subset A_{n+1}^\dag[M] $.
\end{enumerate}
\end{lemma}

\begin{prooft}{Lemma~\ref{lem:inclusion}}
First, we prove (i). Let
\[
\kappa = \sum_{|i|\leq M+1} \#\NN_{i}([0,n])
\]
be the number of particles which are sent from levels $|i|\leq M+1$ until time $n$. We index them by $j=1,\ldots,\kappa$ according to their starting times $0<\tau_1<\ldots<\tau_{\kappa}<n$ (notice that they are a.s.\,all different). Recall that these particles are based on the same random walks for $A_{n}^\dag[M]$ and for $A_{n}^\dag[M+1]$. Some of these particles come from levels $\pm(M+1)$ and only concern the aggregate $A_{n}^\dag[M+1]$. For $j=1,\ldots,\kappa$, we denote  by $A[M,j]$ (resp. $A[M+1,j]$) the aggregate obtained until (or at) time $\tau_j$ with particles from levels $|i|\leq M$ (resp. from levels $|i|\leq M+1$). We set $A[M,0]=A[M+1,0]=\emptyset$. Let $1\leq j\leq\kappa$ and assume that a.s.\,$A[M,j-1]\subset A[M+1,j-1]$. If particle $j$ (which is sent at time $\tau_j$)  comes from level $\pm(M+1)$ then $A[M,j]=A[M,j-1]\subset A[M+1,j-1]\subset A[M+1,j]$. Otherwise, let $y$ be the site added to $A[M,j-1]$ by particle $j$. If $y\notin A[M+1,j-1]$ then
\[
A[M+1,j] = A[M+1,j-1] \cup \{y\} \supset A[M,j-1] \cup \{y\} = A[M,j] ~.
\]
Otherwise, the random walk associated with particle $j$ continues its trajectory till exiting $A[M+1,j-1]$ on a site $y'$. In this case,
\[
A[M+1,j] = A[M+1,j-1] \cup \{y'\} \supset A[M,j-1] \cup \{y\} = A[M,j] ~.
\]
By induction over $j=0,\ldots,\kappa$, we a.s.\,get $A_{n}^\dag[M+1]=A[M+1,\kappa]\supset A[M,\kappa]=A_{n}^\dag[M]$.

Assertion (ii) is easy to check by letting increase $A_{n}^\dag[M]$ with particles from levels $|i|\leq M$ on the time interval $(n,n+1]$.
\end{prooft}

Similarly to the previous subsections, we let
\[A_n^\dag[\infty] = \bigcup_{M\geq 0}A_n^\dag[M].\]
Because the sequences $(A_n^\dag[M])_{M\geq 0}$ and $(A_n^*[M])_{M\geq 0}$ are increasing, it follows from \eqref{eq:samemarginal} that 
\begin{equation} \label{eq:samedistributionl} A_n^\dag[\infty]\overset{\text{law}}{=}A_n^*[\infty]. 
\end{equation} 

{One of our motivations is to define a random forest which is stationary w.r.t.\,to vertical translations (see Section~\ref{sec:forest}). To get it, it is crucial to base our model on a family of PPP's like $A_n^\dag[\infty]$; a random forest which would be based only on a family of Poisson random variables and on the same protocol as $A^*_n[M]$ ({\it i.e.} by sending particles w.r.t.\,the usual order, first from level $0$, then from level $1$, then from level $-1$, \ldots) would be not stationary. However, this usual order makes it possible to build aggregates ``from the center''  and plays an important role in several proofs. Equalities \eqref{eq:samemarginal} and \eqref{eq:samedistributionl} then allow us to deduce results for the aggregates based on the PPP's from {those which concern the aggregates defined w.r.t.\,the usual order}.}

\subsection{First properties of the aggregates\label{subsec:firstprop}}

\subsubsection{Invariance w.r.t.\,symmetries and translations}
Let $k\in\Z$. In what follows, we denote by $\tau_{k}$ (resp.\,$S_{k/2}$) the translation operator w.r.t.\,vector $(0,k)$ (resp. {the symmetry} w.r.t.\,the horizontal axis $\R\times\{k/2\}$). The following proposition claims that the three infinite random aggregates are invariant w.r.t.\,translations and symmetries.

\begin{proposition}
\label{prop:transfo}
Let $n\geq 0$. The following properties hold:
\begin{enumerate}[(i)]
\item  the distributions of $A_n[\infty]$, $A^\ast_n[\infty]$ and $A^\dag_n[\infty]$ are invariant w.r.t.\,$\tau_k$, $k\in \Z$;
\item the distributions of $A_n[\infty]$, $A^\ast_n[\infty]$ and $A^\dag_n[\infty]$ are invariant w.r.t.\,$S_{k/2}$, $k\in \Z$.
\end{enumerate}
\end{proposition}
The above result is intuitively clear for $A^\dag_n[\infty]$ since the Poisson clocks are independent. However, this is not so intuitive for $A_n[\infty]$ and $A^\ast_n[\infty]$ because, in their constructions, the particles are sent in the specific ``usual'' order. 

Notice that $\tau_kA_n[\infty]$ (resp. $S_{k/2}A_n[\infty]$) is the increasing union of the aggregates $\tau_kA_n[M]$ (resp.\,$S_{k/2}A_n[M]$), $M\geq 0$. In distribution, they are obtained by sending $n$ particles per level in the order $k,k+1,k-1,k+2,\dots$ (resp.\,in the order $k,k-1,k+1,k-2,\dots$). The construction is similar for $\tau_kA^\ast_n[\infty]$ and $S_{k/2}A^\ast_n[\infty]$. 

\begin{prooft}{Proposition~\ref{prop:transfo}}
Since $\tau_k=S_{k/2}\circ S_0$, we only need to show (ii). We only give the proof for $A_n[\infty]$. The proof is similar for $A^\ast_n[\infty]$ and implies the result for $A^\dag_n[\infty]$. 

To do so, it suffices to check that the random sets $A_n[\infty]$ and $S_{k/2}A_n[\infty]$ have the same probability to intersect any given compact set (see \cite[Theorem 2.1.3]{SW}). Let $C\subset \R^2$ be a compact and $\varepsilon>0$. Let $M_0$ be such that for all $M\geq M_0$:
\begin{equation*}
\left\vert\PPP{A_n[\infty]\cap C\neq\emptyset}-\PPP{A_n[M]\cap C\neq\emptyset} \right\vert\leq \varepsilon
\end{equation*} 
and
\begin{equation*}
\left\vert\PPP{S_{k/2}A_n[\infty]\cap C\neq\emptyset}-\PPP{S_{k/2}A_n[M]\cap C\neq\emptyset} \right\vert\leq \varepsilon.
\end{equation*}

Now, let us grow the aggregate $A_n[M]$ by sending $n$ particles per level from $M+1$ to $M+k$. We denote by $A_1$ the resulting aggregate. Similarly, let us grow the aggregate $S_{k/2}A_n[M]$ by sending $n$ particles per level from $k-M-1$ to $-M$. We denote by $A_2$ the resulting aggregate. In both $A_1$ and $A_2$, $n$ particles are sent per level from $-M$ to $M+k$, but not in the same order. Nevertheless, by the {\it Abelian property for finite aggregates}, they are equally distributed. So,
\[\PPP{A_1\cap C\neq\emptyset}=\PPP{A_2\cap C\neq\emptyset}.\] 

The aggregates can be coupled in such a way that
\[A_n[M]\subset A_1\subset A_n[M+k]\]
and
\[S_{k/2}A_n[M]\subset A_2\subset S_{k/2} A_n[M+k],\]
 with probability 1, respectively. This implies that 

\begin{equation}\label{eq:sym1}
\left\vert\PPP{A_n[\infty]\cap C\neq\emptyset}-\PPP{A_1\cap C\neq\emptyset} \right\vert\leq \varepsilon
\end{equation} 
and
\begin{equation}\label{eq:sym2}
\left\vert\PPP{S_{k/2}A_n[\infty]\cap C\neq\emptyset}-\PPP{A_2\cap C\neq\emptyset} \right\vert\leq \varepsilon.
\end{equation}
 Indeed, for \eqref{eq:sym1}, we have
\begin{align*}
\PPP{A_n[\infty]\cap C\neq\emptyset}-\varepsilon&\leq \PPP{A_n[M]\cap C\neq\emptyset}\leq\PPP{A_1\cap C\neq\emptyset}\\&\leq\PPP{A_n[M+k]\cap C\neq\emptyset}\leq \PPP{A_n[\infty]\cap C\neq\emptyset}+\varepsilon.
\end{align*} 
 The bound \eqref{eq:sym2} is obtained by the same way.
 
Finally, we get by collecting bounds that:
 \begin{align*}
 &\left\vert\PPP{A_n[\infty]\cap C\neq\emptyset}-\PPP{(S_{k/2}A_n[\infty])\cap C\neq\emptyset} \right\vert\\
 &\qquad\qquad\qquad\qquad\leq\left\vert\PPP{A_n[\infty]\cap C\neq\emptyset}-\PPP{A_1\cap C\neq\emptyset} \right\vert\\&\qquad\qquad\qquad\qquad\qquad+\left\vert\PPP{A_1\cap C\neq\emptyset}-\PPP{A_2\cap C\neq\emptyset} \right\vert\\&\qquad \qquad\qquad\qquad\qquad
+\left\vert\PPP{A_2\cap C\neq\emptyset}-\PPP{(S_{k/2}A_n[\infty])\cap C\neq\emptyset} \right\vert \\
&\qquad\qquad \qquad\qquad\leq 2\varepsilon. 
\end{align*}
\end{prooft}

\subsubsection{Mean size of the aggregates per level}
As a consequence of Proposition~\ref{prop:transfo}, the following result shows that the expected width of $A_n[\infty]$ (resp.\,$A^*_n[\infty]$ and $A^\dag_n[\infty]$) equals $n$ when it is restricted to a horizontal line.

\begin{proposition}
\label{prop:ExpWidth}
For any $i\in \Z$, we have:
\begin{enumerate}[(i)]
\item $\EEE{\# A_n[\infty]\cap (\Z\times\{i\})} = n$;
\item $\EEE{\# A^\ast_n[\infty]\cap (\Z\times\{i\})} = n$;
\item $\EEE{\# A^\dag_n[\infty]\cap (\Z\times\{i\})} = n$.
\end{enumerate}
\end{proposition}

{The above result shows that the expected number of sites in $A_n[\infty]\cap (\Z\times\{i\})$ equals the (expected) number of particles emitted from level $i$.
Roughly, this means that it is, in average, as if the particles contribute to the growth of the level at which they are sent.}

\begin{prooft}{Proposition~\ref{prop:ExpWidth}}
We begin with (i). To do it, for any  $i\in \Z$, we denote by $Q(i)$ the number of sites in $A_n[\infty]$ with ordinate $i$, {\it i.e.}
\[Q(i)=\# A_n[\infty]\cap (\Z\times\{i\}).\]
For $j\in \Z$, we also denote by  $Q(i,j)$ the number of sites in $A_n[\infty]\cap \left(\Z\times \{i\}\right)$ which are created by particles from level $j$. Thus $Q(i)=\sum_{j\in \Z} Q(i,j)$ a.s.. According to Proposition~\ref{prop:transfo}, the random variables $Q(i,j)$ and $Q(j,i)$ have the same distribution. Since $\sum_{j\in \Z} Q(j,i)=n$, we obtain
\[
\EEE{Q(i)} = \sum_{j\in\Z} \EEE{Q(i,j)} = \sum_{j\in\Z} \EEE{Q(j,i)} = n ~.
\]

In a similar way, we obtain (ii), (iii) by noticing that $\EEE{N_i}=\EEE{\#\mathcal{N}_i}=n$.
\end{prooft}

To derive Proposition~\ref{prop:ExpWidth}, we used the fact that $\EEE{Q(i,j)}=\EEE{Q(j,i)}$. Such an equality can be understood as a mass transport principle (see {\it e.g.} \cite{BLPS}). {Proposition~\ref{prop:ExpWidth} will be applied to derive a stabilization result (Theorem \ref{cor:StabWRTO}) which itself implies that the random aggregates $A_n^*[\infty]$ and $A_n^\dag[\infty]$ are a.s.\, countable unions of \textit{finite} connected components (Corollary \ref{cor:finitecc}). In a similar way, the mass transport principle  was  used in the literature to bound component sizes for other growth processes with infinitely many sources (see e.g. Theorem 1.1. in \cite{AP17})}.

\subsubsection{{Weak stabilization}}
The following {theorem} claims that, given a strip $\Z_M$, all the aggregates $A_n[M']$ (resp.\,$A_n^*[M'$]) coincide with  $A_n[\infty]$ (resp.\,$A_n^*[\infty]$) on  $\Z_M$ when $M'$ is large enough. 

\begin{theorem}
\label{th:existence}
Let $n\geq 1$.
\begin{enumerate}[(i)]
\item A.s., for all integer $M$, there exists $M_0$ such that for all $M' \geq M_0$:
\begin{equation*}
A_n[M'] \, \cap \, \Z_M = A_n[\infty] \, \cap \, \Z_M ~.
\end{equation*}
\item A.s., for all integer $M$, there exists $M_0$ such that for all $M' \geq M_0$,
\begin{equation*}
A^*_n[M'] \, \cap \, \Z_M = A_n^*[\infty] \, \cap \, \Z_M ~.
\end{equation*}
\end{enumerate}
\end{theorem}

\begin{prooft}{Theorem~\ref{th:existence}}
First, notice that the random variable $X_{n,M}=\# A_n[\infty]\cap\Z_M$ is a.s.\,finite since, according to Proposition~\ref{prop:ExpWidth}, we have $\EEE{X_{n,M}} = (2M+1)n$. Now, let $z_1,\dots, z_{X_{n,M}}$ be an enumeration of the sites in $A_n[\infty]\cap\Z_M$. Denote by $t_i$ the level from which the particle which creates $z_i$ is emitted. The conclusion readily follows by setting $M_0=\max \{\vert t_i\vert:\,i=1,\dots,X_{n,M}\}+1$. Note that $M_0$ is a.s.\,finite since $X_{n,M}$ is a.s.\,finite. This proves (i). The same arguments can be used to get (ii). 
\end{prooft}
{Theorem \ref{th:existence} is referred to as a \textit{weak} stabilization result since a stronger version is  established in the next section (see Theorem \ref{th:stabilize}).}

In the proof of Theorem~\ref{th:existence}, we strongly used the fact that $A_n[\infty]$ and $A_n^*[\infty] $ are constructed w.r.t.\,the usual order. Indeed, all the sites of $A_n[\infty] \, \cap \, \Z_M$ (resp. $A_n^*[\infty] \, \cap \, \Z_M)$) are produced by particles which are necessarily sent until a \textit{finite} level $\pm M_0$ since all of the particles which are beyond this level are sent \textit{after} those with level $|i|\leq M_0$. This argument does not hold for $A_n^\dag[\infty]$  because $A_n^\dag[\infty]$ is based on a family of independent PPP's. In particular, for any level $i$, an \textit{infinite} number of particles, beyond $i$, are  sent \textit{before} those with level $i$. However, as we will see in Section~\ref{sec:forest}, Theorem~\ref{th:existence} remains true for $A_n^\dag[\infty]$. The proof will require a specific treatment by introducing a notion of chain of changes.

Note also that the conclusions of Theorem~\ref{th:existence} hold if we apply one of the transformations $\tau_k$ or $S_{k/2}$, $k\in \Z$, to the aggregates and strips.

\section{Far particles do not touch central strips} 
\label{sec:hightpart}

The main result of this section (Theorem~\ref{th:stabilize}) claims that upon adding a new site to the aggregate, each particle from levels $(0,\pm M')$, with $M'$ large enough, does not visit the strip $\Z_M$.  In particular, it provides a (more sophisticated) alternative proof of Theorem~\ref{th:existence} since it ensures that, a.s.\,for $M$ large enough, all the aggregates $A_n[M']$ with $M'> M^{\alpha}$ coincide on $\Z_M$.

Not only Theorem~\ref{th:stabilize} is more precise than Theorem~\ref{th:existence}, it also has its own interest and appears as a fruitful tool throughout this paper. First, it is one of the main ingredients to derive Propositions~\ref{prop:DifferentLaw} and~\ref{prop:mixing}, which themselves are fundamental because they imply that the aggregate $A_n^\dag[\infty]$ only consists of finite connected components (see Corollary~\ref{cor:finitecc}). This fact allows us to define the infinite IDLA forest without controlling tricky chains of changes (see Section~\ref{sec:forest}). Such a control would be a challenging question and contains technical difficulties. Another {results which are based on Theorem~\ref{th:stabilize} are shape theorems (Theorems~\ref{th:ShapeTh} and \ref{th:ShapeThPoisson})}. Proving {Theorem \ref{th:stabilize}} requires fine estimates which cannot be deduced directly from {Theorem~\ref{th:existence}}. The arguments in the proof of Theorem~\ref{th:stabilize} are also extensively re-used in Section~\ref{sec:stabOrigine}.

\begin{theorem}Let $n\geq 1$.
\label{th:stabilize}
\begin{enumerate}[(i)]
\item Let $\alpha>1$. The following property holds with probability $1$. There exists a random integer $M_0=M_0(n)\geq 1$ such that, for any $M\geq M_0$, the trajectory of any particle contributing to $A_n[\infty]$  and starting from $(0,i)$, with  $|i| >M^{\alpha}$, does not visit the horizontal strip $\Z_M$. 
\item The same holds for $A_n^\ast[\infty]$.
\end{enumerate}
\end{theorem}

\subsection{Proof of Theorem~\ref{th:stabilize}, (i)\label{sec:proofStab}}

Let $\alpha>1$. For any integers $M\geq 1$, $k\geq 0$, let 
\begin{equation}
\label{eq:defJMK}
 J_{M,k}=\left\{0\right\}\times\left\llbracket\lfloor M^{\alpha}\rfloor+kM+1,\lfloor M^{\alpha}\rfloor+(k+1)M  \right\rrbracket.
 \end{equation}
Given a realization of $A_{n}[M^\alpha+kM]$, we send $nM$ particles from the interval $J_{M,k}$ first from level $\lfloor M^{\alpha}\rfloor+{k} M+1$, then from level $\lfloor M^{\alpha}\rfloor+kM+2$, and so on.  In what follows, the \textit{current aggregate} associated with a particle $P$ denotes the aggregate which is  produced just before sending $P$. Now, let $E_{M,k}$ be the following event:   
\begin{eqnarray}
\label{DefEMk}
E_{M,k} = \left\lbrace \begin{array}{c}
\mbox{ At least one of the $nM$ particles starting from $J_{M,k}$  }\\
\mbox{ visits the strip $\Z_M$ before exiting the current aggregate}.\\
\end{array}
\right\rbrace ~.
\end{eqnarray}  
 The event $E_{M,k}$ describes the unpleasant situation where particles started far away from the origin, more  precisely from $J_{M,k}$, may modify the aggregate close to the origin. The following lemma shows that $E_{M,k}$ occurs with small probability.

\begin{lemma}
\label{le:ProbaEMk}
There exist constants $0<\eta<1$, $c_1,c_2>0$ such that, for {any $M,n\geq 1$ and for any $k\geq 0$
\[\PPP{E_{M,k}} \leq nM (1-\eta)^{\frac{M^{\alpha}+(k-2)M-c_1n}{c_2n^{2}M}}.\]}
\end{lemma}
{Since $\alpha>1$, it follows from Lemma~\ref{le:ProbaEMk} that 
\begin{equation*}
\sum_{M\geq 1}\PPP{\bigcup_{k\geq 0}E_{M,k}}  \leq \sum_{M\geq 1} \sum_{k\geq 0} \PPP{E_{M,k}} < \infty.
\end{equation*}
}
This together with the Borel-Cantelli lemma implies Theorem~\ref{th:stabilize} (i). 

For $x\in\R^2$ and $r>0$, we denote by $B(x,r)$ the intersection of $\Z^2$ and of the Euclidean {(closed)} ball, centered at $x$ with radius $r$. We set $B(r)=B(0,r)$ { and we write $A+B$ for the sum of the sets $A$ and $B$, i.e.\,$A+B=\{a+b:\,a\in A,\, b\in B\}.$}

\begin{prooft}{Lemma~\ref{le:ProbaEMk}  }
First, we  introduce some notation. For any $\ell\geq 0$, let  $B_{\ell}=J_{M,k}{+} B(n+\ell r)$ and $S_{\ell+1}=B_{\ell+1}\setminus B_{\ell}$, where $r>0$ will be specified later. Notice that if a particle starting from $J_{M,k}$ meets the strip $\Z_M$ then necessarily it crosses the annulus $S_\ell$, for any $\ell\leq \ell_{\mbox{\tiny{max}}}$, with 
\begin{equation}
\label{eq:deflmax}
 \ell_{\mbox{\tiny{max}}} = \left\lfloor r^{-1} \left(\lfloor M^{\alpha}\rfloor+(k-1)M+1 - n \right) \right\rfloor.
 \end{equation}

One of the key ingredients to derive Lemma~\ref{le:ProbaEMk} is the following result, referred to as the \textit{crossing lemma}.
{\begin{lemma}[Crossing Lemma]
\label{le:crossinglemma}
There exist $\eta_0,\eta>0$ such that for any $S,V\subset\Z^{2}$ and $r>0$ satisfying $S\subset V{+} B(r)$ and $\# S\setminus V \leq \eta_0 r^{2}$, and for any particle $\xi^{x}$ starting from $x\in V$ and stopped upon exiting $V{+} B(r)$, the following inequality holds:
\begin{equation}
\label{TechnicalIneq}
\PP \big( \xi^{x} \cap \big( (V{+} B(r)) \setminus (S\cup V) \big) \not= \emptyset \big) \geq \eta ~.
\end{equation}
\end{lemma}
{In the above result, with a slight abuse of notation, we have identified the particle $\xi^x$ to its range. Lemma \ref{le:crossinglemma}} is an adaptation of \cite[Lemma 3.2]{DLYY} written in the context of $\Z^2$. It expresses the difficulty for a particle to cross an annulus when the aggregate occupies only a small portion of it.
}

{To apply Lemma \ref{le:crossinglemma}}, the width $r$ is chosen in such a way that  $nr\leq \eta_0 r^{2}$. More precisely, we let {\[r=\frac{4n}{\eta_0}.\]}   Now, to be in the framework of Lemma~\ref{le:ProbaEMk}, we introduce the notion of good (resp. bad) annuli as follows. Given a realization of $A_{n}[M^\alpha+kM]$, we say that  $S_{\ell}$ is \textit{good} if
\[
\# \left( A_{n}[M^\alpha+kM]\cap S_{\ell} \right) \leq \eta_0 r^{2}. 
\] Otherwise, we say that the annulus is \textit{bad}.  Notice that the number $N^{\mbox{\tiny{bad}}}$ of bad annuli is deterministically bounded since a.s.
\begin{equation*}
\label{Ngood}
 N^{\mbox{\tiny{bad}}} \times \eta_0 r^{2} \leq \# A_{n}[M^\alpha+kM] \leq n \left( 2(M^{\alpha}+kM)+1 \right).
\end{equation*}
Denoting by $N^{\mbox{\tiny{good}}}=\ell_{\mbox{\tiny{max}}} - N^{\mbox{\tiny{bad}}} $ the number of good annuli, it follows from the above inequality and from the definitions of $r$, $\ell_{\mbox{\tiny{max}}}$  that 
\begin{equation*}
N^{\mbox{\tiny{good}}} \geq r^{-1}\left( \frac{1}{2}M^\alpha+\frac{1}{2}kM-M\right) -\frac{\eta_0}{4}-\frac{\eta_0}{16n}-1.
\end{equation*}
Thus
\begin{equation}
\label{eq:proportioNgood}
N^{\mbox{\tiny{good}}} \geq  \frac{M^\alpha+(k-2)M}{2r}-c_0,
\end{equation}
where $c_0=1+\frac{5\eta_0}{16}$. 

Now, for $i=1,\ldots,nM$, let us call particle $i$ the $i$-th particle which is sent from the interval $J_{M,k}$. The main idea to prove Lemma~\ref{le:ProbaEMk} is to use the fact that, if particle $i$ hits the strip $\Z_M$ then,   according to Equation \eqref{eq:proportioNgood}, it necessarily has to cross a large number of good annuli. But, {according to} Lemma~\ref{le:crossinglemma}, for each new good annulus that particle $i$ meets, it has probability at least $\eta$ to be stuck inside. Roughly, this will imply that the event $E_{M,k}$ cannot occur with high probability. 

More precisely, given a realization of $A_{n}[M^\alpha+kM]$, let $X_1$ be the number of good annuli which are crossed by particle $1$. Given $X_1,\ldots,X_i$, with $i\geq 1$, let $X_{i+1}$ be the number of good annuli which are crossed by particle $i+1$, \textit{but not} already crossed by particles $1,\ldots,i$. {Notice that after the launches of the first $i$ particles, $i$ new sites are added to the aggregate. With nonnull probability, some of these particles will settle on the first good annuli. Hence, some of these good annuli will contain more than $\eta_0 r^2$ sites of the current aggregate. To deal with this problem, we choose to only focus  on \textit{new good annuli} for the random variable $X_{i+1}$}.   Notice that the random variables $X_i$, $i\geq 1$, are not identically distributed and that the sequence $(\PPP{X_i=0})_{i\geq 1}$ is non-decreasing.  The following lemma shows that the random variable $X_i$ is stochastically dominated by a geometric distribution with parameter independent of $i$.

\begin{lemma}
\label{le:InductionXi}
There exists $\eta\in (0,1)$ such that, for any $i\geq 1$ and for any $s,t\geq 0$, the following properties hold: 
\begin{itemize}
\item[(i)] \label{le:inductiona} $\PP \big( X_1>t \, \big| \, A_{n}[M^\alpha+kM] \big) \leq (1-\eta)^{t}$;
\item[(ii)] \label{le:inducationb} $\PP \big( X_{i+1}>t \, \big| \, A_{n}[M^\alpha+kM] , \, \sum_{j\leq i} X_j \leq s \big) \leq (1-\eta)^{t}$.
\end{itemize}
\end{lemma}

If one of the $nM$ particles starting from $J_{M,k}$ hits the strip $\Z_M$ then the sum of $X_i$'s is necessarily larger than $N^{\mbox{\tiny{good}}}$. Thus, according to \eqref{eq:proportioNgood},  we have 
\[E_{M,k} \subset \left\{  \sum_{i=1}^{nM} X_i \geq \frac{M^\alpha+(k-2)M}{2r}-c_0    \right\}.  \] 
This implies that 
\begin{align}
\label{SumXi}
\PPP{E_{M,k}} \leq \PPP{\sum_{i=1}^{nM} X_i \geq \frac{M^\alpha+(k-2)M}{2r}- c_0 \;\left|\; \sum_{i=1}^{nM-1} X_i \leq \frac{nM-1}{nM} \left(\frac{M^\alpha+(k-2)M}{2r}- c_0 \right)  \right.}\\
 + \PPP{\sum_{i=1}^{nM-1} X_i \geq \frac{nM-1}{nM} \left(\frac{M^\alpha+(k-2)M}{2r}- c_0 \right) }\nonumber . 
\end{align}
Thanks to Assertion (ii) in Lemma~\ref{le:InductionXi}, we can bound the conditional probability in \eqref{SumXi} by $(1-\eta)^{\frac{M^{\alpha}+(k-2)M-2c_0r}{2rnM}}$. By induction and Lemma~\ref{le:InductionXi}, we get
\[
\PP	(E_{M,k}) \leq nM (1-\eta)^{\frac{M^{\alpha}+(k-2)M-2c_0r}{2rnM}}.
\]
This concludes the proof of Lemma~\ref{le:ProbaEMk} by taking $c_1=\frac{8c_0}{\eta_0}$ and $c_2=\frac{8}{\eta_0}$.
\end{prooft}

\begin{prooft}{Lemma~\ref{le:crossinglemma}}
It relies on an adaptation of \cite[Lemma 3.2]{DLYY}. For sake of completeness, we recall the main arguments of \cite{DLYY}.  Let $Y$ be the outer boundary of $V{+} B(r/2)$ in $\Z^2$, {\it i.e.} 
\[
Y = \left\{ y \in \Z^2 \, : \; y\notin V{+} B(r/2) \, \mbox{ and } \, \exists y' \in V{+} B(r/2) \,\mbox{ s.t. } \, |y-y'|_{1}=1 \right\} ~,
\] 
where $|\cdot |_1$ denotes the $1$-norm on $\Z^2$, {\it i.e.} $|z|=|z(1)|+|z(2)|$ for any $z=(z(1),z(2))\in \Z^2$. First, notice that every path from $x\in V$ to the complement of $V{+} B(r)$ must hit $Y$. Thus, by Markov's property, it suffices to prove (\ref{TechnicalIneq}) for random walks starting from $Y$.

Let $y\in Y$, $B=B(y,r/3)$ and $Q=B\setminus S$. According to \cite[Lemma 3.1]{DLYY}, for any $t>0$, we have 
\begin{align}
\label{IneqLemma3.1}
\PPP{ \xi^{y} \cap \big( (V{+} B(r)) \setminus (S\cup V) \big) \not= \emptyset} & \geq \PPP{ \xi^{y} \cap Q \not= \emptyset}\notag \\
 & \geq  \PPP{ B \subset A_{t}(y\mapsto r/3) \big)} \times \frac{\# Q}{t},
\end{align}
where $A_{t}(y\mapsto r/3)$ denotes the aggregate obtained by letting $t$ particles starting from $y$ and stopped upon exiting $B(y,r/3)$.

To deal with \eqref{IneqLemma3.1}, recall that from \cite[Section 2]{DLYY}, there exists $\alpha>0$ such that
\[
\PP\big( B \subset A_{t}(y\mapsto r/3) \big) \geq \alpha
\]
for $t=\# B(y,r/(3\alpha))$. The previous inequality is referred to as the \textit{weaker lower bound} in  \cite{DLYY}. Now, let $c>0$ be such that $\# B \geq c(r/3)^{2}$. Since $B\subset V^c$ and $\#S\setminus V\leq \eta_0 r^2$, we have 
\[
\# Q = \# B - \# B\cap S \geq c(r/3)^{2} - \eta_0 r^{2} = \eta_0 r^{2}, 
\]
with $\eta_0=c/18$. Taking $C>0$ in such a way that $t\leq C(r/(3\alpha))^{2}$, we have
\[\PPP{ B \subset A_{t}(y\mapsto r/3) \big)} \times \frac{\# Q}{t}\geq  \alpha \times \frac{\eta_0 r^{2}}{C(r/(3\alpha))^{2}}=:\eta.\]
\end{prooft}

\begin{prooft}{Lemma~\ref{le:InductionXi}}
{One of the key ingredients to derive Lemma~\ref{le:InductionXi} is Lemma \ref{le:crossinglemma}.}

First, we prove (i).  To do it, we first show that, for any $t\geq 0$, 
\begin{equation}
\label{X1Geom}
\PP \big( X_1>t \, \big| \, A_{n}[M^\alpha+kM] , \, X_1>t-1 \big) \leq 1-\eta ~.
\end{equation}

{Given a realization of $A_{n}[M^\alpha+kM]$, we denote by $T$ the index of the $\lfloor t\rfloor+1$-th good annulus which is reached by particle 1.  Now, let \[S= A_{n}[M^\alpha+kM]\cap B_{T} \quad \text{and} \quad  V\equiv B_{T-1}.\] Notice that $\# S\setminus V \equiv \# A_{n}[M^\alpha+kM]\cap S_{T} \leq \eta_0 r^{2}$ since $S_T$ is a good annulus. Conditional on the event $\{X_1>t-1\}$,  if $X_1>t$, then necessarily particle $1$ crosses the annulus $S_{T}$. According to Lemma~\ref{le:crossinglemma}, the event $\{X_1>t\}$ occurs with probability smaller than $1-\eta$, which proves  \eqref{X1Geom}. By induction and because 
\begin{multline*}
\PPP{ X_1>t \, | \, A_{n}[M^\alpha+kM]} \\
   \leq  \PPP{ X_1>t \, | \, A_{n}[M^\alpha+kM] , \, X_1>t-1} \; \PPP{ X_1>t-1 \, | \, A_{n}[M^\alpha+kM]},
\end{multline*}
we get Assertion (i).}

Now, we prove (ii). We proceed in the same spirit as above. To do it,  it is sufficient to show that for any {$t\geq 0$},
\[
\PPP{ X_{i+1}>t \, \left| \, A_{n}[M^\alpha+kM] , \, \sum_{j\leq i} X_j \leq s , \; X_{i+1}>t-1 \right. } \leq 1-\eta.
\]
Given $X_1,\ldots,X_i$ and $A_{n}[M^\alpha+kM]$, we denote by $T$ the index of the {$\lfloor t\rfloor+1$}-th good annulus which is counted from the {$\lfloor s\rfloor+1$}-th good one. Conditional on the event $\left\{X_{i+1}>t-1\right\}$, if $X_{i+1}>t$, then particle $i+1$ has to cross the good annulus $S_T$ through the current aggregate $\tilde{A}$, which has been augmented from $A_{n}[M^\alpha+kM]$ by exactly $i$ sites corresponding to the first $i$ particles sent from $J_{M,k}$. But, by definition of $X_{i+1}$, the annulus $S_{T}$ has not been visited by the first $i$ particles. {In particular, $\tilde{A}\cap S_T$ equals} $A_{n}[M^\alpha+kM]\cap S_T$. Since $\#A_{n}[M^\alpha+kM]\cap S_T$  is smaller than $\eta_0 r^{2}$, we can apply Lemma~\ref{le:crossinglemma}. The end of the proof follows from the same lines as in (i).  
\end{prooft}

\subsection{Proof of Theorem~\ref{th:stabilize}, (ii)}
\label{sec:proofStab*}

This proof will be sketched because it relies on a simple adaptation of the proof of Theorem~\ref{th:stabilize}, (i). The main difference is that we have to provide estimates for the number of particles which are sent per site.

Let $\alpha>1$, $M,n\geq 1$ and $k\geq 0$. We have to show that the series  $\sum_{M\geq 1}\PPP{\bigcup_{k\geq 0}E_{M,k}^{\ast}}$ is convergent, where the event $E_{M,k}^{\ast}$ is defined in the same spirit as \eqref{DefEMk} by replacing the aggregate {$A_{n}[M^\alpha+kM]$} by {$A_{n}^\ast[M^\alpha+kM]$}. To control the size of the aggregate $A_{n}^{\ast}[M^{\alpha}+kM]$, let 
\begin{equation*}
\label{DefFMk}
F_{M,k} = \left\{ \sum_{|i|\leq M^{\alpha}+kM} N_i \, \leq 2 n \left( 2(M^{\alpha}+kM)+1 \right) \right\} ~.
\end{equation*} 

On the event $F_{M,k}$, we can adapt the main arguments of Section~\ref{sec:proofStab}  by considering good and bad annuli $S_{\ell}$, $\ell\geq 1$, with width $r$ as follows.  First, recall that each particle starting from $J_{M,k}$ has to cross $\ell_{\mbox{\tiny{max}}}$ annuli to hit the strip $\Z_M$, where $J_{M,k}$ and  $\ell_{\mbox{\tiny{max}}}$ are defined in \eqref{eq:defJMK} and \eqref{eq:deflmax}, respectively. Then, conditional on $A_{n}^{\ast}[M^{\alpha}+kM]$, we say that the annulus $S_{\ell}$ is \textit{good} if \[\# \left( A_{n}^{\ast}[M^{\alpha}+kM]\cap S_{\ell}\right)\leq\eta_0 r^{2}\] and we say that it is \textit{bad} otherwise. Since the aggregate $A_{n}^{\ast}[M^{\alpha}+kM]$ potentially contains twice more particles than $A_{n}[M^{\alpha}+kM]$, the width $r$ has to be chosen as twice the width appearing in the proof of (i). Thus we let  $r=8n/\eta_0$. On the event $F_{M,k}$, the number $N^{\mbox{\tiny{bad}}}$ of bad annuli is still deterministically bounded since, a.s., 
\[
N^{\mbox{\tiny{bad}}} \times \eta_0 r^{2} \leq \# A_{n}^{\ast}[M^{\alpha}+kM] \leq 2 n \big( 2(M^{\alpha}+kM)+1 \big)
\]
and thus $N^\text{good}$ can be bounded in the same spirit as \eqref{eq:proportioNgood}.

Now, to control the number of particles starting from $J_{M,k}$, let 
\[G_{M,k} = \left\{ \sum_{(0,i)\in J_{M,k}} N_i \, \leq 2\; {\max\left\{nM, k^\varepsilon\right\}}    \right\},\] {for some $\varepsilon \in (0,1)$.}
Proceeding exactly in the same spirit as in the proof of Lemma~\ref{le:ProbaEMk}, we obtain that, for some constants $c'_1,c'_2>0$: 
\[\PPP{\left. E_{M,k}^{\ast}\right|F_{M,k}\cap G_{M,k}} \leq 2\; {\max\left\{nM, k^\varepsilon\right\}} (1-\eta)^{\frac{M^{\alpha}+(k-2)M-c'_1n}{c'_2n{\times \max\left\{nM, k^\varepsilon\right\}}}}.\]
Moreover, we can easily prove that the series $\sum_{M\geq 1} \sum_{k\geq 0} \PPP{F_{M,k}^{c}}$ and $\sum_{M\geq 1} \sum_{k\geq 0} \PPP{G_{M,k}^{c}}$ are finite. Therefore $\sum_{M\geq 1}\PPP{\bigcup_{k\geq 0}E_{M,k}^{\ast}}$ is finite, which concludes the proof of Theorem~\ref{th:stabilize}, (ii).

\begin{Rk}\label{rk:thinner} To prove Theorem~\ref{th:stabilize}, (ii), we introduced the event $F_{M,k}$, which roughly means that the aggregate $A_{n}^{\ast}[M^{\alpha}]$ is thin. In particular, we used the fact that the thinner the aggregate $A_{n}^{\ast}[M^{\alpha}]$ is (which appears when the sum of Poisson random variables $N_i$, $|i|\leq M^\alpha+kM$, is not big), the bigger is the probability that the particles, which start from levels $|i|>M^{\alpha}$, do not get the strip $\Z_M$. This remark will be used in the proof of Proposition~\ref{prop:DifferentLaw}.
\end{Rk}

\subsection{$A_n^*[\infty]$ avoids $\Z\times\{0\}$ with positive probability \label{sec:DiffLaw}}

The next result ensures that the infinite aggregates $A_n[\infty]$ and $A_n^*[\infty]$ are not identically distributed. The first one contains, by construction, the vertical axis $\{0\}\times\Z$ with probability $1$, whereas the second one does not intersect the axis $\Z\times\{0\}$ with positive probability.

\begin{proposition}
\label{prop:DifferentLaw}
Let $n\geq 1$. With positive probability, the aggregate\, $A_n^*[\infty]$ does not intersect the axis $\Z\times\{0\}$. 
\end{proposition}

\begin{prooft}{Proposition~\ref{prop:DifferentLaw}}
{To prove that $\PPP{A_n^*[\infty]\cap (\Z\times\{0\})=\emptyset}$ is positive, it is sufficient to show that, for $M$ large enough,  
\begin{equation}
\label{eq:aimdifferentlaw}
 \PPP{\bigcap_{k\geq 0}\left(E_{M,k}^{\ast}\right)^c   \cap \left\{\sum_{|i|\leq M^\alpha}N_i=0\right\}}>0,
 \end{equation}
where $E_{M,k}^{\ast}$ is as in the proof of Theorem~\ref{th:stabilize}, (ii). To do it, we recall that the sums $\sum_{M\geq 1} \sum_{k\geq 0} \PPP{E_{M,k}^{\ast}}$ and $\sum_{M\geq 1} \sum_{k\geq 0} \PPP{F_{M,k}^c}$ are finite (see the proof of Theorem~\ref{th:stabilize}, (ii)). 
Thus, for $M$ large enough, we have 
\begin{equation}
\label{CsqTh*}
\PPP{\bigcap_{k\geq 0}\left(E_{M,k}^{\ast}\right)^c\cap  \left\{\sum_{|i|\leq M^\alpha}N_i\leq 2n (2M^\alpha+1)\right\}}>0. 
\end{equation}
Now, recall that the aggregate $A_n^*[\infty]$ is constructed w.r.t.\,two independent random variables:  the first one, say $\omega_1$ (resp. the second one denoted by $\omega_2$) concerns random walks starting from $I(M^\alpha) = \{0\}\times \llbracket- M^\alpha, M^\alpha\rrbracket$ (resp. starting from $(\{0\}\times \Z)\setminus I(M^\alpha)$) and Poisson random variables $N_i$ indexed by the set of $|i|\leq M^\alpha$ (resp. by the set of $|i|>M^\alpha$), {\textit{i.e.}    
\[\omega_1=\left\{ \xi^i_j: \;  |i|\leq M^\alpha, j\leq N_i  \right\}\]
and 
\[\omega_2=\left\{ \xi^i_j: \;  |i|> M^\alpha, j\leq N_i  \right\},\]
where $\{N_i: \;  i\in \Z\}$ is a family (which is independent of the $\xi^i_j$'s) of independent Poisson random variables with parameter $n$ and where $\{\xi^i_j: \; i\in \Z, j\in \N\}$ is a family of independent simple random walks, with $\xi^i_j$ starting at vertex $(0,i)$. 
}   In what follows, for any set $A\subset \Z^2$, we write
{\begin{eqnarray*}
\mathcal{E}(A)= \left\lbrace \begin{array}{c}
\mbox{$\omega_2$: particles associated with $\omega_2$ (starting from $ |i|>M^\alpha$) do not visit }\\
\mbox{the strip $\Z_{M}$ when they are used to grow the initial aggregate $A$ }
\end{array}
\right\rbrace ~.
\end{eqnarray*}}
Notice that \[  \bigcap_{k\geq 0}\left( E^*_{M,k}  \right)^c =\left\{(\omega_1,\omega_2): \omega_2\in \mathcal{E}\left(A_n^*[M^\alpha](\omega_1)\right)\right\}. \]
Remark~\ref{rk:thinner} expresses that $\PPP{\mathcal{E}(A)}$ is decreasing w.r.t.\,$A$, {\it i.e.} for any $A,A'\subset \Z^2$,
\begin{equation}
\label{eq:deacreasingimplication}
A \subset A' \, \Longrightarrow \, \PPP{\mathcal{E}(A)} \geq \PPP{\mathcal{E}(A')}. 
\end{equation}

We are now prepared to prove \eqref{eq:aimdifferentlaw}. Indeed, since $\sum_{|i|\leq M^\alpha}N_i=0$ if and only if  $A_n^*[M^\alpha]=\emptyset$, we have 
\begin{align*}
 \PPP{\bigcap_{k\geq 0}\left(E_{M,k}^{\ast}\right)^c   \cap \left\{\sum_{|i|\leq M^\alpha}N_i=0\right\}} & = \int\ind{\omega_2\in \mathcal{E}\left(A^*_n[M^\alpha](\omega_1)  \right)}  \ind{\sum_{|i|\leq M^\alpha}N_i(\omega_1)=0}\mathrm{d}\omega_1\mathrm{d}\omega_2\\
 & =  \int\ind{\omega_2\in \mathcal{E}\left( \emptyset  \right)}\mathrm{d}\omega_2\times \PPP{\sum_{|i|\leq M^\alpha}N_i=0}.
\end{align*}
The second term of the last equality is positive. For the first one, according to \eqref{eq:deacreasingimplication}, we know that 
\begin{align*}
\int\ind{\omega_2\in \mathcal{E}\left( \emptyset  \right)}\mathrm{d}\omega_2 & \geq \int\ind{\omega_2\in \mathcal{E}\left( A_n^*[M^\alpha](\omega_1)\right)} \ind{\sum_{|i|\leq M^\alpha}N_i(\omega_1)\leq 2n(2M^\alpha+1)} \mathrm{d}\omega_1\mathrm{d}\omega_2\\
& = \PPP{\bigcap_{k\geq 0}\left(E_{M,k}^{\ast}\right)^c\cap  \left\{\sum_{|i|\leq M^\alpha}N_i\leq 2n (2M^\alpha+1)\right\}},
\end{align*}
which is positive according to \eqref{CsqTh*}. This concludes the proof of Proposition  \ref{prop:DifferentLaw}. 
}

\end{prooft}

\section{Central particles do not touch far levels}  
\label{sec:stabOrigine}

In this section, we show that the aggregates $A_n[\infty]$ and $A_n^*[\infty]$  above some random levels a.s.\,do not depend on particles which are sent around the origin. Such a property is one of the key ingredients to prove that the aggregates $A_n[\infty]$, $A_n^*[\infty]$ and $A_n^\dag[\infty]$ satisfy a mixing property (Proposition~\ref{prop:mixing}).

Let $n,M\geq 0$. We define two random aggregates which are coupled to $A_n[\infty]$ as follows. First, we write  $A_{n,M}(M)=A_n[M]$ and  $B_{n,M}(M)=\emptyset$. Then we launch {$n$ particles per level} from levels $M+1$, $M+2$ and so on by using the same random walks as $A_n[\infty]$. In particular, we couple the growth (only from the top) of the aggregates based on $A_n[M]$ and $\emptyset$. For any $t>M$,  we denote by $A_{n,M}(t)$ (resp. $B_{n,M}(t)$) these two aggregates after sending particles from levels $M+1$, thus $M+2$ etc. till level $t$ included.   Using the same arguments than those used in Theorem~\ref{th:stabilize}, we prove that the increasing sequences $(A_{n,M}(t))_{t\geq M}$ and $(B_{n,M}(t))_{t\geq M}$ a.s.\,converge to infinite aggregates denoted respectively by $A_{n,M}(\infty)$ and $B_{n,M}(\infty)$. With these notations, we emphasize the dependence of those aggregates on parameters $n,M$ and the level $t$ is referred to as the time. 

{Now, let $(N_i)_{i\in \Z}$ be a family of independent Poisson random variables with parameter $n$.  In a similar way, we define two aggregates based on $A_{n,M}^*(M)=A_n^*[M]$ and  $B_{n,M}^*(M)=\emptyset$, by launching (the same) particles from levels $M+1$, $M+2$ and so on, with $N_i$ particles per level $i$. We denote the aggregates at time $t>M$ (resp. the limits of the aggregates) by  $A_{n,M}^*(t)$ and $B_{n,M}^*(t)$ (resp. $A_{n,M}^*(\infty)$ and $B_{n,M}^*(\infty)$). 

For any $K\geq 0$, we let $\mathcal{H}_K = \{ (x,y)\in\Z^2 : y\geq K \}$. We are now prepared to state the main result of this section.

\begin{theorem}
\label{cor:StabWRTO}
Let $n\geq 1$. 
\begin{enumerate}[(i)]
\item  A.s.\,there exists a (random) integer $M_0=M_0(n)\geq 0$ such that, for any $M\geq M_0$, there exists $t_0=t_0(n,M)\geq M$ such that, for any $t\geq t_0$,
\begin{equation*}
\label{Stab2bis}
A_{n}[\infty]\cap \mathcal{H}_t = B_{n,M}(\infty)\cap \mathcal{H}_t.
\end{equation*}
\item The same property holds if we replace $A_{n}[\infty]$ and $B_{n,M}(\infty)$ by $A_{n}^*[\infty]$ and $B_{n,M}^*(\infty)$ respectively.
\end{enumerate}
\end{theorem}
Theorem~\ref{cor:StabWRTO}} asserts that the infinite aggregates {$A_{n}[\infty]$ and $A^{\ast}_{n}[\infty]$} restricted to the half-plane $\mathcal{H}_t$ do not depend on particles starting from levels $i<M$. We call this property {\it Stabilization w.r.t.\,the origin}. It is a consequence of the following proposition.

\begin{proposition}
\label{prop:StabWRTO}
Let $n\geq 1$. 
\begin{enumerate}[(i)]
\item  For any $M\geq 1$, a.s.\,there exists $t_0$ such that, for any $t\geq t_0$, we have 
\[A_{n,M}(\infty)\cap \mathcal{H}_t = B_{n,M}(\infty)\cap \mathcal{H}_t.\]
\item The same property holds if we replace $A_{n,M}(\infty)$ and $B_{n,M}(\infty)$ by $A_{n,M}^*(\infty)$ and $B_{n,M}^*(\infty)$ respectively.
\end{enumerate}
\end{proposition}

In what follows, we detail the proof of the above results for the Poisson case because the fact that $A_n^*[\infty]$ does not intersect horizontal lines with positive probability (Proposition~\ref{prop:DifferentLaw}) leads to a more natural strategy. However, we explain how we adapt the proof for the deterministic case when modifications are required.

Obviously, by symmetry, a similar result to Theorem~\ref{cor:StabWRTO} but for negative levels also holds, {\it i.e.} aggregates $A_{n}[\infty]$ and $A^{\ast}_{n}[\infty]$ restricted to the half-plane $\{(x,y)\in\Z^2 : y\leq -t\}$ do not depend on particles starting from levels $i>-M$, for $M,t$ large enough.

\begin{prooft}{Theorem~\ref{cor:StabWRTO} } Let $n\geq 1$. {Recall that the aggregate $A^{\ast}_{n}[\infty]$ (resp. $A^{\ast}_{n,M}(\infty)$) is based on $A^{\ast}_n[M]$ and on particles which are sent from levels $|i|>M$ (resp. from levels $i>M$),  {\it i.e.}\,from the top and from the bottom w.r.t.\,the usual order (resp.  from the top). The aggregates are coupled in the sense that they are based on the same particles.}  Theorem~\ref{th:stabilize} states that a.s.\,there exists $M_0$ such that for any $M\geq M_0$, particles from levels $i<-M^2$ (taking $\alpha=2$ for instance) do not visit the horizontal strip $\Z_M$. This means that both aggregates $A^{\ast}_{n}[\infty]$ and $A^{\ast}_{n,M^{2}}(\infty)$ coincide on the hyperplane $\mathcal{H}_M$. This together with Proposition~\ref{prop:StabWRTO}  implies Theorem~\ref{cor:StabWRTO}. 
\end{prooft} 

The end of this section is devoted to the proof of Proposition~\ref{prop:StabWRTO}.    Thanks to our coupling, the inclusion $B^{\ast}_{n,M}(t) \subset A^{\ast}_{n,M}(t)$ holds a.s.\,for any $t\geq M$, and we denote by $\Delta(t)$ the random (symmetric) difference between these aggregates, {\it i.e.}
\[
\Delta(t) = A^{\ast}_{n,M}(t) \!\setminus\! B^{\ast}_{n,M}(t) ~.
\] 

Let us explain how the difference set $\Delta(t)$ evolves over time. Let $x\in\Delta(t)$.  Let $t'\geq t$ be the first level (if it exists) for which a particle, {say $P$}, starting from level $t'$ reaches the site $x$ and thus exits {the current aggregate}  through some site $y$. Here, ``current aggregate'' denotes the aggregate which is produced just before sending $P$ { in the construction of $A^{\ast}_{n,M}(t')$}. Just after $P$ is sent, the site $x$ is added to $B^{\ast}_{n,M}(t')$ and is no longer a difference between $A^{\ast}_{n,M}(\cdot)$ and $B^{\ast}_{n,M}(\cdot)$, but $y$  becomes a new difference between both.    

 {If there is no particle starting from level $t'\geq t$ which visits $x$, then } $x$ is a {difference} forever, {\it i.e.}\,$x\in \bigcap_{t'\geq t} \Delta(t')$. Although the symmetric difference evolves, its cardinality remains constant over time and equals
\[
\#\Delta(0) = \# \left( A^{\ast}_{n,M}(0) \!\setminus\! B^{\ast}_{n,M}(0) \right) = \# A^{\ast}_{n}[M].
\]
The above integer is random and Poisson distributed with parameter $(2M+1)n$.    

Notice that  $\big( \cup_{t'\geq M} \Delta(t') \big) \cap \mathcal{H}_{t} = \emptyset$  if and only if the aggregates $A^{\ast}_{n,M}(t')$ and $B^{\ast}_{n,M}(t')$ coincide on the hyperplane $\mathcal{H}_{t}$ for any time $t'$. The same holds for their limits. Henceforth, to prove Proposition~\ref{prop:StabWRTO}, (ii) it suffices to prove that for any $M\geq 1$,
\begin{equation}
\label{DeltaHyperK}
\text{a.s.\,, for $t$ large enough,} \; \left( \bigcup_{t'\geq M} \Delta(t') \right) \cap \mathcal{H}_t = \emptyset.
\end{equation}

To do it, our strategy consists in bounding the growth of the aggregate $(A^{\ast}_{n,M}(t))_{t\geq M}$ from the top at some specific (and random) times $t$. In the same spirit as \cite{LS}, {we introduce the \textit{vertical spreading} of the aggregate $A^{\ast}_{n,M}(t)$ (w.r.t. $t$)} as
\begin{equation*}
\label{eq:defheight}
h \left( A^{\ast}_{n,M}(t) \right) = \max \left\{ y\in\Z : \, (x,y) \in A^{\ast}_{n,M}(t) \right\} \, - \, t ~.
\end{equation*}

The above quantity is the {difference, w.r.t.\,the $y$ axis,} between the highest vertex of $A^{\ast}_{n,M}(t)$ and the level $t$ at which the last particles are sent. Roughly speaking, $h(A^{\ast}_{n,M}(t))$ expresses how much the aggregate $A^{\ast}_{n,M}(t)$ drools beyond level $t$. Given $\zeta>0$, let us set $\tau_1(\zeta)$ as the first time $t>M$ at which the vertical spreading of $A^{\ast}_{n,M}(t)$ becomes smaller than $\zeta$:
\[
\tau_1(\zeta) = \inf \{ t>M : \, h \big( A^{\ast}_{n,M}(t) \big) \leq \zeta \} ~.
\]
Thus, by induction, we define a sequence of random times $(\tau_m(\zeta))_{m\geq 1}$ as follows:
\[
\tau_{m+1}(\zeta) = \inf \{ t>\tau_{m}(\zeta) : \, h \big( A^{\ast}_{n,M}(t) \big) \leq \zeta \}, 
\] 
with the convention $\tau_{m}(\zeta)=\infty$ implies $\tau_{m+1}(\zeta)=\infty$. The next result claims that a.s.\,infinitely often, the vertical spreading of the aggregate $A^{\ast}_{n,M}(\cdot)$ is smaller than $\zeta$.  

\begin{proposition}
\label{prop:FiniteTau}
Let $n,M\geq 1$. Then there exists a positive real number $\zeta=\zeta(n)$ such that, a.s., the random times $(\tau_{m}(\zeta))_{m\geq 1}$ are all finite.
\end{proposition}

As we will see in Section~\ref{sect:TauToSigma}, the choice of $\zeta$ depends on $n$ but not on $M$. Proposition~\ref{prop:FiniteTau} is one of the main ingredients to prove \eqref{DeltaHyperK}.    

\subsection{ Proposition~\ref{prop:FiniteTau} implies \eqref{DeltaHyperK}   }
\label{sect:TauToSigma}

 For any $t\geq M$, we denote by $\mathcal{G}_t$ the $\sigma$-algebra generated by the Poisson r.v.'s $N_i$, $i\leq t$, and by the random walks starting from level $i\leq t$. In particular,  $A^{\ast}_{n,M}(t)$ is $\mathcal{G}_t$-measurable.
 
According to  Proposition~\ref{prop:FiniteTau}, we know that the random time $\tau_m$ is finite for each $m\geq 1$.  Let us now extract a subsequence $(\sigma_m)_{m\geq 1}$ from the sequence $(\tau_m)_{m\geq 1}$. First, we let $\sigma_1=\tau_1$. As a by product of the proof of Theorem~\ref{th:stabilize}, we know that the aggregate which is built by using only particles from levels $R$, thus $R+1$ and so on, does not intersect  the horizontal line $\Z\times\{0\}$, with probability $c(n,R)>0$ for $R$ large enough. This event only concerns randomness above level $R$. Now, given $R>0$ so that $c=c(n,R)>0$, let  $\mbox{Gap}_1$ be the following event:  

\begin{itemize}
\item $N_t=0$ for any $\sigma_1 < t < \sigma_1+\zeta+R$;
\item the aggregate built using only particles from levels $\sigma_1+\zeta+R$, thus $\sigma_1+\zeta+R+1$ and so on, does not intersect the horizontal line $\Z\times\{\sigma_1+\zeta\}$.
\end{itemize}
Conditional on $\mathcal{G}_{\sigma_1}$, the probability of $\mbox{Gap}_1$ is larger than:
\[
c' = c'(n,R,\zeta) = \PP\left( \mbox{$N_t=0$ for any $\sigma_1 < t < \sigma_1+\zeta+R$} \right) \times c \, > \, 0.
\]
Although the event $\mbox{Gap}_1$ depends on $\sigma_1$, this is not the case for the lower bound $c'$.

On the event $\mbox{Gap}_1$, the aggregate $A^{\ast}_{n,M}(\sigma_1)$ remains below the line $\Z\times\{\sigma_1+\zeta\}$ and then cannot help particles coming from levels $t\geq \sigma_1+\zeta+R$ to merge both aggregates, namely $A^{\ast}_{n,M}(\sigma_1)$ and the one which is built with particles from levels $t\geq \sigma_1+\zeta+R$. Hence, the aggregates $A^{\ast}_{n,M}(t)$, for any $t\geq \sigma_1$, does not intersect the horizontal line $\Z\times\{\sigma_1+\zeta\}$. This implies that the set $\Delta(t)$ is stuck below $\Z\times\{\sigma_1+\zeta\}$ for any $t\geq M$. To sum up, conditional on the event $\mbox{Gap}_1$, we have
\[
\left( {\bigcup_{t\geq M}} \Delta(t) \right) \cap \mathcal{H}_{\sigma_1+\zeta+1} = \emptyset,
\]
which implies \eqref{DeltaHyperK} since $\sigma_1$ is a.s.\,finite. 

In order to get back independence when $\mbox{Gap}_1$ does not occur, we have to proceed step by step. Let $t>\sigma_1$. Given a realization of $A^{\ast}_{n,M}(t-1)$, we say that the level $t$ \textit{fails} if one of the following properties holds:
\begin{itemize}
\item $N_t\not= 0$ provided that $t<\sigma_1+\zeta+R$;
\item a particle starting from level $t$ touches $\Z\times\{\sigma_1+\zeta\}$ before exiting the current aggregate provided that  $t\geq \sigma_1+\zeta+R$.
\end{itemize}
Thus $\mbox{Gap}_1^c = \left\{\text{there exists $t>\sigma_1$ such that level $t$ fails}    \right\}$.   

Now, we set
\[
\bar{\sigma}_1 = \inf\{ t > \sigma_1 : \, \mbox{level $t$ fails for $\mbox{Gap}_1$} \} ~.
\]
Notice that $\mbox{Gap}_1^c=\{\bar{\sigma}_1<\infty\}$ is $\mathcal{G}_{\bar{\sigma}_1}$ measurable. In particular, $\mbox{Gap}_1^c$ is $\mathcal{G}_{\sigma_2}$ measurable, where $\sigma_2$ is defined as the first random time $\tau_m$ which is larger than $\bar{\sigma}_1$ (such a quantity exists according to Proposition~\ref{prop:FiniteTau}, provided that $\bar{\sigma}_1$ is finite).   Thus we define the event $\mbox{Gap}_2$ in the same  spirit as we did for $\mbox{Gap}_1$ but this time by replacing $\sigma_1$ by $\sigma_2$. Conditional on $\mathcal{G}_{\sigma_2}$, the event $\mbox{Gap}_2$ occurs with probability larger than $c'>0$. Therefore,
\begin{align*}
\PPP{ \mbox{Gap}_{2}^{c} \cap \mbox{Gap}_{1}^{c} \, \mid \, \mathcal{G}_{\sigma_1}} & =  \PPP{ \{ \bar{\sigma}_2 < \infty \} \cap \{ \bar{\sigma}_1 < \infty \} \, \mid \, \mathcal{G}_{\sigma_1}} \\
& =  \EEE{ \PPP{ \bar{\sigma}_2 < \infty \, \mid \, \mathcal{G}_{\sigma_2}} \ind{\bar{\sigma}_1 < \infty} \, \mid \, \mathcal{G}_{\sigma_1}} \\
& \leq  ( 1 - c' ) \PPP{ \bar{\sigma}_1 < \infty \, \mid \, \mathcal{G}_{\sigma_1}} \\
& \leq  ( 1 - c' )^{2} ~.
\end{align*}
Thus we proceed as previously by introducing
\[
\bar{\sigma}_2 = \inf\{ t > \sigma_2 : \, \mbox{level $t$ fails for $\mbox{Gap}_2$} \} ~.
\]
If $\bar{\sigma}_2=\infty$ then $\mbox{Gap}_2$ occurs and ${\bigcup_{t\geq M}} \Delta(t)$ does not overlap $\mathcal{H}_{\sigma_2+\zeta}$, which implies \eqref{DeltaHyperK}. If $\bar{\sigma}_2<\infty$ then we restart the procedure.

By induction, we deduce that ${\bigcap_{m\geq 1}} \mbox{Gap}_{m}^{c}$ has null probability. So, with probability $1$, there exists some (random) number  $m_0$ such that $\mbox{Gap}_{m_0}$ occurs. This means that ${\bigcup_{t\geq M}} \Delta(t)$ does not overlap $\mathcal{H}_{\sigma_{m_0}+\zeta}$,  where the random time $\sigma_{m_0}$ is a.s.\,finite. This concludes the proof of \eqref{DeltaHyperK} follows.\\

\noindent
\textbf{Adaptation of the proof of \eqref{DeltaHyperK} for $A_{n}[\infty]$.} 
{The track is to introduce a family of events in the context of $A_{n}[\infty]$  as we did for  $A_{n}^*[\infty]$ by introducing the events  $\mbox{Gap}_j$, $j\geq 1$. To do it, let $R$ be sufficiently large so that the particles which allow us to construct $A_n[\infty]$ and which are sent from levels $i\geq R$ do not visit the horizontal axis $\Z\times\{1\}$  with probability $\tilde{c}=\tilde{c}(n,R)>0$ (such a quantity $R$ exists according to Theorem~\ref{th:stabilize}). 

Now, conditional on $\mathcal{G}_{\sigma_1}$, let $\mbox{Thine}_1$ be  the following event:
 \begin{itemize}
\item any particle from level $\sigma_1 < t < \sigma_1+\zeta+R$ settles on the segment $\{0\}\times \llbracket \sigma_1+1, \sigma_1+\zeta+R-1\rrbracket$ or goes directly to the site $(0,\sigma_1+\zeta+1)$ and then goes to the right of $(0,\sigma_1+\zeta+1)$ (on the axis $\Z\times \{\sigma_1+\zeta+1\}$)    until exiting the aggregate;
\item the aggregate built from the initial set $\{0\}\times \llbracket \sigma_1+\zeta+1, \sigma_1+\zeta+R-1\rrbracket$ by sending $n$ particles from level $\sigma_1+\zeta+R$, thus $n$ particles from level $\sigma_1+\zeta+R+1$ and so on, does not intersect the horizontal line $\Z\times\{\sigma_1+\zeta+1\}$.
\end{itemize}
Conditional on $\mathcal{G}_{\sigma_1}$, the two above properties are independent. The first one is realized  with positive probability (depending on $n,R,\zeta$). Indeed, the number of steps imposed to the trajectories in the first property only depends on $n,R$ and $\zeta$ (note that no point of the aggregate lies in $\Z\times \{\sigma_1+\zeta+1\}$ before that the particles of level $\sigma_1+1$ are sent). The second one is realized with probability at least $\tilde{c}$ due to the choice of $R$. Thus, the event $\mbox{Thine}_1$ occurs with probability at least $\tilde{c}'(n,R,\zeta)>0$. In a similar way, we can introduce events $\mbox{Thine}_j$, $j\geq 2$, as we did for $\mbox{Gap}_j$. Observe that, on $\mbox{Thine}_j$, a difference in $\Delta (\sigma_j)$ can be relayed by particles sent from levels $\sigma_j+1$ to $\sigma_j+\zeta+R-1$ at level $\sigma_j+\zeta+1$ but not above. Since, on this event, any particle emitted above or at level $\sigma_j+\zeta+R$ settles before it reaches the line $\Z\times\{\sigma_1+\zeta+1\}$, we deduce that on $\mbox{Thine}_j$:
\[
\left( {\bigcup_{t\geq M}} \Delta(t) \right) \cap \mathcal{H}_{\sigma_j+\zeta+2} = \emptyset.
\]
 The rest of the proof works as before.
}

\subsection{Proof of Proposition~\ref{prop:FiniteTau}}
\label{sect:ProofPropFiniteTau}

Let $n,M\geq 1$. Our aim is to determine a threshold $\zeta$ in such a way that $\tau_m(\zeta)$ is finite for any $m\geq 1$. {Let us first focus on the case $m=1$: it will be explained at the end of the section how the general case can be dealt in a similar way.}

For brievety, we write $h_t=h(A^{\ast}_{n,M}(t))$. First remark that the process $(h_t)_{t\geq M}$ is not Markov because $h_{t+1}$ depends on the whole aggregate $A^{\ast}_{n,M}(t)$ and not only on its vertical spreading. To prove Proposition~\ref{prop:FiniteTau}, we first have to state three lemmas. The first one contains the main idea and claims that a \textit{negative drift} holds for $(h_t)_{t\geq M}$ far away from $0$ ({ the $\sigma$-algebra $\mathcal{G}_t$ appearing in this lemma is defined at the beginning of Section \ref{sect:TauToSigma}).}

\begin{lemma}
\label{lem:Negative Drift}
There exists $\zeta_0(n)$ such that, for any $\zeta\geq \zeta_0(n)$, on the event $\{h_t > \zeta\}$, we have a.s.
\[
\E \big[ h_{t+1} - h_t \,\mid\, \mathcal{G}_t \big] \leq -\frac{1}{2} ~.
\]
\end{lemma}

To state the second lemma, let $t\geq M$ be fixed and assume that   $h_t > \zeta$. We denote by  $H_{t,\zeta}$ the following event:
\[
H_{t,\zeta} = \left\{ \begin{array}{c}
\mbox{At least one of the $N_{t+1}$ particles starting from level $t+1$} \\
\mbox{hits the line $\Z\times\{t+\zeta\}$ before exiting the current aggregate}
\end{array} \right\},
\]
with the convention $H_{t,\zeta}=\emptyset$ if $N_{t+1}=0$. The following lemma claims that, conditional on $\mathcal{G}_t$ and \textit{uniformly on $t$}, the event $H_{t,\zeta}$ does not occur with high probability. 
\begin{lemma}
\label{lem:H}
The following limit holds a.s.:
\[
\lim_{\zeta\to\infty} \sup_{t\geq M} \PPP{H_{t,\zeta} \,\mid\, \mathcal{G}_t} = 0 ~.
\]
\end{lemma}
The above result is also one of the main ingredients to derive Lemma~\ref{lem:Negative Drift}. The next one comes from \cite{CSST} and provides finiteness (and also tail decay but we omit this part here) for the hitting time to $0$ for a discrete-time, non-negative valued process $\{Y_t : t \geq 0\}$ which is not necessarily Markov. Only supermartingale structure and moment conditions for increments are assumed. To state it, we denote by $\nu^Y$ the first hitting time to $0$, {\it i.e.} 
\[
\nu^Y=\inf\{t \geq 1 : Y_t = 0\} ~.
\]
\begin{lemma}{\textsc{\cite[Theorem 5.2]{CSST}}. } 
\label{lem:CSST}
Let $\{Y_t : t \geq 0\}$ be a $\{{\cal G}_t : t \geq 0\}$ discrete-time adapted stochastic process taking values in $\R_+$. Suppose that there exist constants $C_0,C_1>0$ such that, for any $t \geq 0$ and a.s.\,on the event $\{Y_t>0\}$, we have:
\begin{enumerate}[(i)]
\item $\EEE{(Y_{t+1}-Y_t) \mid {\cal G}_t} \leq 0$;
\item $\EEE{ (Y_{t+1} - Y_t)^2 \mid {\cal  G}_t } \geq C_0$;
\item $\EEE{ |Y_{t+1} - Y_t|^3 \mid {\cal  G}_t} \leq C_1$.
\end{enumerate}
Then $\nu^Y<\infty$ a.s.. 
\end{lemma}

\begin{prooft}{Proposition~\ref{prop:FiniteTau}}
Let $\zeta\geq \zeta_0(n)$, where $\zeta_0(n)$ is as in Lemma~\ref{lem:Negative Drift}. To prove that  $\tau_1(\zeta)$ is finite a.s., it is sufficient to apply Lemma~\ref{lem:CSST} to the process $\{Y_t : t\geq 0\}$, where 
\[
Y_t = h_{t} \ind{h_t > \zeta} ~,
\]
for any $t\geq 0$. We check below the three assumptions of Lemma~\ref{lem:CSST}. First, on the event $\{Y_t>0\}$, we notice that $h_t>\zeta$ and $Y_t=h_t$. Hence $Y_{t+1}-Y_t\leq h_{t+1}-h_t$ and Assumption (i) immediately follows from Lemma~\ref{lem:Negative Drift}. 

To prove (ii), we consider two cases. First assume that $h_t > \zeta+1$. In this case, we have $h_{t+1}>\zeta$ and $Y_{t+1}=h_{t+1}$. Since $h_{t+1}-h_t=-1$ on $\{h_t>\zeta\}\cap H_{t,\zeta}^{c}$, we get
\[
\EEE{ (Y_{t+1} - Y_t)^2 \, \mid \, {\cal  G}_t } \geq \EEE{ (h_{t+1} - h_t)^2 \ind{H_{t,\zeta}^{c}} \, \mid \, {\cal  G}_t} = \PP\left( H_{t,\zeta}^{c} \,\mid\, \mathcal{G}_t \right),
\]
which is larger than $C_0=1/2$ for $\zeta$ large enough thanks to Lemma~\ref{lem:H}. Now, if $h_t=\zeta+1$, we have $h_{t+1}=\zeta$ on $\{h_t>\zeta\}\cap H_{t,\zeta}^{c}$ and we conclude similarly. This proves (ii). 

Assumption (iii) is easy to check since $-1\leq h_{t+1}-h_t\leq N_{t+1}$. 
\end{prooft}

\begin{prooft}{Lemma~\ref{lem:Negative Drift}} Let us work conditional on $\mathcal{G}_t$ and assume that $h_t>\zeta$. If the event $H_{t,\zeta}^{c}$ occurs, then the highest ordinate which is reached by the aggregate $A^{\ast}_{n,M}(\cdot)$ does not change between times $t$ and $t+1$ since $h_t > \zeta$, and thus $h_{t+1}-h_t=-1$. If not, we bound $h_{t+1}-h_t$ by the Poisson random variable $N_{t+1}$. So, using the Cauchy-Schwarz inequality, we get 
\begin{align*}
\EEE{ h_{t+1} - h_t \,\mid\, \mathcal{G}_t} & \leq  - \PPP{ H_{t,\zeta}^{c} \,\mid\, \mathcal{G}_t} + \EEE{ N_{t+1} \ind{H_{t,\zeta}} \,\mid\, \mathcal{G}_t} \\
 & \leq  - 1 + \PPP{H_{t,\zeta} \,\mid\, \mathcal{G}_t} + C \, \PPP{ H_{t,\zeta} \,\mid\, \mathcal{G}_t} ^{1/2}
\end{align*}
where $C = \E \big[ N_{t+1}^{2} \mid \mathcal{G}_t \big]^{1/2} = \E \big[ N_{t+1}^{2} \big]^{1/2} = \E \big[ N_{0}^{2} \big]^{1/2}$ is finite and independent of $t$. According to Lemma~{\ref{lem:H}},  we can choose $\zeta$ large enough so that, for any $t\geq M$, the expectation $\E[ h_{t+1} - h_t \mid \mathcal{G}_t]$ is smaller than $-1/2$.  
\end{prooft}

{
\begin{prooft}{ Lemma~\ref{lem:H}}  We do not give all the details because it relies on an adaptation of the proof of  Theorem~\ref{th:stabilize}. In particular, we will introduce a notion of good and bad annuli. 

 Let $\zeta>0$. Because $N_{t+1}$ is Poisson distributed with parameter $n$ for any $t\geq M$, we have
\[
\sup_{t\geq M}\PPP{ N_{t+1} > \zeta^{\delta}} \leq \sup_{t\geq M} \zeta^{-\delta} \EEE{N_{t+1}} = \zeta^{-\delta} n ~.
\]
So it is sufficient to prove that, a.s.,
\begin{equation*}
\label{ProbaHcapN}
\lim_{\zeta\to\infty} \sup_{t\geq M} \PPP{ H_{t,\zeta} \cap \left\{ N_{t+1} \leq \zeta^{\delta} \right\} \,\mid\, \mathcal{G}_t} \, = \, 0 ~.
\end{equation*}
On the event $H_{t,\zeta} \cap \left\{ N_{t+1}\leq\zeta^{\delta} \right\}$, there exists $1\leq i\leq \zeta^{\delta}$ such that the $i$-th particle starting from level $t+1$ reaches ordinate $t+\zeta$ before exiting the current aggregate.
To deal with this event, we have to control the sizes of the current aggregates. Because they are a.s.\,included in $A^{\ast}_{n,M}(\infty)$, we introduce the random variable
\[
X(t,\zeta,n) = \# \left( A^{\ast}_{n,M}(\infty) \cap \left\{(x,y)\in\Z^2 : t+1 \leq y\leq t+\zeta \right\} \right) ~.
\]
Since $A^{\ast}_{n,M}(\infty)$ is included in $A^{\ast}_{n}[\infty]$, we have 
\[
\EEE{ X(t,\zeta,n)} \leq \EEE{ \# \left( A^{\ast}_{n}[\infty] \cap \{(x,y)\in\Z^2 : t+1 \leq y\leq t+\zeta \} \right)} = n \zeta,
\]
where the last equality comes from Propositions~\ref{prop:transfo} (i) and~\ref{prop:ExpWidth}. Taking $0<\delta'<1/2$, we have
\[
\sup_{t\geq M}\PPP{ X(t,\zeta,n) > \zeta^{1+\delta'}} \leq \frac{n \zeta}{\zeta^{1+\delta'}} = n \zeta^{-\delta'} ~.
\]
Hence, it suffices to show that the following conditional probability tends to $0$ as $\zeta$ goes to $\infty$, uniformly on $t$:
\begin{equation}
\label{ProbaOneParticle}
\PP\left( \left.\bigcup_{i\leq\zeta^{\delta}} \left\lbrace \begin{array}{c}
\mbox{The $i$-th particle starting from level $t+1$}\\
\mbox{reaches ordinate $t+\zeta$}\\
\mbox{before exiting the current aggregate}
\end{array} \right\rbrace \cap \left\{ X(t,\zeta,n) \leq \zeta^{1+\delta'} \right\} \right| \mathcal{G}_t \right) ~.
\end{equation}

Let $B_0 = A^{\ast}_{n,M}(t)\cap\{(x,y)\in\Z^2 : y\leq t+1\}$ and let us set for any integer $\ell$
\[
B_{\ell} = B_0 {+} B(0,\ell r),
\]
where $r=\zeta^{1/2}$. Thus $S_{\ell+1} = B_{\ell+1}\!\setminus\! B_{\ell}$ is an annulus with width $r$. Hence, a particle starting from the source $(0,t+1)$ has to cross at least $\ell_\text{max}=\lfloor \zeta/r\rfloor$ such annuli to reach the horizontal line with ordinate $t+\zeta$.

An annulus $S_{\ell}$ is said \textit{good} if
\[
\# \left( A^{\ast}_{n,M}(t) \cap S_{\ell} \right) \leq \eta_0 r^{2} ~.
\] 
Otherwise, we say that it is \textit{bad}. Since the aggregate $A^{\ast}_{n,M}(t)$ is $\mathcal{G}_t$-measurable, we know which annuli are good or not conditional on $\mathcal{G}_t$. Let $N^{\mbox{\tiny{good}}}$ and $N^{\mbox{\tiny{bad}}}$ be the numbers of good and bad annuli. We know that
\[
N^{\mbox{\tiny{good}}} + N^{\mbox{\tiny{bad}}} = \ell_{\max} = \lfloor \zeta/r\rfloor.
\]
On the event $\left\{X(t,\zeta,n)\leq\zeta^{1+\delta'}\right\}$, we have 
\[
N^{\mbox{\tiny{bad}}} \times \eta_0 r^{2} \leq \# \left( A^{\ast}_{n,M}(\infty) \cap \left\{(x,y)\in\Z^2 : t+1 \leq y\leq \zeta \right\} \right) \leq \zeta^{1+\delta'} ~.
\]
We deduce from these inequalities that
\[
N^{\mbox{\tiny{good}}} \geq \left\lfloor \frac{\zeta}{r} \right\rfloor - \frac{\zeta^{1+\delta'}}{\eta_0 r^{2}} ~.
\]
The choice of the parameter $r=\zeta^{1/2}$ (and $\delta'<1/2$) implies that $N^{\mbox{\tiny{good}}}$ is larger than $c \zeta^{1/2}$, for some universal constant $c>0$.

The sequel is very close to the proof of Theorem~\ref{th:stabilize}. Conditional on $\mathcal{G}_t$, we denote by $X_1$ the number of good annuli crossed by particle $1$. Thus, given $X_1,\ldots,X_i$ with $1\leq i < \zeta^{\delta}$, we denote by $X_{i+1}$ the number of good annuli crossed by particle $i+1$, \textit{but not} already crossed by particles $1,\ldots,i$. On the one hand, the Crossing Lemma (Lemma~\ref{le:crossinglemma}) claims that, conditional on $\mathcal{G}_t$ and $X_1,\ldots,X_i$, the r.v. $X_{i+1}$ is stochastically dominated by a geometric distribution with parameter $1-\eta$, where $\eta$ is as in Lemma~\ref{le:crossinglemma}. On the other hand, the fact that at least one particle starting from level $t+1$ reaches the ordinate $t+\zeta$ means that $\sum_{i\leq \zeta^\delta} X_i$ is larger than $N^{\mbox{\tiny{good}}}$. In particular, there exists $i\leq \zeta^\delta$ such that $X_i\geq c\zeta^{1/2-\delta}$. Therefore, the probability appearing in \eqref{ProbaOneParticle} is bounded by $\zeta^{\delta} (1 - \eta )^{c \zeta^{1/2-\delta}}$, which tends to $0$ uniformly on $t$. This concludes the proof of Lemma~\ref{lem:H}.
\end{prooft}
}

{Let us justify that the above proof of the a.s. finiteness of $\tau_1(\zeta)$ also works for $\tau_m(\zeta)$, for any integer $m$. On the one hand the null limit of $\PPP{H_{t,\zeta} \,\mid\, \mathcal{G}_t}$, as $\zeta\to\infty$, is uniform on $t$. In particular, the width of the current aggregate $A^{\ast}_{n,M}(t)$ between levels $t+1$ and $t+\zeta$ is controled uniformly on $t$ (via the r.v. $X(t,\zeta,n)$, see the proof of Lemma \ref{lem:H}). On the other hand, the negative drift (Lemma \ref{lem:Negative Drift}) only uses the vertical spreading $h_t$.}

\section{Mixing property for $A_n[\infty]$, $A^\ast_n[\infty]$ and $A^\dag_n[\infty]$}
\label{sec:mixing}

Combining the stabilization results obtained in the previous sections (Theorems~\ref{th:stabilize} and~\ref{cor:StabWRTO}), we get a mixing property for aggregates $A_n[\infty]$, $A^\ast_n[\infty]$ and $A^\dag_n[\infty]$.  An alternative (longer) proof of Proposition~\ref{prop:mixing} could be provided, without using the latter theorems, by introducing once again  good and bad annuli as described in Section~\ref{sec:hightpart}. 

\begin{proposition}
\label{prop:mixing}
Let $n\geq 1$.
\begin{enumerate}[(i)]
\item The distribution of $A_n[\infty]$ is mixing w.r.t.\,vertical translations, {\it i.e.} 
\[\lim_{|k|\rightarrow\infty}\PPP{A_n[\infty]\in \mathcal{A}, A_n[\infty]\in \tau_k\mathcal{B}} = \PPP{A_n[\infty]\in \mathcal{A}}\PPP{A_n[\infty]\in \mathcal{B}}\]   for any events $\mathcal{A}, \mathcal{B}$.
\item The same holds for $A_n^\ast[\infty]$ and $A_n^\dag[\infty]$.
\end{enumerate}
\end{proposition}

The next result is a direct consequence of the above mixing property and Proposition~\ref{prop:DifferentLaw}. It claims that the infinite aggregates $A^*_n[\infty]$ and $A_n^\dag[\infty]$ have a.s.\,infinitely many finite connected components. Here, we say that $C\subset\Z^2$ is {\it connected} if for any $z,z'\in C$ there exist $z_0=z,z_1, \dots,z_n=z'\in C$ with $|z_j-z_{j-1}|_1=1$ for every $j\in\llbracket 1,n\rrbracket$. Corollary~\ref{cor:finitecc} will be used to define the IDLA forest in Section~\ref{sec:forest}. 

\begin{corollary}
\label{cor:finitecc}
\begin{enumerate}[(i)]
\item Let $n\geq 1$. With probability 1, for any integer $M$ there are (infinitely many) levels $i\geq M$ and $j\leq -M$ such that the aggregate\, $A_n^*[\infty]$ does not intersect the axes $\Z\times\{i\}$ and $\Z\times\{j\}$. In particular, a.s.\,$A_n^*[\infty]$ has only finite connected components included in disjoint strips. 
\item The same holds for $A_n^\dag[\infty]$.
\end{enumerate}
\end{corollary}

\begin{prooft}{Corollary~\ref{cor:finitecc}} Notice that from Propositions~\ref{prop:transfo} and~\ref{prop:DifferentLaw}, we have, for any integer $m\in\Z$,
\[
\PP\big( A_n^*[\infty] \cap (\Z\times \{m\}) = \emptyset \big) = \PP\big( A_n^*[\infty] \cap (\Z\times \{0\}) = \emptyset \big) > 0 ~.
\]
Because the infinite aggregate $A_n^*[\infty]$ is mixing w.r.t.\,vertical translations (Proposition~\ref{prop:mixing}) and thus ergodic, we then get that, with probability $1$, for any $M$, there exist levels $i>M$ and $j<-M$ such that $A_n^*[\infty]$ avoids the axes $\Z\times\{i\}$ and $\Z\times\{j\}$. The same holds for $A_n^\dag[\infty]$ because $A_n^*[\infty]$ and $A_n^\dag[\infty]$ are equally distributed.
\end{prooft}

\begin{prooft}{Proposition~\ref{prop:mixing}}
We only prove (i). The proof is similar for $A^\ast_n[\infty]$ and so for $A^\dag_n[\infty]$. 

It is enough to prove that
\begin{equation*}
\label{Mixing}
\lim_{k\to \infty} \PPP{A_n[\infty]\cap \left(C_1\cup\tau_k C_2\right)=\emptyset}=\PPP{A_n[\infty]\cap C_1=\emptyset}\PPP{A_n[\infty]\cap C_2=\emptyset}
\end{equation*}
for all compact sets $C_1,C_2\subset \Z^2$ (see {\it e.g.\,}\cite[Theorem 9.3.2]{SW}). The case where $k\to-\infty$ works as well.

Let $\varepsilon>0$. By Theorem~\ref{th:stabilize} (with $\alpha=2$), we can choose $M$ large enough so that $C_1\subset \Z_M$ and
\[
\PP\big( A_n[\infty] \cap \Z_M = A_n[M^2] \cap \Z_M \big) \geq 1 - \varepsilon/2 ~.
\]
Furthermore, Theorem~\ref{cor:StabWRTO} asserts that
\[
\lim_{M\to\infty} \lim_{t\to\infty} \PPP{ A_{n}[\infty] \cap \mathcal{H}_t = B_{n,M}(\infty) \cap \mathcal{H}_t} = 1,
\]
where $B_{n,M}(\infty)$ denotes the aggregate, coupled with $A_{n}[\infty]$ and included in $A_n[\infty]$,  which is built by   sending particles from levels $M+1$, $M+2$ and so on (see Section~\ref{sec:stabOrigine}). Taking $M$ sufficiently large, we can choose $t\geq M^2$ large enough so that
\[
\PPP{ A_{n}[\infty] \cap \mathcal{H}_t = B_{n,M^2}(\infty) \cap \mathcal{H}_t }\geq 1 - \varepsilon/2 ~.
\]
Now, let $F_{M,t}$ be the following event:
\[F_{M,t} = \left\{A_n[\infty] \cap \Z_M = A_n[M^2] \cap \Z_M \text{ and } A_{n}[\infty] \cap \mathcal{H}_t = B_{n,M^2}(\infty) \cap \mathcal{H}_t  \right\}.\] Let $k$ be such that $\tau_k C_2$ is included in $\mathcal{H}_t$. Since $\PPP{F_{M,t}}\geq 1-\varepsilon$, we have 
\begin{multline*}
\PP\left( A_n[\infty] \cap C_1 = \emptyset \, , \, A_n[\infty] \cap \tau_k C_2 = \emptyset \right) \\
\begin{split}
& = \PP\left( \left\{A_n[M^2] \cap C_1 = \emptyset \, , \, B_{n,M^2}(\infty) \cap \tau_k C_2 = \emptyset \right\}\cap   F_{M,t} \right) \pm \varepsilon\\
& = \PP\left( A_n[M^2] \cap C_1 = \emptyset, B_{n,M^2}(\infty) \cap \tau_k C_2 = \emptyset \right) \pm 2\varepsilon\\
& = \PP\left( A_n[M^2] \cap C_1 = \emptyset \right) \, \PP\left( B_{n,M^2}(\infty) \cap \tau_k C_2 = \emptyset \right) \pm 2\varepsilon,\\
\end{split}
\end{multline*}
where the last line comes from the fact that  $A_n[M^2]$ and $B_{n,M^2}(\infty)$ are based on  independent particles. In the above equation, the notation $p=p'\pm \varepsilon$ means that $|p-p'|\leq \varepsilon$ for any $p,p'\in \R$.  Using once again that $\PPP{F_{M,t}}\geq 1-\varepsilon$,  we have 
\begin{multline*}
\PP\left( A_n[M^2] \cap C_1 = \emptyset \right) \, \PP\left( B_{n,M^2}(\infty) \cap \tau_k C_2 = \emptyset \right) \\
\begin{split}
& =  \PP\left(\left\{ A_n[M^2] \cap C_1 = \emptyset\right\}\cap  F_{M,t} \right) \, \PP \left(\left\{ B_{n,M^2}(\infty) \cap \tau_k C_2 = \emptyset \right\} \cap  F_{M,t} \right) \pm 2\varepsilon \\
& =  \PP\left( A_n[\infty] \cap C_1 = \emptyset \right) \, \PP\left( A_n[\infty] \cap \tau_k C_2 = \emptyset \right) \pm 4\varepsilon \\
& =  \PP\left( A_n[\infty] \cap C_1 = \emptyset \right) \, \PP\left( A_n[\infty] \cap C_2 = \emptyset \right) \pm 4\varepsilon
\end{split}
\end{multline*}
since the distribution of $A_n[\infty]$ is invariant w.r.t.\,vertical translations. This concludes the proof of Proposition~\ref{prop:mixing}.
\end{prooft}

\section{Shape theorems}\label{sec:shapeTh}

Proposition~\ref{prop:ExpWidth} claims that, in expectation, the infinite aggregates $A_n[\infty]$, $A_n^*[\infty]$ and $A_n^\dag[\infty]$ can be approximated with the rectangle $R_{n/2}$, where
\begin{equation*} \label{eq:defrectangle}  R_r=\left\llbracket  - r, r\right\rrbracket \times \Z.  
 \end{equation*}
 for any $r>0$.  The following result ensures that the fluctuations of these aggregates, when they are restricted to a strip, are at most logarithmic.

  {
 \begin{theorem}
  \label{th:ShapeTh}
 There exists $A>0$ such that, for any $\alpha>0$, a.s.\,there exists $N\geq 1$ such that for any $n\geq N$, 
  \[R_{n/2 - A\log(n)}\cap \Z_{n^\alpha} \subset   A_n[\infty]\cap \Z_{n^\alpha}\subset R_{n/2 + A\log(n)}\cap \Z_{n^\alpha}.\]
  \end{theorem}
  }
  
{
When the number of particles which is sent from each site of the vertical axis $I(\infty)$ is Poisson or based on Poisson clocks, we obtain the following shape theorem.
\begin{theorem}
\label{th:ShapeThPoisson}
For any $0 < \varepsilon < 1/2$ and for any $\alpha > 0$, a.s.\,there exists $N\geq 1$ such that for any $n\geq N$, 
\[
R_{n/2 - n^{1/2+\varepsilon}} \cap \Z_{n^\alpha} \subset A_n^*[\infty] \cap \Z_{n^\alpha} \subset R_{n/2 + n^{1/2+\varepsilon}} \cap \Z_{n^\alpha} ~.
\]
The same holds for $A_n^\dag[\infty]$.
\end{theorem}   
Theorem \ref{th:ShapeThPoisson} is less precise than Theorem \ref{th:ShapeTh} since the fluctuations are polynomial  whereas they are logarithmic for the aggregate $A_n[\infty]$. Roughly, this can be explained by the fact that a Poisson random variable of parameter $n$ can have values lower than $n-n^{\beta}$, $0<\beta<1/2$, with positive probability (independent of $n$). In Remark \ref{rk:shapetheorempoisson}, we explain why the approach which used to prove Theorem \ref{th:ShapeTh} cannot directly applied to derive Theorem \ref{th:ShapeThPoisson}.
  }

As an illustration of Theorem~\ref{th:ShapeTh}, {Figure~\ref{fig:shapetheorem}} shows that $A_n[\infty]$, when it is restricted to a strip, is very close to a rectangle. { Notice that, for any integer $n$, there exists a.s.$\,$a random level $i$ such that $(n,i)$ belongs to the aggregate. So, the intersection with the horizontal strip $\Z_M$  is unavoidable in Theorem~\ref{th:ShapeTh}.} 

 Theorem~\ref{th:ShapeTh} is classical in the sense that many results dealing with the shape of aggregates were established for various IDLA models, see {\it e.g.} \cite{AG13long, AG13short,BDCKL,LS}, but these results cannot  be applied in our context. The proof of Theorem~\ref{th:ShapeTh} relies on an adaptation of \cite{AG13long,AG13short}. {This adaptation is not direct at some places because} our random aggregates are based on an infinite number of particles and are \textit{not isotropic}. The counterpart is that our models are invariant w.r.t\, vertical translations and satisfy a mass transport principle (see Proposition~\ref{prop:ExpWidth}). These properties will be extensively used in our proof (the main differences w.r.t.\,\cite{AG13long,AG13short} will be stressed). 
 
 Notice that Theorem~\ref{th:ShapeTh} provides upper bounds for the fluctuations of our random aggregates in terms of Hausdorff distance. In the same spirit as \cite{AG13short}, we could show that the fluctuations are \textit{at least} logarithmic. 

With a slight abuse of notation, lower and upper integer parts are omitted along this section.

\subsection{The lower bound for $A_n[\infty]$} 
\label{subsec:lowerbound}

In this section, we prove the following result (referred to as the lower bound for the shape theorem {of $A_n[\infty]$}): there exists $A>0$ such that, for any $\alpha>0$, a.s.\,there exists $N\geq 1$ such that for any $n\geq N$, 
  \begin{equation}\label{eq:lowerbound}
  R_{n/2 - A\log(n)}\cap \Z_{n^\alpha} \subset   A_n[\infty]\cap \Z_{n^\alpha}.
  \end{equation}
To do it, we mainly follows the strategy developed in \cite{AG13long,AG13short} and, then we will make several references to these articles. In order to make clearer these references, we adopt -- only for Section~\ref{sec:shapeTh} -- their notation $A(\eta)$ to denote the aggregate generated by an initial configuration $\eta$. The notation $\eta=n\ind{C}$, where $n\geq 1$ and $C\subset \Z^2$, means that we send $n$ particles from each site of $C$. In particular, we have $\A{M}=A_n[M]$ and  $\A{\infty}=A_n[\infty]$.
 
 For each integer $k$, we let 
  \[S_k=\left(R_{(k+1)\log(n)}\setminus R_{k\log(n)}\right)\cap \Z_{n^\alpha}. \] The set $S_k$ is referred to as a \textit{shell}. Now, for each $z\in \partial R_{k\log(n)}$, with \[\partial R_{k\log(n)}=\{-\lfloor k\log(n)\rfloor,\lfloor k\log(n)\rfloor\}\times \Z,\] we define the so-called \textit{tile} and \textit{cell} centered at $z$ as
\begin{equation*}\label{eq:defTileAndCell}
\tau(z) = B(z,\log(n)/2)\cap \partial R_{k\log(n)} \qquad\mbox{and}\qquad C(z) = B(z,\log(n))\cap R_{k\log(n)}^c.
\end{equation*}
 Notice that  \[S_k\subset \underset{z\in \partial_{k,n}}{\bigcup} C(z),\]
where $\partial_{k,n}=\partial R_{k\log(n)}\cap \Z_{n^\alpha}$. Given $\eta$ and $B\subset R_{k\log(n)}$,  the number of particles (resp. random walks), with initial configuration $\eta$, hitting $B$ before or when they exit $R_{k\log(n)}$, are  denoted by $W_{k\log(n)}(\eta, B)$ (resp. $M_{k\log(n)}(\eta, B)$). With a slight abuse of notation, when $\eta=\ind{C}$, we simply write $W_{k\log(n)}(C, B)$ (resp. $M_{k\log(n)}(C, B)$) instead of $W_{k\log(n)}(\ind{C}, B)$ (resp. $M_{k\log(n)}(\ind{C}, B)$). Recall that we use the word \textit{particle} for a random walk which is stopped when exiting the aggregate and adding site; thus trajectories of particles depend on the aggregate while random walks do not.

We say that a set $B$ is {\it not covered} if  $B \not\subset \A{\infty}$. According to the Borel-Cantelli lemma, it is sufficient to prove that there exists $A$ such that, for any $L>0$, $n\geq 1$ and $k\leq \frac{n}{2\log(n)}-A$, we have:
\[\PPP{S_k \text{ is not covered}}\leq c n^{-L}.\]  

As in \cite{AG13short,Lawler95}, it is useful to stop the particles when they reach $\partial R_{k\log(n)}$. {The strategy can be divided into two steps. Roughly, it consists in, first, showing that each tile $\tau$ of $\partial R_{k\log(n)}$ is likely to capture many particles and then arguing that, if many particles exit $\partial R_{k\log(n)}$ from $\tau$, then they are likely to cover the corresponding cell $C$.} To do it, let $A$ be fixed and let $k\leq \frac{n}{2\log(n)}-A$. We write 
\begin{equation}
\label{eq:noncoveredstrip}
\PPP{S_k \text{ is not covered}} \leq p(n,k)+q(n,k),
\end{equation}
where
\[p(n,k) = \PPP{\exists \tau\in \mathcal{T}_{k\log(n)}, W_{k\log(n)}(n\ind{I(\infty)},\tau)\leq \frac{1}{3}\mu(\tau) }\]
and 
\[q(n,k) = \PPP{\forall \tau\in \mathcal{T}_{k\log(n)} , W_{k\log(n)}(n\ind{I(\infty)},\tau)\geq \frac{1}{3}\mu(\tau) \text{ and $S_k$ is not covered} }.\] The term $\mu(\tau)$ appearing in the above equations will be defined later in \eqref{eq:defmu}, and the set $\mathcal{T}_{k\log(n)}$ denotes the family of tiles, {\it i.e.} \[\mathcal{T}_{k\log(n)} = \{\tau(z): z\in \partial_{k,n}\}.\]

We prove below that $p(n,k)$ and $q(n,k)$ are smaller than any power of $n^{-1}$ when $A$ is large enough. 
 
\subsubsection{Upper bound for $p(n,k)$}
 Since $\#\mathcal{T}_{k\log(n)}\leq 4n^\alpha+2$, it is sufficient to prove that

\begin{equation*}\label{eq:LB1}
\PPP{W_{k\log(n)}(n\ind{I(\infty)},\tau)\leq \frac{1}{3}\mu(\tau)}
\end{equation*}  is lower than any power of $n^{-1}$, for any tile $\tau \in\mathcal{T}_{k\log(n)}$. To do it, we will apply an analog of Lemma 2.4 of \cite{AG13short} that deals with series of Bernoulli random variables instead of sums of Bernoulli variables. It is obtained by straightforward modifications of the proof of \cite[Lemma 2.4]{AG13short} and is  stated as follows.
 
\begin{lemma} \label{Le:AGshortL2.4}
Suppose that a sequence of random variables $\{W_n, M_n,L_n,\widetilde{M}_n; n\geq 0\}$ and a sequence of real numbers $(c_n)_{n\geq 0}$ satisfy for any $n\geq 0$:
\begin{equation*}
W_n+L_n+c_n\geq \widetilde{M}_n\qquad\mbox{and}\qquad \widetilde{M}_n\overset{\text{law}}{=}M_n.
\end{equation*}
Assume that $W_n$ and $L_n$ are independent and that $L_n$ and $M_n$ are both series of independent Bernoulli random variables with finite first moment. Assume also that 
\begin{enumerate}
\item[\bf (H1)]\label{H1} the Bernoulli variables $\{Y_1^{(n)}, Y_{2}^{(n)}, \dots\}$ whose series is $L_n$ satisfy for some $\kappa>1$:\[\underset{n}{\sup}~\underset{i}{\sup}\,\EEE{Y_i^{(n)}}<\dfrac{\kappa -1}{\kappa};\]
\item[\bf (H2)]\label{H2} $\mu_n=\EEE{M_n-L_n}\geq 0$.
\end{enumerate}
Then, for any $n\geq 0$ and $\xi_n\in\R$, we have for any $\lambda\geq 0$,
\begin{equation*}
\PPP{W_n<\xi_n}\leq \exp\left(-\lambda\left(\mu_n-\xi_n-c_n\right)+\frac{\lambda^2}{2}\left(\mu_n+\kappa\sum_{i=1}^{\infty}\EEE{Y_i^{(n)}}^2\right)\right).
\end{equation*}
\end{lemma} 

The desired upper bound for~\ref{eq:LB1} will be obtained thanks to Lemma~\ref{Le:AGshortL2.4}. Hence, we must check that hypotheses of Lemma~\ref{Le:AGshortL2.4} are satisfied. First of all, note that, following the strategy initiated in \cite{LBG92}, similar arguments to \cite{AG13short} show that
\begin{equation}\label{eq:stochdom}
W_{k\log(n)}(n\ind{I(\infty)}, \tau)+M_{k\log(n)}(R_{k\log(n)},\tau) \geq \widetilde{M}_{k\log(n)}(n\ind{I(\infty)}, \tau),
\end{equation}
with
\[\widetilde{M}_{k\log(n)}(n\ind{I(\infty)}, \tau)=W_{k\log(n)}(n\ind{I(\infty)}, \tau)+M_{k\log(n)}(A_{k\log(n)}(n\ind{I(\infty)}),\tau),\]
{where $A_r(\eta)$ denotes the aggregate produced by particles with initial configuration $\eta$ stopped when they leave $R_r$.}

For $z\in \partial R_{k\log(n)}$, let $\mathcal{Z}=\mathcal{Z}(z,b,n)$ be the set
\[\mathcal{Z}=\{z'\in R_{k\log(n)}:\, d(z',\tau(z))\leq b\log (n)\}.\]
The set $\mathcal{Z}$ will be useful to check that {\bf (H1)} holds, for some $b>0$. Set $c_n=\vert \mathcal{Z}\vert\leq c(b\log(n))^2$,  so that \eqref{eq:stochdom} leads to:   
\begin{equation*}\label{eq:stochdom2}
W_{k\log(n)}(n\ind{I(\infty)}, \tau)+M_{k\log(n)}(R_{k\log(n)}\setminus\mathcal{Z},\tau)+c_n \geq \widetilde{M}_{k\log(n)}(n\ind{I(\infty)}, \tau).
\end{equation*}

We will apply Lemma~\ref{Le:AGshortL2.4} with $W_n=W_{k\log (n)}(n\ind{I(\infty)}, \tau)$, $M_n=M_{k\log (n)}(n\ind{I(\infty)}, \tau)$, $L_n=M_{k\log(n)}(R_{k\log(n)}\setminus\mathcal{Z},\tau)$ and $\widetilde{M}_n=\widetilde{M}_{k\log (n)}(n\ind{I(\infty)},\tau)$. So, $\mu(\tau)$ in \eqref{eq:noncoveredstrip} must be defined as:
\begin{equation}
\mu(\tau)=\EEE{M_{k\log(n)}(n\ind{I(\infty)}, \tau)}-\EEE{M_{k\log(n)}(R_{k\log(n)}\setminus\mathcal{Z},\tau)}.\label{eq:defmu}
\end{equation}

\paragraph{Verification of {\bf (H1)}}
Assumption {\bf(H1)} is ensured  for some $b>0$ by the following lemma which is the counterpart in our context of \cite[Lemma 5.1]{AG13long}. Its proof is very similar to the one of Lemma 5.1 in \cite{AG13long} and is omitted here.

\begin{lemma}
\label{Le:probhitting}
There exists $\kappa>0$ such that, for any $r$, any $y\in R_r$ and any $x\in \partial R_r\setminus\{y\}$, we have
\[\PP_y(S(H_r) = x) \leq \frac{\kappa}{|x-y|},\]
where $H_r$ denotes the hitting time of $\partial R_{r}$ for the simple random walk $(S(t))_{t\geq 0}$.
\end{lemma} 
 
\paragraph{Lower bound for $\mu(\tau)$ and verification of {\bf (H2)}}
We show that, if $z$ is at a distance at least $A\log(n)$ from $\partial R_{\frac{n}{2}}$, then 
\begin{equation}
\label{eq:LBmutau} \mu(\tau)\geq cA\log(n)^2
\end{equation}
for some $c>0$. Note that this ensures {\bf(H2)}. {In what follows, given a multiset $\eta$ and $r\leq r'$, we denote by $M_{r'}(\eta,\tau)$ the number of random walks starting from $\eta$ that exit $R_{r'}$ in a site of $\tau$. Inequality \eqref{eq:LBmutau}  will be derived from the following lemma. Its proof relies on the invariance w.r.t\, vertical translations of our model and cannot be adapted from \cite{AG13long,AG13short}.

\begin{lemma}\label{Le:eqnegligible} Let $r\leq r'$ and let $\tau \subset\partial R_{r'}$ be finite. Then 
\[\EEE{M_{r'}(R_{r},\tau)}=\frac{2r+1}{2}\#\tau.\]
\end{lemma}
In particular, as a consequence of Lemma \ref{Le:eqnegligible} applied to $R_0=I(\infty)$, we have for any $r'\geq 1$,
\begin{equation}
 \label{eq:expectationbis}
 \EEE{M_{r'}(I(\infty), \tau)} = \frac{\#\tau}{2}.
\end{equation}

\begin{prooft}{Lemma~\ref{Le:eqnegligible}}
With the notation of Lemma~\ref{Le:probhitting}, we have:
\begin{align*}
\EEE{M_{r'}(R_r,\tau)}&=\sum_{y\in \tau} \sum_{z\in R_r}\mathbb{P}_z\left(S(H_{r'})=y\right)=\#\tau \sum_{z\in R_r}\mathbb{P}_z\left(S(H_{r'})=(r',0)\right),
\end{align*}  
where we used invariance w.r.t.\,vertical translations and the symmetry w.r.t.\,$I(\infty)$ in the second equality. Hence,
\begin{align*}
\EEE{M_{r'}(R_r,\tau)}
&=\#\tau \sum_{i=-r}^r\sum_{j\in \Z}\mathbb{P}_{(i,j)}\left(S(H_{r'})=(r',0)\right)\\
&=\#\tau \sum_{i=-r}^r\sum_{j\in \Z}\mathbb{P}_{(i,0)}\left(S(H_{r'})=(r',j)\right)\\
&=\#\tau \sum_{i=-r}^r\mathbb{P}_{(i,0)}\left(S(H_{r'})\in \{r'\}\times \Z\right)\\
\end{align*}
The result then follows since, by symmetry w.r.t.\,$I(\infty)$, $$\mathbb{P}_{(i,0)}\left(S(H_{r'})\in \{r'\}\times \Z\right)+\mathbb{P}_{(-i,0)}\left(S(H_{r'})\in \{r'\}\times \Z\right)=1.$$
\end{prooft}}

{Now, to get~\eqref{eq:LBmutau}, we write
\begin{align*}
\mu(\tau) & \geq \EEE{M_{k\log(n)}(n\ind{I(\infty)}, \tau)-M_{k\log(n)}(R_{k\log(n)},\tau)}\\
& = n\EEE{M_{k\log(n)}\left(I(\infty), \tau\right)} - \EEE{M_{k\log(n)}(R_{k\log(n)},\tau)}\\
\end{align*}
According to \eqref{eq:expectationbis} and Lemma \ref{Le:eqnegligible}, we have
\begin{align*}
 n\EEE{M_{k\log(n)}\left(I(\infty), \tau\right)} - \EEE{M_{k\log(n)}(R_{k\log(n)},\tau)} & = \frac{1}{2}\left(n-2k\log(n)-1\right)\#\tau\\
& \geq cA\log(n)^2.
\end{align*}
This gives \eqref{eq:LBmutau}. 

  }
    { 
 
\paragraph{Second-order estimate}
In order to exploit the upper bound given by Lemma~\ref{Le:AGshortL2.4}, we have to control $\sum_i\EEE{Y_i^{n}}^2$.
Due to the definition of $\mathcal{Z}=\mathcal{Z}(z,b,n)$ and $\tau=\tau(z)$, and according to Lemma~\ref{Le:probhitting}, we have for all $y\in R_{k\log(n)}\setminus \mathcal{Z}$ {(corresponding to index $i$ in the sum defining $L_n$)}:
\begin{align*}
\EEE{Y_i^{n}}=\PP_y(S(H_{k\log(n)})\in \tau) & \leq \#\tau \max_{x\in \tau}\PP_y(S(H_{k\log(n)})=x)%
\leq \#\tau \max_{x\in \tau}\frac{\kappa}{|x-y|}%
\leq
c\frac{\#\tau}{|z-y|}.
\end{align*}
Summing over $y$, it follows that
\begin{align*}
\sum_{y\in R_{k\log(n)}\setminus\mathcal{Z}}\PP_y(S(H_{k\log(n)})\in \tau)^2 
&\leq c \#\tau^2\sum_{y\in R_{k\log(n)}\setminus\mathcal{Z}}\frac{1}{|z-y|^{2}}\\
&\leq  \#\tau^2\left( c+2 \sum_{j=1}^{2k\log(n)}\int_1^\infty \frac{1}{j^{2}+x^2}\operatorname{d} x \right)\nonumber
\leq
c \log(n)^3.
\end{align*}

\paragraph{Application of Lemma~\ref{Le:AGshortL2.4}}
Lemma~\ref{Le:AGshortL2.4} and similar computations as in \cite[Section 3.1.2]{AG13short} imply that
\begin{equation}\label{eq:lowerboundpart1} \PPP{W_{k\log(n)}(n\ind{I(\infty)},\tau)\leq \frac{1}{3}\mu(\tau)}\leq \exp\left( -cA^2\log(n)\right).
\end{equation} 
The last term converges to 0 faster than any power of $n^{-1}$ for $A$ large enough.

\subsubsection{Upper bound for $q(n,k)$} 

We prove below that $q(n,k)$ is bounded by any power of $n^{-1}$.

{Recall that \cite[Lemma 1.3]{AG13short} roughly says that if many particles (w.r.t\,its radius) initially lie in the middle of a ball, then the aggregate they produce is very likely to cover the ball. } By using  \eqref{eq:LBmutau} and \cite[Lemma 1.3]{AG13short}, one has:
\begin{align*}
&\PPP{\forall \tau\in \mathcal{T}_{k\log(n)} , W_{k\log(n)}(n\ind{I(\infty)},\tau)\geq \frac{1}{3}\mu(\tau) \text{ and $S_k$ is not covered} }\\
&\qquad\qquad\leq \PPP{\exists z\in  \partial_{k,n}, C(z) \text{ is not covered }  \left| \forall \tau\in \mathcal{T}_{k\log(n)} , W_{k\log(n)}(n\ind{I(\infty)},\tau)\geq \frac{1}{3}\mu(\tau)\right.}\\
&\qquad\qquad\leq cn^{\alpha}\exp\left( -c \frac{A\log(n)^2}{\log(\log(n))} \right).
\end{align*}
The last term converges to 0 faster than any given power $n^{-1}$. This together with \eqref{eq:noncoveredstrip} and \eqref{eq:lowerboundpart1} concludes the proof of \eqref{eq:lowerbound}.

\subsection{The upper bound for $A_n[\infty]$}
\label{subsec:upperbound}

In this section, we prove the following result (referred to as the upper bound for the shape theorem of  $A_n[\infty]$}):  there exists $A>0$ such that for any $\alpha>0$, there  a.s.\,exists $N\geq 1$ such that for any $n\geq N$, 
  \begin{equation}\label{eq:upperbound}
  \A{\infty}\cap \Z_{n^\alpha}\subset R_{ n/2 + A\log(n)}\cap \Z_{n^\alpha}.
  \end{equation}
  
{ Let $\alpha>0$ be fixed and let $A>0$ (which will be chosen sufficiently large later). One can extract from the proof of Theorem~\ref{th:stabilize} that, for any $L>0$, for any $\gamma>\alpha$, and for $n$ large enough, 
 \begin{equation}
 \label{eq:finiteinfinite}
  \PPP{\A{\infty} \cap  \Z_{n^{\alpha}} \neq A(n\ind{\{0\}\times \llbracket -n^\gamma, n^\gamma\rrbracket}) \cap  \Z_{n^{\alpha}}  } \leq n^{-L}.
  \end{equation} 
The above inequality is crucial to derive the upper-bound because it reduces our problem to an aggregate which is generated by a \textit{finite} number of particles. According to the Borel-Cantelli lemma, it is sufficient to prove that, for $n$ large enough,
\begin{equation}\label{eq:UB3}\PPP{A(n\ind{\{0\}\times \llbracket -n^\gamma, n^\gamma\rrbracket} ) \cap \Z_{n^\alpha}   \not\subset  R_{n/2+A\log(n)}\cap \Z_{n^\alpha}} 
\end{equation} 
is smaller than any power of $n^{-1}$. To do it, we bound \eqref{eq:UB3} by
\begin{equation*}
\PPP{t_n\geq\frac{n}{2}+A\log (n)},
\end{equation*}
where $t_n=\max\{|z(1)|:\,z\in A(n\ind{\{0\}\times \llbracket -n^\gamma, n^\gamma\rrbracket})\}$. Remark that $t_n$ is a.s.\,smaller than $n(2n^\gamma+1)$. 
} Taking the supremum over the point $z\in\Z_{n^\alpha}\cap\{z:\frac{n}{2}+A\log (n)\leq\vert z(1)\vert\leq n(2n^\gamma+1)\}$, it is enough to prove that 
\begin{align*}
\sup_z \PPP{z\in A(n\ind{\{0\}\times \llbracket -n^\gamma, n^\gamma\rrbracket}), |z(1)|=t_n}
\end{align*}
is lower than any power of $n^{-1}$. We do so in the same spirit as \cite{AG13short} and we recall several arguments included in \cite{AG13short} in order to make the paper self-contained. Let  $z\in \Z_{n^\alpha}$, with $\frac{n}{2}+A\log(n)\leq |z(1)|\leq n(2n^{\gamma}+1)$. In what follows, we set {\[h(n)=|z(1)|-\frac{n}{2}.\]} First, we write
\begin{align}
&\PPP{z\in A(n\ind{\{0\}\times \llbracket -n^\gamma, n^\gamma\rrbracket}), |z(1)|=t_n}\nonumber\\
 &\qquad\leq \PPP{\#(A(n\ind{\{0\}\times \llbracket -n^\gamma, n^\gamma\rrbracket})\cap B(z,h(n)))>\beta h^2(n), |z(1)|=t_n}\label{eq:UBeq1} \\
&\qquad\quad+ \PPP{z\in A(n\ind{\{0\}\times \llbracket -n^\gamma, n^\gamma\rrbracket}), \#(A(n\ind{\{0\}\times \llbracket -n^\gamma, n^\gamma\rrbracket})\cap B(z,h(n)))\leq \beta h^2(n)},\label{eq:UBeq2}
\end{align}
{where $\beta$ is chosen as in Lemma~\ref{Le:Jerison} below in order to ensure that \eqref{eq:UBeq2} vanishes quickly. The later lemma is an adaptation of \cite[Lemma A]{JLS12}. It allows us to avoid the use of an analog of the {\it flashing process } introduced in \cite{AG13long,AG13short} to derive {Theorem \ref{th:ShapeTh}}.  {However, such a process remains useful to prove that the fluctuations in the shape theorem are of correct order and could be used in our context.}} 

  \begin{lemma}
\label{Le:Jerison}
Let $\gamma>0$ and let $A(\eta)$ be the aggregate with initial configuration $\eta$, such that $\eta$ has support in $\{0\}\times \llbracket -n^\gamma, n^\gamma\rrbracket$ and $\#\eta<\infty$.  Then there exist positive universal constants $\beta$, $C_0$ and $c$ such that for any real number $m>0$ and all $z\in \Z^2$ with $|z(1)|>m$,
\[ \PPP{z\in A(\eta), \#(A(\eta)\cap B(z,m))\leq \beta m^2}\leq C_0 e^{-cm^2/\log m}.\]
\end{lemma} 

According to Lemma~\ref{Le:Jerison}, Equation \eqref{eq:UBeq2} is bounded by any power of $n^{-1}$. To deal with \eqref{eq:UBeq1}, we need to introduce some notation. In what follows, for any set $\Gamma\subset\Z^2$, we denote by $M^*_{n/2+h(n)}\left(\eta, \Gamma  \right)$ the number of \textit{random walks}, with initial configuration $\eta$,  satisfying the following two properties:
\begin{itemize}
\item{  the particle associated with the random walk hits $\partial R_{n/2}$ before exiting the aggregate;
\item the random walk intersects  $\Gamma$ before exiting $R_{n/2+h(n)}$.}
\end{itemize}
In particular, on the event $\left\{|z(1)|=t_n\right\}$ and by definition of $h(n)$,  we  have \[\#\left(A(n\ind{\{0\}\times \llbracket -n^\gamma, n^\gamma\rrbracket})\cap B(z,h(n))\right) \leq M^*_{n/2+h(n)}\left( n\ind{\{0\}\times \llbracket -n^\gamma, n^\gamma\rrbracket }, B(z,h(n))  \right)\] since no particle can escape $R_{n/2+h(n)}$. Therefore
\begin{multline*}
\PPP{\#\left(A(n\ind{\{0\}\times \llbracket -n^\gamma, n^\gamma\rrbracket})\cap B(z,h(n))\right)>\beta h^2(n), |z(1)|=t_n}\\
\leq \PPP{M^*_{n/2+h(n)}\left( n\ind{\{0\}\times \llbracket -n^\gamma, n^\gamma\rrbracket }, B(z,h(n))  \right)>\beta h^2(n)}.
\end{multline*}
It is clear that $M^*_{n/2+h(n)}\left( n\ind{\{0\}\times \llbracket -n^\gamma, n^\gamma\rrbracket }, B(z,h(n))  \right)$ is a.s.\,finite. Conditional on the fact that this random variable equals some integer $k$, the trajectories of the $k$ random walks \textit{are not independent} since they have to satisfy the first property mentioned above. However, one of the key ingredients is to remark that, conditional on the later event, these random walks \textit{are independent after they reach $\partial R_{n/2}$}, and so they evolve independently after reaching $\partial B(z,h(n))$. Now, recall that a random walk starting from some point $x\in \partial B(z,h(n))$ has a probability at least $\rho>0$, {wich does not depend on  $n$}, to hit $B(z,2h(n))\cap (R_{n/2+h(n)})^c$ when it exits $B(z,2h(n))$. In particular, it has a probability at least $\rho$ to hit the tile \[\tau(z)=B(z,2h(n))\cap \partial R_{n/2+h(n)}.\] Therefore, conditional on the fact that $M^*_{n/2+h(n)}\left( n\ind{\{0\}\times \llbracket -n^\gamma, n^\gamma\rrbracket }, B(z,h(n))  \right)$ is larger than $\beta h^2(n)$, the random variable $M^*_{n/2+h(n)}(n\ind{\{0\}\times \llbracket -n^\gamma, n^\gamma\rrbracket}, \tau(z))$ stochastically dominates a binomial distribution with parameters $(\beta h^2(n), \rho)$. Thus, there exists $I>0$ such that
\begin{multline*}
\PPP{\left.M^*_{\frac{n}{2}+h(n)}(n\ind{\{0\}\times \llbracket -n^\gamma, n^\gamma\rrbracket}, \tau(z)) \leq  \frac{\beta h^2(n)}{2}\right|M^*_{\frac{n}{2}+h(n)}(n\ind{\{0\}\times \llbracket -n^\gamma, n^\gamma\rrbracket}, B(z,h(n)))>\beta h^2(n)}\\
\leq \exp(-I h^2(n)).
\end{multline*}
The right-hand side is negligeable compared to any power of $n^{-1}$. Thus, it is sufficient to prove that for $n$ large enough, for any $L>0$,  
 \begin{equation}
\label{eq:aimupperbound}
 \PPP{M^*_{n/2+h(n)}(n\ind{\{0\}\times \llbracket -n^\gamma, n^\gamma\rrbracket}, \tau(z)) >  \frac{\beta}{2}h^2(n)} \leq n^{-L}.
\end{equation} 
Since
\[  M^*_{n/2+h(n)}(n\ind{\{0\}\times \llbracket -n^\gamma, n^\gamma\rrbracket}, \tau(z)) \leq M^*_{n/2+h(n)}(n\ind{I(\infty)}, \tau(z)),\]
we have to prove that $\PPP{M^*_{n/2+h(n)}(n\ind{I(\infty)}, \tau(z)) >  \frac{\beta}{2}h^2(n)}$ is lower than any power of $n^{-1}$. 

We will use the following adaptation of \cite[Lemma 2.5]{AG13short} stated in the context of series of Bernoulli random variables.

\begin{lemma}
\label{Le:AGchebychev}
Let $W_n, L_n, \widetilde{M}_n$ be independent random variables and let $\mathcal{A}_n$ be an event independent of $L_n$. Assume that $L_n$ and $M_n\overset{\text{law}}{=}\widetilde{M}_n$ are series of independent Bernoulli random variables with finite expectations such that $\mu_n=\EEE{M_n}-\EEE{L_n}\geq 0$   and write $L_n=\sum_{i\geq 0}Y_i^{(n)}$. If the following  holds:
\[(W_n+L_n)\ind{\mathcal{A}_n}\underset{\text{sto}}{\leq}\widetilde{M}_n,\]
then, for all $n\geq 0$, $\xi_n\in \R$ and $\lambda\in [0,\log 2]$, 
\[\PPP{W_n\geq \xi_n, \mathcal{A}_n} \leq \exp\left( -\lambda(\xi_n-\mu_n)+\lambda^2\left(\mu_n + 4\sum_{i\geq 0}\EEE{Y_i^{(n)}}^2   \right)  \right).\]
\end{lemma}

Similarly to (2.4) in \cite{AG13short}, it follows from the definition of $M^*_{n/2+h(n)}(n\ind{I(\infty)},\tau(z))$ that 
\[M^*_{n/2+h(n)}(n\ind{I(\infty)},\tau(z)) + M_{n/2+h(n)}(A_{n/2}(n\ind{I(\infty)}), \tau(z))\overset{\text{law}}{=}  \widetilde{M}_{n/2+h(n)}(n\ind{I(\infty)}, \tau(z)), \]
where $A_{n/2}(n\ind{I(\infty)})$ stands for the positions of settled particles before the associated random walks leave $R_{n/2}$ and where $\widetilde{M}_{n/2+h(n)}(n\ind{I(\infty)}, \tau(z))$ is an independent copy of ${M}_{n/2+h(n)}(n\ind{I(\infty)}, \tau(z))$.   
Let us denote by $\delta_I(n)$ the internal error of $A_{n/2}(n\ind{I(\infty)})$ on the strip $\Z_{n^\alpha}$, {\it i.e.}
\[\delta_I(n) = \max\left\{\frac{n}{2} - |z(1)|: z\in \left( R_{n/2}\setminus A_{n/2}\left(n\ind{I(\infty)}\right) \right)\cap \Z_{n^\alpha} \right\},\]
with the convention $\max\emptyset=0$. Notice that in the proof of the lower bound (see Section \ref{subsec:lowerbound}), no particle exits $R_{n/2}$; hence we actually have $\delta_I(n)\leq C\, \log(n)$ w.h.p.\,provided that $C$ is large enough. Here, we say that a sequence of events $(E_n)_{n\geq 0}$ occurs w.h.p.\,if $\PPP{E_n^c}$ is lower than any power of $n^{-1}$ for $n$ large enough.

Since $R_{n/2-\delta_I(n)}\cap \Z_{n^\alpha}\subset \A{\infty}$,  it follows that 
 \begin{align*}
\left(M^*_{n/2+h(n)}(n\ind{I(\infty)},\tau(z)) + M_{n/2+h(n)}\left(R_{n/2-\alpha_2\frac{h(n)}{2A}}{,\tau(z)}\right)\right)&\ind{\delta_I(n)\leq \alpha_2h(n)/2A} \\&\underset{\text{sto}}{\leq} \widetilde{M}_{n/2+h(n)}(n\ind{I(\infty)}, \tau(z)),
\end{align*}
where  $\alpha_2$ is some constant which will be chosen later. {We are now ready to apply Lemma~\ref{Le:AGchebychev} with}
 \[W_n=M^*_{n/2+h(n)}\left(n\ind{I(\infty)},\tau(z)\right),\qquad L_n=M_{n/2+h(n)}\left(R_{n/2-\alpha_2\frac{h(n)}{2A}},\tau(z)\right)\]\[\widetilde{M}_n=\widetilde{M}_{n/2+h(n)}(n\ind{I(\infty)}, \tau(z))\qquad\mbox{and}\qquad \mathcal{A}_n=\left\{\delta_I(n)\leq \alpha_2\frac{h(n)}{2A}\right\}.\]
Note that $L_n$ is independent of $\mathcal{A}_n$ by definition. It remains to control $\mu_n$ and $\sum_i\mathbb{E}[Y_i^{(n)}]^2$.
According to \eqref{eq:expectationbis}  and Lemma~\ref{Le:eqnegligible}, we know that
\begin{multline*}
\EEE{ \widetilde{M}_{n/2+h(n)}(n\ind{I(\infty)}, \tau(z))} - \EEE{ M_{n/2+h(n)}\left(R_{n/2-\alpha_2\frac{h(n)}{2A}}, \tau(z)\right)}\\
\begin{split}
& = n\EEE{M_{n/2+h(n)}(I(\infty), \tau(z))} - \EEE{ M_{n/2+h(n)}\left(R_{n/2-\alpha_2\frac{h(n)}{2A}}, \tau(z) \right)} \\
&  = \frac{\alpha_2}{A}h^2(n) + O(h(n)).
\end{split}
\end{multline*}
Moreover, it follows from the same computations as in Section~\ref{subsec:lowerbound} that
\[\sum_{y\in R_{n/2}}\PP_y(S(H(\Sigma_1))\in \tau(z))^2 \leq ch(n)^3. \]
Taking successively $A$ and $\alpha_2$ large enough, and proceeding exactly in the same spirit as in   \cite[Section 4.3]{AG13long}, it follows from Lemma~\ref{Le:AGchebychev} that  $\PPP{M^*_{n/2+h(n)}(n\ind{I(\infty)},\tau(z)) >  \frac{\beta}{2}h^2(n)}$ is negligeable compared to any power of $n^{-1}$. This concludes the proof of \eqref{eq:upperbound}.

{
{
\subsection{Proof of Theorem \ref{th:ShapeThPoisson} }
Since $A_n^*[\infty]$ and $A_n^\dag[\infty]$ are equally distributed, we only deal with $A_n^*[\infty]$. Recall that the aggregate $A_n^*[\infty]$ is generated by random walks starting from $\{0\}\times \Z$ and by a family of independent Poisson random variables $(N_i, i\in \Z)$ with parameter $n$.

Let $\gamma > \alpha > 0$ and $0 < \varepsilon < 1/2$ be fixed. The proof of Theorem \ref{th:ShapeThPoisson} mainly relies on Theorem \ref{th:ShapeTh} and on two other ingredients that we describe below. First, a.s., for $n$ large enough, all the events 
\[
\mathcal{E}_n = \mathcal{E}_n(\varepsilon,\gamma) := \bigcap_{|i|\leq n^\gamma} \big\{ | N_i - n | \leq n^{1/2+\varepsilon} \big\},
\] occur. This property is a consequence of the Borel-Cantelli lemma and of the classical concentration inequality
\[
\PPP{ | X - \lambda | \geq x } \leq 2 e^{-\frac{x^2}{2(\lambda+x)}},
\]
where $X$ is a Poisson random variable with parameter $\lambda$. Secondly, similarly to \eqref{eq:finiteinfinite}, we have
\[
\PPP{ A\left(\sum_{i\in \Z} N_i\ind{\{(0,i)\}} \right) \cap \Z_{n^\alpha} \neq A\left(\sum_{|i|\leq n^\gamma} N_i\ind{\{(0,i)\}}   \right) \cap \Z_{n^\alpha} } \leq  n^{-L}
\]
for any $L>0$ and for $n$ large enough. This implies that, a.s., for  $n$ large enough, 
\[
A \left( \sum_{i\in \Z} N_i\ind{\{(0,i)\}} \right) \cap \Z_{n^\alpha} =  A \left( \sum_{|i|\leq n^\gamma} N_i\ind{\{(0,i)\}} \right) \cap \Z_{n^\alpha}.
\]

Now, let us prove the upper bound in Theorem \ref{th:ShapeThPoisson}. A.s., for $n$ large enough, we have 
\begin{align*}
A_n^*[\infty] \cap \Z_{n^\alpha} & =  A \left( \sum_{|i|\leq n^\gamma} N_i\ind{\{(0,i)\}} \right) \cap \Z_{n^\alpha}\\
& \subset  A \left( m(n) \ind{\{0\}\times \llbracket -n^\gamma,n^\gamma\rrbracket} \right) \cap \Z_{n^\alpha} \hspace{0.5cm} \mbox{with $m(n):=n + n^{1/2+\varepsilon}$}\\
& \subset  A \left( m(n) \ind{I(\infty)} \right) \cap \Z_{m(n)^\alpha} \\
& \subset  R_{m(n)/2 + A\log(m(n))} \hspace{0.5cm} \mbox{ by Theorem \ref{th:ShapeTh}} \\
& \subset  R_{n/2+n^{1/2+\varepsilon}} ~.
\end{align*}

Similar arguments work for the lower bound. Indeed, a.s., for $n$ large enough,
\begin{align*}
A_n^*[\infty] \cap \Z_{n^\alpha} & =  A \left( \sum_{|i|\leq n^\gamma} N_i\ind{\{(0,i)\}} \right) \cap \Z_{n^\alpha}\\
& \supset  A \left( m(n) \ind{\{0\}\times \llbracket -n^\gamma,n^\gamma\rrbracket} \right) \cap \Z_{n^\alpha} \hspace{0.5cm} \mbox{with $m(n):=n - n^{1/2+\varepsilon}$}\\
& \supset  A \left( m(n) \ind{\{0\}\times \llbracket -m(n)^\gamma,m(n)^\gamma\rrbracket} \right) \cap \Z_{m(n)^{\alpha'}},
\end{align*}
where the parameter $\alpha'=\alpha'_n$ is chosen in such a way that $m(n)^{\alpha'}=n^\alpha$. To show that the aggregates $A(m(n) \ind{\{0\}\times \llbracket -m(n)^\gamma,m(n)^\gamma\rrbracket})$ and $A(m(n) \ind{I(\infty)})$ coincide on the strip $\Z_{m(n)^{\alpha'}}$, we have to be careful because the power $\alpha'_n$ depends on $n$. However, it is possible to do it since the random integer $n_0(\alpha)$ from which
\[
 A \left( m(n) \ind{\{0\}\times \llbracket -m(n)^\gamma,m(n)^\gamma\rrbracket} \right) \cap \Z_{m(n)^{\alpha}} = A \left( m(n) \ind{I(\infty)} \right) \cap \Z_{m(n)^{\alpha}}
\]
is increasing w.r.t. $\alpha$ and, for any $n$ large enough, $\alpha'_n \in (\alpha,2\alpha)$. Henceforth, a.s and for $n$ large enough,
\begin{align*}
A_n^*[\infty] \cap \Z_{n^\alpha} & \supset  A \left( m(n) \ind{I(\infty)} \right) \cap \Z_{m(n)^{\alpha'}}\\
& \supset  R_{m(n)/2 - A\log(m(n))} \cap \Z_{m(n)^{\alpha'}} \hspace{0.5cm} \mbox{ by Theorem \ref{th:ShapeTh}} \\
& \supset  R_{n/2 - n^{1/2+\varepsilon}} \cap \Z_{n^\alpha} ~.\\
\end{align*}

\begin{Rk}
\label{rk:shapetheorempoisson}
Let us shortly explain why the approach used to obtain logarithmic fluctuations in the shape theorem for $A_n[\infty]$ in Subsections \ref{subsec:lowerbound} and \ref{subsec:upperbound} fails when considering $A^*_n[\infty]$ (and therefore $A^\dag_n[\infty]$) instead of $A_n[\infty]$. In what follows,  we denote by $\eta$ the initial configuration of $A^*_n[\infty]$, {\it i.e.}
\[\eta=\sum_{i\in \Z}N_i\ind{\{(0,i)\}},\]
and we simply write $A(\eta)$ instead of $A_n^*[\infty]$ to make explicit the dependence in $\eta$.  Given a realization of $\eta$, we recall that $M_r(\eta, B)$ (resp. $W_r(\eta,B)$) denotes the number of  random walks (resp. particles) with initial configuration $\eta$, that hit $B$ before or when they exit $R_r$.   
  
In order to proceed in the same spirit as in Section~\ref{subsec:lowerbound}, one has to consider, for any tile $\tau\in \mathcal{T}_{k\log(n)}$, the random variables $M_{k\log(n)}(\eta, \tau)$ and $W_{k\log(n)}(\eta, \tau)$ instead of $M_{k\log(n)}(n\ind{I(\infty)}, \tau)$ and $W_{k\log(n)}(n\ind{I(\infty)}, \tau)$, respectively. 
Then, to apply Lemma~\ref{Le:AGshortL2.4}, the number $\mu(\tau)$, as defined in \eqref{eq:defmu}, has to be replaced with the random variable
\[\mu(\tau|\eta)=\EE\left[\left.M_{k\log(n)}(\eta, \tau)\right|\eta\right]-\EEE{M_{k\log(n)}(R_{k\log(n)}\setminus\mathcal{Z},\tau)},\]
where $\EE[\cdot|\eta]$ denotes the conditional expectation w.r.t.\,$\eta$. One needs to show that the counterpart of \eqref{eq:LBmutau}:
\begin{equation}
\label{eq:fails} \mu(\tau|\eta)\geq cA\log(n)^2
\end{equation}
holds w.h.p., for some $c>0$, for every tile centered at a point at a distance at least $A\log(n)$ from $\partial R_{\frac{n}{2}}$.

The expectation of $\mu(\tau | \eta)$ is $\mu(\tau)$. But, because the variance of $N_i$, $i\in \Z$,  equals $n$ the variance of $\mu(\tau | \eta)$ is at least $c.n$ for some constant $c>0$. By looking carefully at the proof of \eqref{eq:LBmutau}, we can see that \eqref{eq:fails} must fail with asymptotically non negligible probability. Thus, the approach fails to establish the lower bound for $A^*_n[\infty]$. Since it relies on the lower bound with logarithmic fluctuations, the argument presented in Subsection \ref{subsec:upperbound} to derive the upper bound cannot be adapted either.   
\end{Rk} 
}

\section{Construction of the Directed IDLA Forest}
\label{sec:forest}

In this section, we introduce a new random forest $\FF_\infty$ spanning all $\Z^2$ and based on the IDLA protocol with sources on $I(\infty)=\{0\}\times\Z$. For any  $n,M$, we first build in Section~\ref{sect:ForetnM} a random forest ${\FF_n^\dag[M]}$ w.r.t.\,the aggregate $A_{n}^\dag[M]$. By letting $M\to\infty$ (vertical limit), we define a random forest $\FF_n$ as the limit of the sequence $({\FF_n^\dag[M]})_{M\geq 0}$. The existence of $\FF_n$ is based on a non-trivial stabilization result (Proposition~\ref{prop:existenceforest}) which is stated in the same spirit as Theorem~\ref{th:existence}. The key argument to define the limiting forest $\FF_n$ lies on the fact that the infinite aggregate $A_{n}^\dag[\infty]$ is made up with finite connected components. This property prevents the existence of \textit{chain of changes} (see Section~\ref{sect:ChainChanges}) which could come from far levels and perturbate the evolution of $({\FF_n^\dag[M]})_{M\geq 0}$ in the neighborhood of the origin. It will be shown that the sequence $(\FF_n)_{n\geq 1}$ is consistent (Lemma~\ref{lem:Fn}). By letting $n\to\infty$ (horizontal limit), this fact allows us to define easily our random directed IDLA forest $\FF_\infty$. This section ends with some properties of  $\FF_\infty$ (Theorem~\ref{th:forestIDLA}).

\subsection{Each aggregate $A_n^\dag[M]$ generates a forest}
\label{sect:ForetnM}

Let $n,M\geq 1$ be fixed and let $\kappa = \sum_{|i|\leq M} \#\NN_{i}([0,n])$.  As in the proof of Lemma~\ref{lem:inclusion}, we index all the particles starting from level $|i|\leq M$ by some integer $j=1,\ldots,\kappa$ according to their starting times $0<\tau_1<\ldots<\tau_{\kappa}<n$. For $j=1,\ldots,\kappa$, we denote by $A[j]$ the aggregate obtained until (or at) time $\tau_j$. In particular, we have  $A[0]=\emptyset$ and $A[\kappa]=A_{n}^\dag[M]$. 

We define a (finite) random forest with vertices in $A_{n}^\dag[M]$ inductively as follows. First, we let $\FF_n[M,0]=(\emptyset,\emptyset)$. Then, for some $1\leq j\leq\kappa$, assume that a graph  $\FF_n[M,j-1]$ is built, with set of vertices and edges denoted by $V[j-1]$ and $E[j-1]$, respectively. Let $z$ be the site which is added by particle $j$, {\it i.e.}\,$A[j]=A[j-1] \cup \{z\}$.
\begin{itemize}
\item[$\bullet$] If $z=(0,i)$ and if particle $j$ is (the first one) which is sent from level $i$ then $z$ is the root of a new tree in the graph. In this case, we set $\FF_n[M,j]=(V[j],E[j])$, where
\[
V[j] = V[j-1] \cup \{z\} \; \mbox{ and  } \; E[j] = E[j-1]~.
\]
\item[$\bullet$] Otherwise, the site of $\{0\}\times\Z$ from which particle $j$ starts already belongs to $A[j-1]$ and we let $z'$ as the last site of $A[j-1]$ which is visited by particle $j$ before reaching $z$. Then, we set $\FF_n[M,j]=(V[j], E[j])$, where
\[
V[j] = V[j-1] \cup \{z\} \; \mbox{ and  } \; E[j] = E[j-1] \cup \{(z',z)\} ~.
\]
\end{itemize}
In other words, from $\FF_n[M,j-1]$ to $\FF_n[M,j]$ we merely add the new site created by particle $j$ and the (directed) edge from which this new site is reached. In what follows, we set
\[
{\FF_n^\dag[M]} = \FF_n[M,\kappa] ~.
\]
Figure~\ref{fig:finiteforest} gives a realization of ${\FF_n^\dag[M]}$. {Observe that ${\FF_n^\dag[M]}$ is a (finite) random forest with vertices in $A_n^\dag[M]$.

\begin{Rk}
\label{lem:BuiltForest}
Let $n,M\geq 1$. The following properties hold a.s.
\begin{enumerate}[(i)]
\item The set of vertices of ${\FF_n^\dag[M]}$ is $A_{n}^\dag[M]$; 
\item The edges of ${\FF_n^\dag[M]}$ are edges of the square lattice $\Z^2$ plus a direction;
\item The random graph ${\FF_n^\dag[M]}$ is a (finite) union of directed trees with roots in $I(\infty)$. 
\end{enumerate}

Properties (i)-(ii) are satisfied by construction. {Property (iii) comes from the fact that the random graph ${\FF_n^\dag[M]}$ contains no loop since a site cannot be added twice to the aggregate}.

\end{Rk}
}
\begin{center}
\begin{figure}
\begin{center}
\begin{tabular}{c}
\includegraphics[height=2.8cm]{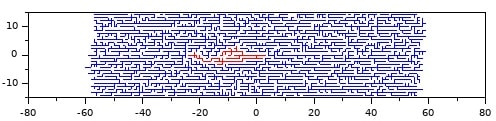}\\
\includegraphics[height=2.8cm]{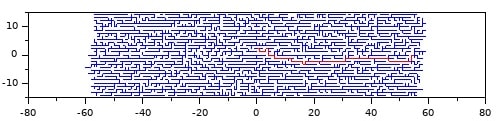}
\end{tabular}\begin{tabular}{c}
\includegraphics[height=5.6cm]{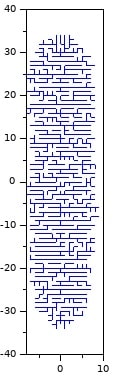}
\end{tabular}
\end{center}
\caption{The top left corner depicts a realization of the random forest ${\FF_{120}^\dag[120]}$ with particles starting from levels $|i|\leq 120$ and during the time interval $[0,120]$, viewed through the strip $\Z_{15}$. The tree of ${\FF_{120}^\dag[120]}$ containing the origin is in red. A second realization of ${\FF_{120}^\dag[120]}$ is given on the bottom left corner. The branch passing through $(55,0)$ (red) remains close to the $x$-axis and comes from the source $(0,2)$. A realization of ${\FF_{15}^\dag[30]}$ is depicted on the right. One can imagine that the vertical edges at the top of ${\FF_{15}^\dag[30]}$ are due to border effects and will not be present in the limiting forest $\FF_{15}$.}\label{fig:finiteforest}
\end{figure}
\end{center}

\subsection{Absence of infinite chain of changes and stabilization}
\label{sect:ChainChanges}

Let $n,M\geq 0$. For both constructions of aggregates $A_n[M]$ and $A^*_n[M]$, recall that the particles are sent w.r.t.\,the usual order, {\it i.e.}\,from levels $0$, thus $\pm 1$, $\pm 2$ and so on by moving away from the origin step by step. Hence, these aggregates are first built around the origin and thus grow mainly from their upper and lower parts. In Theorem~\ref{th:existence}, we proved that the particles which come from far levels cannot visit a neighborhood of the origin in this setting.

The situation is quite different for aggregates $A_{n}^\dag[M]$, $M\geq 0$. Indeed, particles are sent according to the clocks given by $(\NN_{i})_{i\in\Z}$ and some pathological situations, {as described below}, may occur (see Figure \ref{fig:discrepanciesforest}).  For $M'> M\geq 0$, any particle sent from a level $M<i\leq M'$ works for (the growth of) $A_{n}^\dag[M']$ but not for $A_{n}^\dag[M]$. Hence, this particle may create several discrepancies between the forests ${\FF_n^\dag[M]}$ and ${\FF_n^\dag[M']}$ through a mechanism called a \textit{chain of changes} that we describe now.  

\begin{center}
\begin{figure}
\begin{center}
\includegraphics[width=4cm,height=7cm]{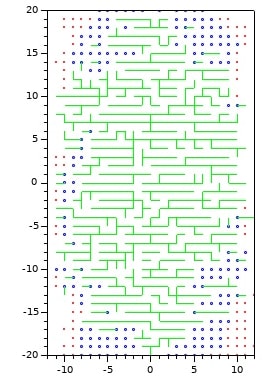}
\end{center}
\caption{\label{fig:discrepanciesforest} Realizations of the forests ${\FF_{20}^\dag[20]}$ and ${\FF_{20}^\dag[50]}$, defined on the same time interval $[0,20]$, with different sets of sources, and restricted to the strip $\Z_{20}$, are depicted. The associated aggregates are coupled in the sense that they are based on the same clocks and random walks {with level $|i|\leq 20$}.  In particular, $A^\dag_{20}[20]$ is included in $A^\dag_{20}[50]$. The edges created in both forests by the same particles are depicted in green. The red points are vertices of $A^\dag_{20}[50]\!\setminus\!A^\dag_{20}[20]$. The blue circles represent vertices in $A^\dag_{20}[20]$ (and then also in $A^\dag_{20}[50]$) which are reached by different particles in both aggregates and whose corresponding edges may differ in both forests ${\FF_{20}^\dag[20]}$ and ${\FF_{20}^\dag[50]}$. These blue vertices are possible discrepancies generated by chains of changes between forests ${\FF_{20}^\dag[20]}$ and ${\FF_{20}^\dag[50]}$.}
\end{figure}
\end{center}

Assume that a particle, referred to as particle 1, starts at time $t_1\in(0,n)$ (from a level $M<|i_1|\leq M'$) and adds a site $z_1$ to $A_{t_{1}-}^\dag[M']$. The aggregate at time $t_1$ becomes
\[
A_{t_{1}}^\dag[M'] = A_{t_{1}-}^\dag[M']\cup\{z_1\}
\]
while $A_{t_{1}}^\dag[M]$ remains unchanged. In the above equation, the set  $A_{t_{1}-}^\dag[M']$ denotes the (current) aggregate produced just before sending particle 1. The site $z_1$ is a \textit{discrepancy} at time $t_1$ between aggregates $A_{t_{1}}^\dag[M]$ and $A_{t_{1}}^\dag[M']$. If there is no other particles starting from a level $|i|\leq M$, at time $t\in (t_1,n)$ and going through $z_1$ then, at the final time $n$, the site $z_1$ constitutes a discrepancy (created by  particle 1) between the aggregates $A_{n}^\dag[M]$ and $A_{n}^\dag[M']$. It also defines a discrepancy between the forests ${\FF_n^\dag[M]}$ and ${\FF_n^\dag[M']}$. Otherwise, we set
\[
t_2 = \min \big\{ t \in (t_1,n) : \, \mbox{a particle, starting from a level $|i| \leq M$ at time $t$, goes through $z_1$} \big\} ~.
\]
The particle starting from time $t_2$ is referred to as particle 2. This particle works for both aggregates. By definition, it adds the site $z_1$ to $A_{t_{2}-}^\dag[M]$, so that the aggregate at time $t_2$ becomes  $A_{t_{2}}^\dag[M]=A_{t_{2}-}^\dag[M]\cup\{z_1\}$. Thus, it continues its trajectory until adding a site $z_2$ (but only) to $A_{t_{2}-}^\dag[M']$ which then becomes $A_{t_{2}}^\dag[M']=A_{t_{2}-}^\dag[M']\cup\{z_2\}$. At this time:
\begin{itemize}
\item the site $z_1$ now belongs to both aggregates but it could be reached via two different edges respectively in $A_{t_{2}}^\dag[M]$ and $A_{t_{2}}^\dag[M']$ so that the forests ${\FF_n^\dag[M]}$ and ${\FF_n^\dag[M']}$ may differ at the edge leading to $z_1$;
\item the site $z_2$ is become a discrepancy between both aggregates at time $t_2$. This discrepancy is generated via a relay between particle 1 and particle 2.  
\end{itemize}
Thus, we iterate this step while the current discrepancy is visited by a new particle starting from a level $|i|\leq M$. After a random number $\ell$ of steps (a.s.\,finite), we finally get the set of possible discrepancies between the forests ${\FF_n^\dag[M]}$ and ${\FF_n^\dag[M']}$, generated by particle 1. This set consists of edges leading to $z_1,\ldots,z_{\ell}$ and the final vertex $z_{\ell}$ itself. The mechanism producing this set of discrepancies is called a chain of changes, initiated by particle 1, between the forests ${\FF_n^\dag[M]}$ and ${\FF_n^\dag[M']}$. Notice that the aggregates $A_{n}^\dag[M]$ and $A_{n}^\dag[M']$ may have other chains of changes initiated by other particles starting from levels $M<|i|\leq M'$.

Roughly speaking, the existence of an infinite chain of changes involving an infinite number of relaying particles and initiated by a ``Big Bang particle'', {\it i.e.} a particle coming from a level arbitrarily far from the origin and born arbitrarily early, could modify infinitely often (in $M$) the forests ${\FF_n^\dag[M]}$, for $M\geq 0$, in the neighborhood of the origin. Proving that such infinite chain of changes does not exist with probability $1$ leads to the next stabilization result and to the existence of the random forest $\FF_n$.

\begin{proposition}
\label{prop:existenceforest}
Let $K\geq 1$. Then, a.s.\,there exists some (random) integer $M_0(K)$ such that, for any $M'>M\geq M_0(K)$, we have
\[
{\FF_n^\dag[M]} \cap \Z_K = {\FF_n^\dag[M']} \cap \Z_K ~.
\]
\end{proposition}

Proposition~\ref{prop:existenceforest} allows us to define a.s.\,the random forest $\FF_n$ as the increasing union
\[
\FF_n = \bigcup_{K\geq 0} \! \uparrow {\FF_n^\dag[M_0(K)]} \cap \Z_K,
\]
with vertex set $V(\FF_n)=A_n^\dag[\infty]$.

\begin{prooft}{Proposition~\ref{prop:existenceforest}}
Let $K\geq 1$. Corollary~\ref{cor:finitecc} provides the almost sure existence of an infinite subset $\mathcal{L}\subset\Z$ such that, for any $i\in\mathcal{L}$,
\[
A_n^\dag[\infty] \cap (\Z\times \{i\}) = \emptyset ~.
\]
Now, let us define $M_0(K)$ as the following (random) integer:
\[
M_0(K) = \max \Big\{ \min\{i\in\mathcal{L} : \, i \geq K\} , -\max\{i\in\mathcal{L} : \, i \leq -K\} \Big\} ~.
\]
For any $M'>M\geq M_0(K)$, there is no chain of changes initiated by a particle starting from some level $M<|i|\leq M'$ between the forests ${\FF_n^\dag[M]}$ and ${\FF_n^\dag[M']}$ which visits the strip $\Z_K$. This implies that, for any $M'>M\geq M_0(K)$, the forests ${\FF_n^\dag[M]}$ and ${\FF_n^\dag[M']}$ coincide on the strip $\Z_K$. 

\end{prooft}

The following lemma is a direct consequence of our construction.

\begin{lemma}
\label{lem:Fn}
The sequence of forests $(\FF_n)_{n\geq 1}$ is consistent in the sense that
\begin{equation}
\label{consistent}
\mbox{a.s. }\, \forall n \geq 1, \quad  V(\FF_n)\subset V(\FF_{n+1}) \quad \text{and} \quad E(\FF_n)\subset E(\FF_{n+1}) ~.
\end{equation}
Moreover, for any $n\geq 1$, $\FF_n$ is a.s.\,made up with infinitely many directed trees rooted at $I(\infty)$.\end{lemma}

\begin{prooft}{Lemma~\ref{lem:Fn}}
Inclusion (\ref{consistent}) is an immediate consequence of our construction (the same argument was used in the proof of Lemma~\ref{lem:inclusion} (ii)). The fact that $\FF_n$ is made up with infinitely many trees is due to Corollary~\ref{cor:finitecc} which asserts that $A^\dag_n[\infty]=V(\FF_n)$ is itself made up with an infinite number of disjoint connected components.
\end{prooft}

\subsection{The directed infinite-volume IDLA forest}
\label{subsec:forestidla}

Equation \eqref{consistent} allows us to define a.s.\,the \textit{directed infinite-volume IDLA forest}
\[
\FF_\infty = \bigcup_{n\geq 1} \! \uparrow \FF_n.
\]
The following theorem states the main properties of $\FF_\infty$. 

\begin{theorem}
\label{th:forestIDLA}
The directed infinite-volume IDLA forest $\FF_\infty$ satisfies the following properties:
\begin{enumerate}[(i)]
\item A.s. the random forest $\FF_\infty$ spans the whole set $\Z^2$, {\it i.e.} $V(\FF_\infty)=\Z^2$, and its edge set $E(\FF_\infty)$ is made up with edges of the square lattice $\Z^2$ plus a direction;
\item A.s. $\FF_\infty$ is a countable (infinite) union of directed trees rooted at $I(\infty)$;
\item the distributions of $\FF_\infty$ and $\FF_n$, $n\geq 1$ are invariant w.r.t.\,vertical translations;
\item the distributions of $\FF_\infty$ and $\FF_n$, $n\geq 1$ are mixing w.r.t.\,vertical translations;
\item the distributions of $\FF_\infty$ and $\FF_n$, $n\geq 1$ are symmetric invariant w.r.t.\,the $y$-axis and w.r.t.\,$S_{k/2}$, $k\in \Z$. 
\end{enumerate}
\end{theorem}

Item (iii) means that the forests $\FF_\infty$ and $\FF_n$, $n\geq 1$, are no longer sensitive to the border effects whereas the ${\FF_n^\dag[M]}$'s are (see the right-hand side of Figure~\ref{fig:finiteforest}). This result is one of the original motivation of this paper. Let us also remark that Item (iii) is not an immediate consequence of the translation invariance of the aggregates  $A_{n}^\dag[\infty]$, $n\geq 1$, (see Proposition~\ref{prop:transfo}) since $\FF_\infty$ actually is a richer object including edges which depend on the trajectories of particles.

\begin{prooft}{Theorem~\ref{th:forestIDLA}} The fact that $E(\FF_\infty)$ is made up with edges of the square lattice $\Z^2$ plus a direction and Item  (ii) are direct consequences of Proposition~\ref{prop:existenceforest} and Lemma~\ref{lem:Fn}. To prove that $V(\FF_\infty)=\Z^2$, we use the lower bound appearing in  {Theorem \ref{th:ShapeThPoisson}}. Indeed, {Theorem~\ref{th:ShapeThPoisson}} implies that, for any finite subset $S$ of $\Z^2$,
\begin{equation}
\label{eq:convshapetheorem}
\PPP{S\subset A_n^\dag[\infty]} \conv[n]{\infty}1 ~.
\end{equation}
Since $\PPP{S\subset V(\FF_\infty)} \geq \PPP{S\subset V(\FF_n)} = \PPP{S\subset A_n^\dag[\infty]}$ for any $n$, we get $\PPP{S\subset V(\FF_\infty)}=1$. This concludes the proof of (i). 

Let us prove (iii). Let $K\subset\R^2$ be a compact set and let $k\geq 0$. Recall that $\tau_k \FF_\infty$ (resp. $\tau_k \FF_n$) denotes the directed IDLA forest $\FF_\infty$ (resp. $\FF_n$) translated w.r.t.\,the vector $(0,k)$. According to \cite[Theorem 2.1.3]{SW}, it is is sufficient to prove that both forests $\FF_\infty$ and $\tau_k \FF_\infty$ (resp. $\FF_n$ and $\tau_k \FF_n$) have the same probability to intersect $K$, where all these graphs are seen as subsets of $\R^2$. First, assume that this holds for the $\FF_n$'s, {\it i.e.}
\begin{equation}
\label{InvarianceTranslationFn}
\PPP{\FF_n \cap K \neq \emptyset} = \PPP{\tau_k \FF_n \cap K \neq \emptyset} ~.
\end{equation}
{Notice that if $(C+B(0,1))\cap\mathbb{Z}^2\subset A^\dag_n[\infty]$, then $\mathcal{F}_n\cap C=\mathcal{F}_\infty\cap C$.} Now, let $\varepsilon>0$.    Thanks to \eqref{eq:convshapetheorem} applied to the set $(C{+} B(0,1))\cap\Z^2$, with $C=K$ and $C=\tau_{-k}K$, there exists an integer $n_0$ such that  
\[
\left| \PPP{\FF_\infty \cap K \neq \emptyset} - \PPP{\FF_{n_0} \cap K \neq \emptyset} \right| \; \mbox{ and } \; \left| \PPP{\tau_k \FF_\infty \cap K \neq \emptyset} - \PPP{\tau_k \FF_{n_0} \cap K \neq \emptyset} \right|
\]
are both smaller than $\varepsilon$. By \eqref{InvarianceTranslationFn}, we deduce the translation invariance in distribution for the limiting forest $\FF_\infty$.

It then remains to prove \eqref{InvarianceTranslationFn} for any $n\geq 1$. Let $n\geq 1$ and $M$ large enough so that $K\cap\Z^2 \subset\Z_M$. By Proposition~\ref{prop:existenceforest}, there exists a (deterministic) integer $M_0=M_0(n,\varepsilon)$ such that, with probability at least $1-\varepsilon$, we have for any $M'\geq M_0$,
\begin{equation}
\label{Stab1-dag}
\FF_n \cap \Z_M = {\FF_n^\dag[M']} \cap \Z_M \; \mbox{ and } \; \tau_k \FF_n \cap \Z_M = \tau_k {\FF_n^\dag[M']} \cap \Z_M ~.
\end{equation}
Then,
\begin{multline*}
| \PPP{\FF_n \cap K \neq \emptyset} - \PPP{\tau_k \FF_n \cap K \neq \emptyset} | 
 \leq \, \left| \PPP{ {\FF_n^\dag[M_0]}  \cap K \neq \emptyset} - \PPP{\tau_k  {\FF_n^\dag[M_0]}   \cap K \neq \emptyset} \right| + 2 \varepsilon ~.
\end{multline*}
Now let us increase the forests ${\FF_n^\dag[M_0]}$ and $\tau_k {\FF_n^\dag[M_0]}$ as follows. Let $\mathfrak{F}_1$ be the random forest obtained by sending the particles used to build ${\FF_n^\dag[M_0]}$ plus those from levels $M_0+1,\ldots,M_0+k$ (according to their own PPP's). Besides, let $\mathfrak{F}_2$ be the translation w.r.t.\,the vector $(0,k)$ of the forest induced by the particles used to build ${\FF_n^\dag[M_0]}$ plus those from levels $-M_0-k,\ldots,-M_0-1$. On the one hand, using \eqref{Stab1-dag}, we know that the forests $\mathfrak{F}_1$ and $\mathfrak{F}_2$ coincide respectively with ${\FF_n^\dag[M_0]}$ and $\tau_k {\FF_n^\dag[M_0]}$ on the strip $\Z_M$, with probability at least $1-\varepsilon$. This implies
\[
| \PPP{\FF_n \cap K \neq \emptyset} - \PPP{\tau_k \FF_n \cap K \neq \emptyset} | \leq | \PPP{\mathfrak{F}_1 \cap K \neq \emptyset} - \PPP{\mathfrak{F}_2 \cap K \neq \emptyset} | + 4 \varepsilon ~.
\]
On the other hand, the forests $\mathfrak{F}_1$ and $\mathfrak{F}_2$ are produced by the same IDLA protocol from the same levels $i=-M_0,\ldots,M_0+k$, {with i.i.d. Poisson clocks}, during the time interval $[0,n]$. {In particular, the forests $\mathfrak{F}_1$ and $\mathfrak{F}_2$ are identically distributed}. This concludes the proof of (iii).

Let us prove (iv). By \cite[Theorem 9.3.2]{SW}, it is enough to check that
\begin{equation*}
\label{MixingCond}
\lim_{k\to \infty} \PP\left(\FF_\infty\cap\left(C_1\cup\tau_k C_2\right) = \emptyset\right) = \PP\left(\FF_\infty\cap C_1 = \emptyset\right) \PP\left(\FF_\infty\cap C_2 = \emptyset\right)
\end{equation*}
holds for any compact sets $C_1$ and $C_2$ in $\R^2$. Let $\varepsilon>0$ and $C_1, C_2$ be two compact sets in $\R^2$. Let $r>0$ be  such that $C_1 \cup C_2$ is included in the ball $B(0,r-1)$. By \eqref{eq:convshapetheorem} we can choose $n$ so that 
\[
\PP\left( B(0,r) \cap \Z^2 \subset A^\dag_n[\infty] \right) \geq 1 - \varepsilon ~.
\]
Replacing $r-1$ with $r$ in the above probability allows us to take into account all the edges incident to the vertices of $(C_1 \cup C_2)\cap \Z^2$. On the above event, $\FF_\infty \cap C_i$ and $\FF_n \cap C_i$ are equal, for $i\in\{1,2\}$. The translation invariance in distribution of $A^\dag_n[\infty]$ (Proposition~\ref{prop:transfo}) implies that the probability of the event $\{ \tau_k (C_1 \cup C_2) \cap \Z^2 \subset A^\dag_n[\infty] \}$ is also larger than $1-\varepsilon$, for any integer $k$ (which does not depend on $\varepsilon$). Hence, with probability larger than $1-2\varepsilon$, $\FF_\infty \cap (C_1\cup\tau_k C_2)$ and $\FF_n \cap (C_1\cup\tau_k C_2)$ are equal. Henceforth,
\begin{multline*}
| \PP\left(\FF_\infty \cap (C_1\cup\tau_k C_2) = \emptyset \right) - \PP\left( \FF_\infty \cap C_1 = \emptyset \right) \PP\left( \FF_\infty \cap C_2 = \emptyset \right) | \\
 \leq | \PP\left(\FF_n \cap (C_1\cup\tau_k C_2) = \emptyset \right) - \PP\left(\FF_n \cap C_1 = \emptyset \right) \PP\left(\FF_n \cap C_2 = \emptyset\right) | + 4 \varepsilon ~.
\end{multline*}
Let us consider the event $D_{n,k}$ defined by ``there is no connected component of $A^\dag_n[\infty]$ overlapping simultaneously $C_1$ and $\tau_k C_2$''. By Corollary~\ref{cor:finitecc}, given $n$, the probability of $D_{n,k}$ is larger than $1-\varepsilon$ for any $k$ large enough. For any such integer $k$, the events $\{\FF_n \cap C_1=\emptyset\}$ and $\{\FF_n \cap \tau_k C_2=\emptyset\}$ are independent on $D_{n,k}$. Thus
\[
\PP\left( \{\FF_n \cap (C_1\cup\tau_k C_2) = \emptyset\} \cap D_{n,k} \right)  
 = \PP\left(\{\FF_n \cap C_1 = \emptyset\} \cap D_{n,k} \right) \PP\left(\{\FF_n \cap \tau_k C_2 = \emptyset\} \cap D_{n,k}\right) ~.
\]
Since $\FF_n$ is invariant in distribution w.r.t.\,vertical translations, we have 
\[
| \PP\left(\FF_n \cap (C_1\cup\tau_k C_2) = \emptyset \right) - \PP\left(\FF_n \cap C_1 = \emptyset \right) \PP\left(\FF_n \cap C_2 = \emptyset\right) | \leq 3 \varepsilon.
\]
Therefore
\[
| \PP\left(\FF_\infty \cap (C_1\cup\tau_k C_2) = \emptyset \right) - \PP\left(\FF_\infty \cap C_1 = \emptyset \right) \PP\left(\FF_\infty \cap C_2 = \emptyset\right) | \leq 7 \varepsilon ~.
\]

(v) By construction, the forests $\FF_n$ and $\FF_\infty$ are invariant under symmetries w.r.t.\, $x$-axis and $y$-axis. Since $S_{k/2}=\tau_k\circ S_0$ for any $k\in \Z$, it follows from (iii) that $\FF_n$ and $\FF_\infty$ are invariant under any horizontal symmetries.
\end{prooft}

\subsection{{Open questions}}
\label{sect:conjectures}
{We end this section with four open questions about the directed IDLA forest $\FF_\infty$.}

\subsubsection{Question 1}
{The main contribution of our paper is the construction of a stationary (w.r.t. vertical translations) random forest based on an Internal DLA protocol.  In the context of the External one, a similar object is introduced in \cite{PYZ20}. This object is an infinite stationary DLA on the upper half planar lattice growing from an infinite line. The construction of such a model is based on a central object which is a stationary version of the harmonic measure describing the growth of the aggregate.  The strategy developped in \cite{PYZ20} should allow us to define a dynamical infinite model $(\mathcal{F}_t)$ according to an Internal DLA protocol and with particles starting from an infinite line. This model should be based on a suitable harmonic measure, namely $\mathcal{H}_B(x,y)=\sum_{z\in I(\infty)}\PP_z(\xi_{\tau_B}=x, \xi_{\tau_B-1}=y)$, where $\tau_B$ is the first time that a random walk $\xi^z$ starting at $z$ hits $B\subset \Z^2$. Following the beginning of Section 3 in \cite{PYZ20}, a graphical thinned representation using the measure $\mathcal{H}_B$ should lead to a stabilization result (in the same spirit as Theorem 4 in \cite{PYZ20}), which is the main ingredient to prove the existence of the model. }

{The approach in \cite{PYZ20} fits very well in the context of External DLA because the harmonic measure is the natural tool to study its growth. But this approach requires a big preparative work on the harmonic measure (see \cite{PYZ21,PZ19}) and appears to be more artificial in the context of Internal DLA.   
 However, provided that  the dynamical infinite model based on an Internal DLA exists, it is a natural question to know if it coincides with our random forest. }

\subsubsection{{Question 2}} Consider the infinite IDLA tree $\TT_{\infty}$ rooted at the origin (as described in the Introduction) and focus on its branches only through the ball $B((n,0),R)$, where $R$ is fixed and where $n$ is intended to go to infinity.  The idea is that the radial character of its branches restricted to $B((n,0),R)$ should fade away as $n\to\infty$ since $R$ is constant. Hence, the infinite IDLA tree $\TT_{\infty}$ restricted to the ball $B((n,0),R)$ should look like to a directed forest with direction the vector $(-1,0)$; as if, roughly speaking, the root of the tree is sent to infinity according to the vector $(-1,0)$. Approximating the distribution of a tree (locally and far away from the root) by the distribution of a directed forest is classical in the literature, see \textit{e.g.}\,\cite{BB07,C18}. We {think} that the directed IDLA forest $\FF_\infty$ is the natural candidate to approximate the distribution of the tree $\TT_{\infty}$. More precisely, {is it true} that the infinite IDLA tree $\TT_{\infty}$ and the directed IDLA forest $\FF_\infty$, both restricted to $B((n,0),R)$, are asymptotically equally distributed?\\

\subsubsection{{Question 3}} Directed forests in $\R^d$ may coalesce or not according to the dimension $d$, see \textit{e.g.}\,\cite{CT,FLT,SS}. But whatever the dimension, in the backward sense ({\it i.e.} in the opposite direction to the one in which the branches coalesce), branches are finite (for the models previously cited) so that the directed forests do not contain bi-infinite branches. In regards to the directed IDLA forest $\FF_\infty$, branches coalesce when they get closer to the source axis $I(\infty)$ so that the backward sense is moving away from $I(\infty)$. {It is a natural question to show} that the same holds for the directed IDLA forest $\FF_\infty$, i.e. a.s. all the (infinitely many) trees making up the directed IDLA forest and rooted on the axis $I(\infty)$ are finite. See, for instance, the tree associated with the origin in Figure \ref{fig:finiteforest}. {The classical strategy to deal with such question consists in adapting a Burton-Keane type argument. In our context, the main difficulties to do it are the lack of finite energy property (allowing to locally change the forest; see Definition 2 in \cite{BK89} for a precise statement of this property) and the fact that $\FF_\infty$ is not invariant in distribution w.r.t.\,\textit{horizontal} translations.} 

Notice that, in \cite{BKP19}, the authors proved in that all the trees are a.s.\,finite in the forest they defined. The counterpart of this result for our IDLA forest is an open and challenging question.

\subsubsection{{Question 4}} A challenging question about the infinite IDLA tree $\TT_{\infty}$ is the existence of (many) infinite branches with asymptotic directions. Following the strategy initiated by Howard and Newman \cite{HN01}, the key point would be to control the fluctuations w.r.t.\,the segment $[0,z]$ (in $\R^2$) of the branch in $\TT_{\infty}$ joining the root $0$ to any given vertex $z\in\Z^2$, with $|z|_2 \gg 1$. This question is difficult for various reasons. First, any branch $\gamma$ of the IDLA tree $\TT_{\infty}$ is not produced by a single particle but by many particles, each of them adding exactly one edge depending on the shape of the current aggregate. Moreover, this random subgraph of the lattice $\Z^2$ is radial since its branches are directed to the origin and then it does not satisfy any useful invariance properties in distribution.  

An intermediate step would be to control the fluctuations of branches in the directed IDLA forest $\FF_\infty$ which presents the advantage to be invariant in distribution (and even mixing) w.r.t.\,vertical translations. Let $n\geq 1$ and let $(z_i)_{0\leq i\leq\kappa_n}$ be the branch joining a source $z_0$ on $I(\infty)$ to $z_{\kappa_n}=(n,0)$. As an illustration,  Figure \ref{fig:finiteforest} depicts the branch associated with $(55,0)$. Denote by $\Delta_n$  the maximal fluctuation (or maximal deviation) of the branch $(z_i)_{0\leq i\leq\kappa_n}$ w.r.t.\,the $x$-axis, {\it i.e.} 
\[
\Delta_n = \max_{0\leq i\leq\kappa_n} z_i(2). 
\]
{Is it true} that, with probability tending to $1$ with $n$, the maximal deviation $\Delta_n$ is negligible w.r.t.\,$n$? {From a macroscopic point of view, dealing with this question allows us to know if} the branch $(z_i)_{0\leq i\leq\kappa_n}$ asymptotically merges with the horizontal axis.

Similar open questions are also considered in  \cite[Section 5]{M17}  for a directed and external DLA model.

\section*{Acknowledgements}
This work was partially supported by the French ANR grant ASPAG (ANR-17-CE40-0017), by the French RT GeoSto (RT-3477), and by the French PEPS-JCJC 2019. {{We thank our PhD student, Keenan Penner, for pointing us a mistake in our paper and}  two anonymous referees for suggestions (in particular for pointing us the references \cite{MPZ19,PYZ20}) and improvements of the manuscript. }




\begin{thebibliography}{4}

{
\bibitem{AP17}
\textsc{T.~Antunovi\'c and E.~B. Procaccia} (2017).
\textit{Stationary Eden model on Cayley graphs.}
Ann. Appl. Probab., \textbf{27(1)}:517–-549.
}

\bibitem{AG13long}
\textsc{Asselah, A. and Gaudilli\`ere, A.} (2013). \textit{From logarithmic to subdiffusive polynomial fluctuations for internal DLA and related growth models}.  Ann. Probab. \textbf{41(3A)}:1115--1159.


\bibitem{AG13short}
\textsc{A.~{Asselah} and A.~{Gaudilli\`ere}} (2013). \textit{Sublogarithmic fluctuations for internal DLA.} Ann. Probab., \textbf{41(3A)}:1160--1179.

\bibitem{AG14}
\textsc{A.~{Asselah} and A.~{Gaudilli\`ere}} (2014).
\textit{Lower bounds on fluctuations for internal DLA.}
{Probab. Theory Related Fields}, \textbf{158(1-2)}:39--53.


\bibitem{AR16}
\textsc{A.~Asselah and H.~Rahmani} (2016).
\textit{Fluctuations for internal {DLA} on the comb.}
Ann. Inst. Henri Poincar\'{e} Probab. Stat., \textbf{52(1)}:58--83.

\bibitem{BB07}
\textsc{F.~Baccelli and C.~Bordenave} (2007).
\textit{The radial spanning tree of a {P}oisson point process.}
Ann. Appl. Probab., \textbf{17(1)}:305--359.


\bibitem{BDCKL}
\textsc{I.~Benjamini, H.~Duminil-Copin, G.~Kozma, and C.~Lucas} (2020).
\textit{Internal diffusion-limited aggregation with uniform starting points.}
Ann. Inst. Henri Poincar\'{e} Probab. Stat., \textbf{56(1)}:391--404.

\bibitem{BLPS}
\textsc{I.~Benjamini, R.~Lyons, Y.~Peres, and O.~Schramm} (1999).
\textit{Critical percolation on any nonamenable group has no infinite
  clusters.}
Ann. Probab., \textbf{27(3)}:1347--1356.

\bibitem{BKP19}
\textsc{N.~Berger, J.~J. Kagan, and E.~B. Procaccia} (2014).
\textit{Stretched {IDLA}.}
ALEA Lat. Am. J. Probab. Math. Stat., \textbf{11(1)}:471--481.

\bibitem{B04}
\textsc{S.~Blach\`ere} (2004).
\textit{Internal diffusion limited aggregation on discrete groups of
  polynomial growth.}
Random walks and geometry, pages 377--391. Walter de
  Gruyter, Berlin.
  
\bibitem{BlB07}
\textsc{S.~Blach\`ere, and S.~Brofferio} (2007).
\textit{Internal diffusion limited aggregation on discrete groups having exponential growth.}
Probab. Theory Related Fields, \textbf{137(3-4)}:323--343.

\bibitem{BK89}
\textsc{R. M.~Burton and M. S.~Keane} (1989).
\textit{Density and uniqueness in percolation.}
Comm. Math. Phys., \textbf{121}:501--505.

\bibitem{C18}
\textsc{D.~Coupier} (2018).
\textit{Sublinearity of the number of semi-infinite branches for geometric
  random trees.}
Electron. J. Probab., \textbf{23}:{Paper No.} 37, 33pp.

\bibitem{CSST}
{\textsc{D.~Coupier, K.~Saha, A.~Sarkar, and C.~Tran} (2021).
\textit{The 2d-directed spanning forest converges to the brownian web.}
Ann. Probab. \textbf{49(1):} 435-484. }

\bibitem{CT}
\textsc{D.~Coupier and V.~C. Tran} (2013).
\textit{The 2{D}-directed spanning forest is almost surely a tree.}
Random Structures Algorithms, \textbf{42(1)}:59--72.

\bibitem{DF91}
\textsc{P.~Diaconis and W.~Fulton} (1993).
\textit{A growth model, a game, an algebra, {L}agrange inversion, and
  characteristic classes.}
Rend. Sem. Mat. Univ. Politec. Torino, \textbf{49(1)}:95--119.
Commutative algebra and algebraic geometry, II (Italian) (Turin,
  1990).

\bibitem{DLYY}
\textsc{H.~{Duminil-Copin}, C.~{Lucas}, A.~{Yadin}, and A.~{Yehudayoff}} (2013).
\textit{Containing internal diffusion limited aggregation.}
Electron. Commun. Probab., \textbf{18}: {Paper No} 50, 8pp.

\bibitem{FLT}
\textsc{P.A. Ferrari, C.~Landim, and H.~Thorisson} (2004).
\textit{Poisson trees, succession lines and coalescing random walks.}
Ann. Inst. Henri Poincar\'{e} Probab. Stat., \textbf{40(2)}:141--152.

\bibitem{HN01}
\textsc{C.~Douglas Howard and Charles~M. Newman} (2001).
\textit{Geodesics and spanning trees for {E}uclidean first-passage
  percolation.}
Ann. Probab., \textbf{29(2)}:577--623.

\bibitem{H08}
\textsc{W.~Huss} (2008).
\textit{Internal diffusion-limited aggregation on non-amenable graphs.}
Electron. Commun. Probab., \textbf{13}:272--279.

\bibitem{HS12}
\textsc{W.~Huss and E.~Sava} (2012).
\textit{Internal aggregation models on comb lattices}.
Electron. J. Probab., \textbf{17}: {Paper No} 30, 21pp.

\bibitem{JLS12}
\textsc{D.~Jerison, L.~Levine, and S.~Sheffield} (2012).
\textit{Logarithmic fluctuations for internal {DLA}.}
J. Amer. Math. Soc., \textbf{25(1)}:271--301.

\bibitem{JLS13}
\textsc{D.~Jerison, L.~Levine, and S.~Sheffield} (2013).
\textit{Internal {DLA} in higher dimensions.}
Electron. J. Probab., \text{18}: {Paper No 98}, 14pp.

\bibitem{JLS14}
\textsc{D.~Jerison, L.~Levine, and S.~Sheffield} (2014).
\textit{Internal {DLA} and the {G}aussian free field.}
Duke Math. J., \textbf{163(2)}:267--308.

\bibitem{JLS14b}
\textsc{D.~Jerison, L.~Levine, and S.~Sheffield} (2014).
\textit{Internal {DLA} for cylinders.}
Advances in analysis: the legacy of {E}lias {M}. {S}tein, Princeton Math. Ser., Vol. 50 189--214. Princeton Univ. Press, Princeton, NJ.


\bibitem{Lawler95}
\textsc{G.~F. {Lawler}} (1995).
\textit{Subdiffusive fluctuations for internal diffusion limited
  aggregation.}
Ann. Probab., \textbf{23(1)}:71--86.

\bibitem{LBG92}
\textsc{G.~F. Lawler, M.~Bramson, and D.~Griffeath} (1992)
\textit{Internal {D}iffusion {L}imited {A}ggregation.}
Ann. Probab., \textbf{20(4)}:2117--2140.

\bibitem{LL10}
\textsc{G.~F. {Lawler} and V.~{Limic}} (2010).
\textit{Random walk: A modern introduction.}, Vol. 123.
\newblock Cambridge: Cambridge University Press.

\bibitem{LP10}
\textsc{L.~Levine and Y.~Peres} (2010).
\textit{Scaling limits for internal aggregation models with multiple sources.}
J. Anal. Math., \textbf{111}:151--219.

\bibitem{LS}
\textsc{L.~Levine and V.~Silvestri} (2019).
\textit{How long does it take for internal {DLA} to forget its initial
  profile?}
Probab. Theory Related Fields, \textbf{174(3-4)}:1219--1271.

\bibitem{L14}
\textsc{C.~Lucas} (2014).
\textit{The limiting shape for drifted internal diffusion limited aggregation
  is a true heat ball.}
Probab. Theory Related Fields, \textbf{159(1-2)}:197--235.


\bibitem{M17}
\textsc{S.~{Martineau}} (2017).
\textit{{Directed diffusion-limited aggregation}.}
{ALEA, Lat. Am. J. Probab. Math. Stat.}, \textbf{14(1)}:249--270.

\bibitem{MD86}
\textsc{P.~Meakin and J.M. Deutch} (1986).
\textit{The formation of surfaces by diffusion-limited annihilation.}
J. Chem. Phys., \textbf{85(4)}.


{\bibitem{MPZ19} 
\textsc{Y.~Mu, E.~B. Procaccia, and Y.~Zhang} (2019).
 \textit{Scaling limit of DLA on a long line segment.}
 Available at \url{https://arxiv.org/pdf/1912.02370.pdf}.}

{\bibitem{PYZ20}
\textsc{E.~B. Procaccia, J.~Ye, and Y.~Zhang} (2020).
\textit{Stationary DLA is well defined.}
J. Stat. Phys., \textbf{181(4)}:1089--1111.}


{\bibitem{PYZ21}
\textsc{E.~B. Procaccia, J.~Ye, and Y.~Zhang} (2021).
\textit{Stationary Harmonic Measure as the Scaling Limit of Truncated Harmonic Measure.}
ALEA, Lat. Am. J. Probab. Math. Stat., \textbf{18}:1529--1560.}


{\bibitem{PZ19}
\textsc{E.~B. Procaccia, J.~Ye, and Y.~Zhang} (2019).
\textit{Stationary Harmonic Measure and DLA in the Upper Half Plane.}
J. Stat. Phys., \textbf{176}:946--980.}


\bibitem{SS}
\textsc{K.~Saha and A.~Sarkar} (2016).
\textit{Random directed forest and the brownian web.}
Ann. Inst. Henri Poincar\'{e} Probab. Stat., \textbf{52(3)}:1106--1143.

\bibitem{SW}
\textsc{R.~Schneider and W.~Weil} (2008)
\textit{Stochastic and integral geometry}.
Probability and its Applications (New York). Springer-Verlag, Berlin.

\bibitem{Shellef}
\textsc{E.~{Shellef}} (2010).
\textit{{IDLA} on the supercritical percolation cluster.}
Electron. J. Probab., \textbf{15:} {Paper No} 24, 723--740.

\bibitem{S19}
{
\textsc{V.~Silvestri} (2020).
\textit{Internal {DLA} on cylinder graphs: fluctuations and mixing.}
Electron. Commun. Probab. \textbf{25:} 1-14.}

\end{thebibliography}
\end{document}